%% file: main.tex
\documentclass[11pt]{article}
\usepackage[T1]{fontenc}
\usepackage{newtxtext}
\usepackage[letterpaper,margin=1in,bottom=1.06in,includefoot]{geometry}

\usepackage[dvipsnames]{xcolor}
\usepackage{tikz}
\usepackage{pgfplots}
\usepackage{placeins}
\usepgfplotslibrary{fillbetween}
\pgfplotsset{compat=1.18}
\usepackage{amsmath,amsthm,mathtools}
\usepackage{newtxmath}
\usepackage{bm}
\usepackage{enumitem}
\usepackage[normalem]{ulem}

\usepackage{natbib}
\bibpunct[, ]{(}{)}{,}{a}{}{,}%

\usepackage{xr}
\makeatletter
\DeclareRobustCommand{\gref}[1]{%
  \@ifundefined{r@G-#1}{\textbf{??}}{\expandafter\expandafter\expandafter\gref@extract\csname r@G-#1\endcsname}%
}
\newcommand{\gref@extract}[5]{#1}
\makeatother

\usepackage[linktocpage=true,pagebackref=true,colorlinks=true,allcolors=blue,bookmarks=true,bookmarksopen,bookmarksnumbered]{hyperref}
\newcommand{\mathbbm}[1]{\mathbf{#1}}
\newcommand{\mcc}{\mathrm{mcc}}

\usepackage{setspace}

\input{TEX-Macros}

\newcommand{\support}[1]{\text{supp}(#1)}

\newcommand\cvar[0]{\text{CVaR}}

\newcommand\opt[0]{\text{OPT}}

\newtheorem{theorem}{Theorem}
\newtheorem{lemma}{Lemma}
\newtheorem{proposition}{Proposition}

\newtheorem{claim}{Claim}

\newtheorem{definition}{Definition}

\makeatletter
\renewenvironment{proof}[1]{%
  \par
  \pushQED{\qed}%
  \normalfont \topsep6\p@\@plus6\p@\relax
  \trivlist
  \item[\hskip\labelsep\itshape #1\@addpunct{.}]\ignorespaces
}{%
  \popQED\endtrivlist\@endpefalse
}
\makeatother
\newcommand{\Halmos}{}

\newcommand{\FIGURE}[3]{%
  \begin{figure}[ht]
    \centering
    #1%
    \caption{#2}%
    #3%
  \end{figure}%
}

\definecolor{denimlight}{RGB}{155,175,240}
\definecolor{denimblue}{RGB}{140,165,195}
\definecolor{deepblue}{RGB}{90,115,150}
\definecolor{softpeach}{RGB}{240,190,150}
\definecolor{burntpeach}{RGB}{220,140,100}
\definecolor{coolgray}{RGB}{110,120,130}
\definecolor{dustygreen}{RGB}{140,165,130}
\definecolor{mauve}{RGB}{170,140,170}

\title{Beyond Expected Profit: How Risk Aversion Alters Learning in the Data-Driven Newsvendor}
\author{Omar Besbes and Alon Rieger}
\date{}

\begin{document}

\onehalfspacing

\maketitle

\begin{abstract}
We study how risk aversion changes learning and the value of data in the classical newsvendor problem. We evaluate data-driven policies through the Conditional Value-at-Risk (CVaR) of realized loss on the combined draw of an i.i.d.\ historical sample and an independent future demand from the same distribution. We measure worst-case regret over demand distributions on a bounded interval relative to an oracle that knows the distribution. We show that, for every fixed positive risk level, Sample Average Approximation (SAA) can retain regret bounded away from zero as the sample size \(n\) grows, even at a fixed demand distribution, and characterize its asymptotic worst-case regret. We trace this failure to persistent variability in the learned order, which contributes to tail risk under this combined-risk evaluation. We establish that risk aversion also slows the minimax learning rate over all data-driven policies from \(\Theta(n^{-1/2})\) under expectation to \(\Theta(n^{-1/4})\). We construct policies that average nearby order statistics and attain this optimal rate. For finite samples, we reduce worst-case analysis for mixtures of convex combinations of order statistics to Bernoulli demand on the endpoints of the interval. We use this reduction to formulate a linear program for computing minimax policies within this class to any accuracy. In numerical experiments, we find that minimizing this objective alone can leave substantial average CVaR regret when each learned order is evaluated over future demand alone. We therefore develop a two-stage linear program that minimizes the worst-case value of this average regret while allowing a prescribed increase in the original objective. This yields improvements across a broad range of scenarios. Our results show that risk aversion changes both the rate at which data improve decisions and the policies needed to realize that improvement.
\end{abstract}

\begin{samepage}
\par
\end{samepage}

\input{TEX-Intro}

\input{TEX-ProblemAbstract}

\input{TEX-Computation}

\input{TEX-SAA}

\bibliographystyle{informs2014}
\bibliography{sample}

\appendix

\input{TEX-Appendix}

\end{document}

%% file: TEX-Macros.tex
\newcommand{\qstar}{q^*}

\newcommand{\mcost}{m}
\newcommand{\RegretCVaR}[2]{\mathcal{R}_{\delta,n}\left(#2, #1\right)}
\newcommand{\RegretCVaRDelta}[3]{\mathcal{R}_{#3,n}\left(#2, #1\right)}

\newcommand{\supp}{\text{supp}}

%% file: TEX-Intro.tex
\section{Introduction}\label{sec:intro}

The newsvendor problem is a canonical model in operations management and revenue management: a decision maker selects an inventory or capacity level before demand is realized, balancing the costs of ordering too little against those of ordering too much. In its classical form, the objective is to maximize expected profit (equivalently, to minimize expected cost). In many applications, however, managers care not only about average performance but also about exposure to rare, adverse outcomes. 
Weak demand can leave a firm with unsold inventory and large write-offs, even when its inventory decisions generate attractive average profits. These considerations motivate performance criteria that emphasize downside exposure rather than average outcomes.

Conditional Value-at-Risk (CVaR) is a widely used measure of downside risk that focuses on tail losses. 
Classical CVaR formulations assume that the decision maker knows the demand distribution. In practice, decision makers typically do not know this distribution and must learn inventory or capacity decisions from historical observations. Even when data are plentiful in aggregate, the number of relevant observations may be limited because products, markets, and demand conditions change over time.

When a decision maker learns an inventory decision from data, different historical samples can lead the same rule to select different orders. The resulting profit therefore depends on both the sample used to choose the order and the future demand. We evaluate a decision rule before observing the sample, applying CVaR to the realized loss while accounting for both sources of randomness simultaneously. We call this criterion combined-risk CVaR, following the total-CVaR perspective of~\cite{LuShanthikumarShen2015}. It captures the tail risk of repeatedly learning and deploying the same rule, as in the case of comparable products, locations, or seasons. Under this criterion, variability in the learned order can itself contribute to tail risk. How a policy uses data is therefore part of the risk-management problem.

Under an expectation objective, Sample Average Approximation (SAA), which selects the order that is optimal for the empirical demand distribution, enjoys strong learning guarantees: relatively few observations can already yield effective decisions, and worst-case regret vanishes as the sample size grows~\cite{BesbesM23}. We ask three questions under combined-risk CVaR. Does SAA retain its favorable learning guarantees? How does risk aversion change the best achievable rate of learning over all data-driven policies? Which policies attain that rate and provide effective decisions from finite data? To address these questions, we compare the tail-risk performance of a data-driven policy with that of an oracle that knows the demand distribution, and study the worst-case gap over all demand distributions on a bounded interval.

\subsection{Main contributions}

\emph{Failure of sample average approximation.}
We show that optimizing against the empirical demand distribution need not control the risk arising from variability in the learned order.
Under the expectation benchmark, SAA has vanishing distribution-free worst-case regret. Under combined-risk CVaR, for every fixed positive risk level, we exhibit a fixed two-point demand distribution under which SAA's regret converges to a strictly positive constant (Theorem~\ref{thm:saa_nonvanishing}). Its worst-case regret therefore cannot vanish as the sample size grows. We also characterize its asymptotic worst-case regret exactly (Theorem~\ref{thm:saa_exact_worst_case}).

\emph{Optimal learning rates and window policies.}
We establish that, for every fixed positive risk level, combined-risk CVaR slows the minimax learning rate from $\Theta(n^{-1/2})$ under expectation to $\Theta(n^{-1/4})$ over all data-driven policies (Theorem~\ref{thm:global_minimax_rate}). Risk aversion therefore changes the best achievable rate of learning, beyond the failure of SAA. To attain this rate, we construct explicit window policies that average order statistics around the empirical rank selected by SAA. Averaging across nearby ranks reduces the sensitivity of the learned order to sampling fluctuations. We choose uniform windows of order $n^{3/4}$ to balance this stabilization against the bias from averaging across ranks, and prove that their worst-case regret is $\Theta(n^{-1/4})$. We prove a matching lower bound over all data-driven policies, including randomized policies that need not be based on order statistics. Thus, these simple windows attain the optimal worst-case rate. Stabilizing the learned order makes it possible to realize the value of additional data at the best attainable rate.

\emph{The effect of decision variability.} To explain these findings, we use the regret decomposition of~\cite{BesbesRiegerGeneral}. They decompose combined-risk regret into an expected post-sample regret and a sample-risk premium, and derive the asymptotic order of magnitude of each term across problem classes. The first averages, over samples, the regret of the selected order when risk is evaluated over future demand alone. The second measures the additional risk from evaluating the sample and future demand simultaneously. We specialize that decomposition to the newsvendor profit. 
Under regularity and concentration conditions, the sample-risk premium vanishes more slowly than the expected post-sample regret.
We derive explicit expressions for these leading terms. Variability across samples can therefore dominate regret even when the learned order converges to the oracle.
At the fixed distribution demonstrating SAA's failure, the contrast is even starker: SAA continues to switch between the two extreme orders, with probabilities approaching one half. The post-sample regret is zero, while the sample-risk premium converges to a strictly positive constant.

{\emph{Finite-sample minimax policies within MCC.}}
{To characterize the value of finite data, we study mixtures of convex combinations of order statistics (MCC policies). This class contains SAA and the window policies, and allows us to optimize the weights on the order statistics for combined-risk CVaR.} The original minimax problem is infinite-dimensional, since the adversary ranges over all bounded demand distributions. We prove that, for MCC policies, {the worst-case value over $\mathcal F_T$ equals the supremum over Bernoulli demand on $\{0,T\}$ (Theorem~\ref{thrm:worst_case_bernoulli})}. {We show that minimax optimization over MCC can be restricted to deterministic convex combinations of order statistics, for which the Bernoulli reduction yields a linear programming formulation. We bound the discretization error explicitly, so we can compute policies whose worst-case regret is within any prescribed tolerance of the minimax value. The optimal weights do not depend on the support bound, so the resulting policy uses only the observed demands.} {Using this reduction, we quantify SAA's finite-sample worst-case regret and show that it is highly sensitive to the risk level and can remain very large at practically relevant sample sizes.} 
{With one additional observation, we apply a policy designed for the original sample size to each subsample obtained by leaving out one observation, and average the resulting orders. This gives a policy in the same class whose combined-risk regret is no greater than that of the original policy at any demand distribution.}

We develop a worst-case reduction that remains valid when regret is not separable in the demand CDF. In the expectation setting, worst-case regret reductions exploit pointwise separability in the demand CDF: the adversary can optimize $F(y)$ independently for each $y$. Under CVaR, this approach breaks down because regret depends on the ordering of losses generated by the combined randomness and is no longer pointwise separable across CDF values. 
We first approximate bounded demand distributions by discrete ones. We then exploit a convexity property of the regret induced by convex combinations of order statistics to remove interior support points successively. This argument reduces the worst-case search to Bernoulli demand on the endpoints of the bounding interval.

\emph{Refining finite-sample policy design.}
In our numerical experiments, we find that policies selected to minimize worst-case combined-risk regret on a finite grid can perform poorly on typical demand distributions. To guide policy design, we examine how these policies balance the two components of regret. We find that these policies often reduce the sample-risk premium at the expense of a larger expected post-sample regret. We therefore ask whether allowing a small increase in worst-case combined-risk regret can substantially reduce worst-case expected post-sample regret. 
We propose a two-stage linear programming procedure that minimizes worst-case expected post-sample regret subject to keeping worst-case combined-risk regret within a prescribed tolerance of the minimax value in our policy class. 
In our experiments, at the grid distribution that maximizes combined-risk regret, expected post-sample regret is substantially smaller, with little change in combined-risk regret. 
The gains over SAA are substantial across a broad range of demand distributions and sample sizes in our experiments, although the policy does not improve on SAA in every case. The two-stage program provides a computational prescription at a given sample size; window policies attain the optimal learning rate.

\subsection{Related literature}
We organize the discussion around five threads: (i) data-driven newsvendor learning under expectation objectives; (ii) risk-averse inventory with known demand distributions; (iii) total CVaR and estimation inside a risk evaluation; (iv) data-driven robustness and empirical tail objectives; and (v) quantile estimation and policies based on order statistics.

The first thread studies how historical demand observations should be mapped to inventory decisions when performance is evaluated through expected cost. The general SAA framework of~\cite{KleywegtSW02} and the sampling-based inventory policies of~\cite{LeviRS07} provide  guarantees for this approach. For the single-period newsvendor,~\cite{LeviPerakisUichanco2015} establish distribution-free bounds for SAA together with refinements that depend on properties of the demand distribution. \cite{BesbesM23} give an exact finite-sample characterization of SAA's worst-case (relative) expected regret, reduce the adversarial distribution search to Bernoulli distributions, and derive a globally minimax policy over the unrestricted class of data-to-decision mappings. 
Their analysis exploits the separability of expected regret. Under combined-risk CVaR this separability fails, so our two-point reduction requires a different argument. We characterize finite-sample minimax performance within MCC and establish the optimal asymptotic rate over all data-driven policies (Theorem~\ref{thm:global_minimax_rate}). In sequential versions of the data-driven newsvendor,~\cite{LinHKU22} analyze cumulative SAA regret under global and local separation conditions, while~\cite{LyuYYZ24} establish matching bounds under related strong-convexity conditions for more general inventory costs. \cite{ChenMa26} provide a unified analysis across regret criteria and distribution classes. Contextual extensions include~\cite{BanR19}, who derive finite-sample guarantees for feature-based empirical and kernel methods, and~\cite{BesbesMaMouchtaki2023}, who exactly characterize the context-dependent worst-case expected regret of weighted empirical policies and reduce the adversarial distribution problem to a line search. A distinct heterogeneity model is studied by~\cite{BesbesMaMouchtaki2025}, where historical observations may arise from distributions near the future demand distribution. Across this literature, sampling variability and future-demand uncertainty are averaged through expectation; combined-risk CVaR instead places both sources of randomness inside the tail-risk evaluation.

The second thread studies risk-averse inventory when the demand distribution is known. Expected-utility models of the risk-averse newsvendor date at least to~\cite{EeckhoudtGollierSchlesinger1995}. Coherent risk measures~(\cite{ArtznerDEH99}) and tractable representations of CVaR~(\cite{RockafellarU02}) provide the foundations for this literature. In the single-period newsvendor,~\cite{GotohT07} derive optimal ordering rules under CVaR, while~\cite{ChoiRuszczynski2008} study the broader class of law-invariant coherent risk measures. \cite{ChenXuZhang2009} incorporate CVaR into joint pricing and ordering decisions, and~\cite{ArikanFichtinger2017} classify the structural effects of spectral risk measures across alternative newsvendor formulations. Extensions include coherent-risk inventory models over multiple periods~(\cite{AhmedCS07}), utility-based formulations for multiperiod inventory and pricing~(\cite{ChenSSS07}), and multiproduct newsvendor models under law-invariant coherent risk measures~(\cite{ChoiRuszczynskiZhao2011}). In the papers just cited, the decision is chosen using the known demand distribution; the randomness inside the risk measure therefore stems exclusively from operational outcomes rather than from learning the decision.

The third thread studies total $\cvar$ and estimation inside a risk evaluation. 
\cite{LuShanthikumarShen2015} introduce total $\cvar$ to evaluate sampling and future-demand uncertainty in combination, and optimize positively homogeneous operational statistics for a repeated newsvendor in a location-scale family. \cite{BesbesRiegerGeneral} provide a general regret decomposition and characterize fixed-distribution asymptotics under regularity conditions. Building on these perspectives, we study worst-case regret over all demand distributions on a bounded interval. We use our two-point reduction to evaluate finite-sample performance and compute minimax policies within the MCC class up to any prescribed accuracy. We also construct window policies that attain the optimal worst-case learning rate over all data-driven policies.

Other papers bring estimation in by a different route. Some adjust a plug-in capital or quantile, or build an uncertainty set from a risk measure, so that a fixed decision meets a stated tail standard~\cite{BertsimasBrown2009,Pflug2006,GerrardTsanakas2011,FrancioniHerzog2012}. \cite{BignozziTsanakas2016} and \cite{PiteraSchmidt2018} place a random capital estimator inside the risk measure. Their objective is calibrated capital. The objective here is the regret of a learned inventory decision, uniformly over bounded demand distributions, and the associated minimax rate. Some put an outer risk on a posterior, a predictive distribution, or a second modeled source of uncertainty~\cite{TaoDelageXu2025,SadatsafaviLeeGustafson2022,JaiswalHonnappaRao2023,KimKang2021,BodnarLindholmNiklassonThorsen2022,BodnarLindholmNiklassonThorsen2020,MaChenXu2026,LinRenZhou2022,QianWangWen2019}. 
Other work studies empirical tail-risk optimization and its post-sample performance~\cite{JiangTYZ26,LiuZhu2024,Chen2024,LinChenLiShen2022,BanholzerFW22,GuiguesKS18,DentchevaPenevRuszczynski2017,Shapiro13,MengSunGoh2010}.
In this paper the demand distribution is fixed and unknown, and the training sample enters the tail through the order it produces. 
{Consistency for the post-sample risk objective does not by itself imply vanishing combined-risk regret.}

The fourth thread addresses limited data by robustifying the optimization problem itself. General data-driven robust optimization constructs ambiguity sets from statistical tests~(\cite{BertsimasCST18}) or Wasserstein neighborhoods of the empirical distribution~(\cite{EsfahaniK18}) and optimizes against the worst distribution in the resulting set. In the newsvendor setting,~\cite{SaghafianTomlin2016} combine demand observations with moment and tail information,~\cite{XuZhengJiang2022} construct ambiguity sets from nonparametric sample information, and~\cite{LeeKimMoon2021} derive explicit solutions for Wasserstein ambiguity sets, including a risk-averse CVaR formulation. A closely related empirical-tail approach is~\cite{BertsimasThiele05}, who optimize a one-sided trimmed-profit objective directly over historical observations. Data-driven distributionally robust CVaR objectives have also been studied in portfolio optimization~(\cite{KangLLZ19,PunWangYan2023}). 
These approaches use data to construct a sample-dependent ambiguity set and optimize against its worst-case distribution. Combined-risk CVaR instead evaluates the ex-ante tail risk of the random sample-dependent policy itself, without localizing the unknown demand distribution into a sample-specific ambiguity set.

The fifth thread connects our policy class to quantile estimation. The demand--supply mismatch component of the classical risk-neutral newsvendor objective is a scaled check loss underlying quantile regression~(\cite{KoenkerBassett1978}), so empirical expected-cost minimization selects a sample quantile, or equivalently an order statistic. 
In our profit-CVaR formulation, the empirical optimizer is also a sample quantile by the CVaR newsvendor characterization of~\cite{GotohT07}. More generally, deterministic convex combinations of order statistics are classical $L$-estimators~(\cite{Stigler1974}), while MCC policies also permit ex-ante randomization over such estimators. Weighted and smoothed quantile estimators in this class include the distribution-free estimator of~\cite{HarrellDavis1982}, the subsample-averaged estimator of~\cite{KaighLachenbruch1982}, and smooth or kernel estimators that spread weight across nearby order statistics~(\cite{Yang1985,SheatherMarron1990}). These estimators choose order-statistic weights to estimate a population quantile accurately. 
We choose weights to control the worst-case CVaR regret of the inventory decision.

%% file: TEX-ProblemAbstract.tex
\section{Problem Formulation}\label{sec:problem}

We consider a newsvendor problem in which a decision maker chooses an order quantity $x \in \mathbb{R}_+$ before demand $D \in \mathbb{R}_+$ is realized.
Excess capacity incurs an overage cost $h>0$ per unit, while unmet demand incurs an underage cost $b>0$ per unit. 
The newsvendor profit is given by
\[
    \Phi\left(x,D\right) ~~=~~ bD - b \cdot \left(D-x\right)^+ - h\left(x-D\right)^+ \ ,
\]
and we define the loss as negative profit,\footnote{We evaluate downside risk of realized profit, hence use its negative rather than demand-supply mismatch cost. The profit and mismatch cost formulations have the same optimizer under expectation but can differ under $\cvar$ because their difference, $-bD$, is random; see Appendix~\ref{apx:cost_profit_cvar}.}
\begin{equation}\label{eqn:neg_profit}
    c\left(x,D\right) ~~=~~ -\Phi\left(x,D\right) ~~=~~  b \cdot \left(D-x\right)^+ + h\left(x-D\right)^+ - bD  \ .
\end{equation}
Equivalently,
\[
    c\left(x,D\right) ~=~ \max\left\{-bx,\, hx-\left(h+b\right)D\right\}.
\]

\paragraph{CVaR.}
Performance is evaluated by Conditional Value-at-Risk ($\cvar$).
\begin{definition}[CVaR]
Let $\delta \in [0,1)$, and let $X$ be integrable with distribution function $F_X$. Then
\begin{equation}\label{eqn:cvar_definition}
\cvar_\delta[X] ~~=~~ \frac{1}{1-\delta} \cdot \int_{\delta}^{1} F_X^{-1}\left(p\right)\,dp \ .
\end{equation}
where $F_X^{-1}(q) := \inf\{x: F_X(x) \geq q\}$ is the lower quantile function of $X$.
\end{definition}
For $0<\delta<1$, Value-at-Risk at level $\delta$ is the $\delta$-quantile of the loss. We write it $\mathrm{VaR}_\delta[X]$,
\[
    \mathrm{VaR}_\delta[X] := F_X^{-1}(\delta).
\]
The same quantity is the smallest loss level $\ell$ such that $\mathbb P(X\leq\ell)\geq\delta$.

In words, $\cvar_\delta[X]$ is the average loss over the upper $(1-\delta)$ quantile range. 
If $0<\delta<1$ and the distribution puts no mass at $\mathrm{VaR}_\delta[X]$, this is the conditional expectation
    $\cvar_\delta[X] ~~=~~ \mathbb{E}[X \mid X \geq F_X^{-1}\left(\delta\right)]$.
By \cite[Theorem 10]{RockafellarU02}, the same quantity has the variational form
\begin{equation}\label{eqn:cvar_representation}
    \cvar_\delta[X] ~~=~~ \inf_{t \in \mathbb{R}}\left\{t + \frac{1}{1-\delta} \cdot \mathbb{E}\left[\left(X-t\right)^+\right] \right\}  \ .
\end{equation}
For $0<\delta<1$ and integrable $X$, $\mathrm{VaR}_\delta[X]$ attains the infimum, and no smaller threshold does. At $\delta=0$ the criterion is the expectation. 
For the newsvendor loss at $0<\delta<1$, we call this attaining level the cutoff, so the cutoff of an order $x$ is $\mathrm{VaR}_\delta[c(x,D)]$. 
The same representation applies when the order quantity is random. Draw $x$ from a distribution $H$, independently of demand $D \sim F$. Then
\begin{equation}\label{eqn:cvar_representation_2}
\cvar_\delta[c\left(x,D\right)] ~~=~~ \inf_{t \in \mathbb{R}}\left\{t + \frac{1}{1-\delta} \cdot \mathbb{E}_{x \sim H, D \sim F}\left[\left(c\left(x,D\right)-t\right)^+\right] \right\} \ .
\end{equation}

Demand follows an unknown distribution $F$. All that is observed is an i.i.d.\ sample $D_1^n := (D_1,\ldots,D_n)$ of size $n$. A \emph{data-driven policy} turns that sample into a distribution over order quantities. Let $\Pi$ denote this class, with the dependence on $n$ suppressed. Each $\pi \in \Pi$ is a map
\[
    \pi:D_1^n \mapsto H_{D_1^n} \ , 
\]
where $H_{D_1^n} \in \mathcal{P}(\mathcal{X})$ and $\mathcal{P}(\mathcal{X})$ is the set of distributions on the feasible orders $\mathcal{X}=\mathbb{R}_+$. When demand is supported on $[0,T]$, projecting an order onto $[0,T]$ does not raise its $\cvar$.
The order for a realized sample is drawn from $\pi(D_1^n)$. The policy is called \emph{deterministic} when that distribution is a point mass. The placed order $x^\pi$ is nevertheless random, since it depends on the sample: $D_1^n \sim F$ produces the sample, $x^\pi \sim \pi(D_1^n)$ is the resulting order, and a fresh demand $D \sim F$ is drawn independently of both, so the loss $c(x^\pi, D)$ is random.

\paragraph{Combined-randomness performance.} 
If demand is known to follow $F$, the best $\cvar$ a fixed order can achieve is the \emph{oracle} value
\[   \opt_\delta(F) ~~=~~ \inf_{x \in \mathcal{X}} \cvar_\delta[c\left(x, D\right)] \ . \]
The oracle chooses a nonrandom order, so the only randomness inside $\cvar_\delta[c(x,D)]$ is the demand draw $D\sim F$.
This value is achieved by the order quantity 
\begin{equation}\label{eqn:gotoh_takano_profit}
    x^*_{F} ~~=~~  F^{-1} \left(\qstar \right) \, \qquad \mbox{with} \quad \qstar ~~=~~ \frac{b \cdot \left(1-\delta\right)}{h+b}.
\end{equation}
We point the reader to \cite[Corollary 3.4]{GotohT07} for a derivation of this closed-form solution. We refer to $\qstar$ as the \emph{critical quantile}. The formula also holds for discrete demand: the optimality condition is $F(x-)\leq\qstar\leq F(x)$, which is satisfied by the lower quantile $F^{-1}(\qstar)$.

In contrast to the oracle, the placed order is random through the sample.
\cite{BesbesRiegerGeneral} score a data-driven policy by the $\cvar$ of its loss over the combined draw of sample, policy, and future outcome, applying the total-CVaR perspective of~\cite{LuShanthikumarShen2015} to upper-tail loss.
For this loss, that criterion is
\begin{equation}\label{eqn:policy_joint_cvar}
    \cvar_\delta[c\left(x^\pi, D\right)] ~:=~ \inf_{t \in \mathbb{R}}\left\{t + \frac{1}{1-\delta} \cdot \mathbb{E}_{D_1^n \sim F, x^\pi \sim \pi(D_1^n), D \sim F}\left[\left(c\left(x^\pi,D\right)-t\right)^+\right] \right\} \ .
\end{equation}
We require the maps below to be measurable and the expectations they induce to be finite, whether the $\cvar$ is computed at a fixed order, conditionally on the sample, or over the sample and demand together. No bound on the order set or on the loss is imposed yet.
The $\cvar$ in equation~\eqref{eqn:policy_joint_cvar} is taken before the sample is seen. A bad draw can be a history that produces a poor order, an adverse demand given that order, or a coin flip inside $\pi$; 
conditioning on both the sample and the realized order leaves only future-demand uncertainty. The decision maker chooses the rule $\pi$ up front, randomization included, rather than one order. That is the right object when the same rule is retrained across products, locations, markets, or successive seasons: each deployment records one loss, and those losses differ both because the training demands differ and because the demand that follows differs. Equation \eqref{eqn:policy_joint_cvar} is the tail risk of that full procedure. 
Post-sample $\cvar$ evaluates tail risk over future demand after the sample and any policy randomization have determined the order.

The regret of a policy $\pi$ is the difference between its combined-risk $\cvar$ and the oracle value,
\begin{equation}\label{eqn:regret_cvar_def}
    \RegretCVaR{F}{\pi} ~~=~~  \cvar_\delta[c\left(x^\pi, D\right)] - {\opt_\delta(F)}  \ .
\end{equation}
For worst-case evaluation, we fix $T>0$ and let $\mathcal F_T$ denote the class of demand distributions supported on $[0,T]$. We measure a policy's worst-case regret by $\sup_{F\in\mathcal F_T}\RegretCVaR{F}{\pi}$.

\cite{BesbesRiegerGeneral} show that the regret is nonnegative. For any threshold, the objective in the CVaR formula~\eqref{eqn:cvar_representation}, evaluated at the realized order, is at least its minimum over thresholds. That minimum is the outcome-only $\cvar$ of the order, which is at least the oracle value. Averaging over samples and then minimizing over one threshold common to every order preserves the inequality. The next section specializes their decomposition of this regret to the newsvendor.

\section{Decomposition of combined-risk regret}\label{sec:decomposition}

With $D_1^n$ observed and $x^\pi$ drawn, tail risk can be computed in future demand alone. Write $\Psi_n^\pi$ for that post-sample outcome-only $\cvar$,
\begin{equation}\label{eqn:conditional_psi}
\Psi_n^\pi
~:=~
\inf_{t \in \mathbb{R}}
\left\{
t+\frac{1}{1-\delta}\,
\mathbb{E}_{D \sim F}\!\left[
\left.\left(c\left(x^\pi,D\right)-t\right)^+\,\right|\,D_1^n,\,x^\pi
\right]
\right\} \ .
\end{equation}
Given $(D_1^n,x^\pi)$, the quantity $\Psi_n^\pi$ is a number, the outcome-only $\cvar_\delta$ of $c(x^\pi,D)$. As the sample and the induced order vary, $\Psi_n^\pi$ is random.

\paragraph{Comparison and decomposition.}
\cite[eqs.~\gref{eqn:regret_psi_decomp} and~\gref{eqn:sr_variational}]{BesbesRiegerGeneral} isolate the two components of combined-risk regret. One shared tail threshold prices the sample together with the future outcome; an action-specific threshold prices the future outcome after the sample is known. The gap between those two evaluations is the sample-risk premium, and it is nonnegative. Adding and subtracting the expected post-sample value in equation~\eqref{eqn:regret_cvar_def} gives the identity
\begin{equation}\label{eqn:regret_psi_decomp}
\RegretCVaR{F}{\pi}
=
\underbrace{\Bigl(
\mathbb{E}_{D_1^n,\,x^\pi}\!\left[\Psi_n^\pi\right]
-
\opt_\delta\left(F\right)
\Bigr)}_{\text{$\mathrm{PS}_{\delta, n}(\pi,F)$-- expected post-sample regret}}
+
\underbrace{\Bigl(
\cvar_\delta[c\left(x^\pi,D\right)]
-
\mathbb{E}_{D_1^n,\,x^\pi}\!\left[\Psi_n^\pi\right]
\Bigr)}_{\mathrm{SR}_{\delta, n}(\pi, F) \text{-- sample-risk premium}}.
\end{equation}
The first term, $\mathrm{PS}_{\delta,n}(\pi,F)$, averages, over samples and policy randomization, the outcome-only tail gap between the order just chosen and the oracle. The second term is the sample-risk premium. For $\delta\in(0,1)$, \cite[Proposition~\gref{prop:endpoints}]{BesbesRiegerGeneral} show that the premium is zero exactly when one threshold minimizes the conditional objective almost surely, and that the expected post-sample regret is zero exactly when the selected order is oracle-optimal almost surely.

\paragraph{Why the distinction matters.}
Under expected cost, the tower property makes the sample-risk premium zero. In the variational representation from equation~\eqref{eqn:policy_joint_cvar}, combined-risk $\cvar$ uses a single threshold for the distribution of losses. The post-sample evaluation in equation~\eqref{eqn:conditional_psi} allows this threshold to depend on the realized order. The premium measures the resulting increase in tail risk.

\paragraph{The leading term for the newsvendor loss.}
We specialize the results of~\cite{BesbesRiegerGeneral} to show that, under the conditions below, expected post-sample regret is quadratic in the scale of order error, while the sample-risk premium is linear. 
For an order $x$, the loss equals $-bx$ on $\{D\geq x\}$, with probability $1-F(x-)$, and is strictly larger on $\{D<x\}$. Whenever this probability is positive, the loss has an atom whose location moves with the order. The cutoff sits on that atom at the oracle, where its mass exceeds $\delta$. The statement below records that fact.
\begin{lemma}[Cutoff of the newsvendor loss]\label{lem:nv_atom}
Let $\delta\in(0,1)$. Then $c(x,D)=-bx$ for $D\geq x$ and $c(x,D)>-bx$ for $D<x$. At the oracle order $x^*_{F}=F^{-1}(\qstar)$,
\[
    \mathrm{VaR}_\delta\bigl[c(x^*_{F},D)\bigr] ~~=~~ -b\,x^*_{F},
\]
and the atom at $-bx^*_{F}$ strictly straddles level $\delta$, in the sense that $\mathbb P(c(x^*_{F},D)<-bx^*_{F})<\delta$ and $\mathbb P(c(x^*_{F},D)\leq-bx^*_{F})>\delta$.
\end{lemma}
The proof is in Appendix~\ref{prf:nv_atom}. No condition on $F$ is imposed, so the conclusion covers discrete demand. 
For the asymptotic results below, fix $\delta\in(0,1)$ and a continuous demand distribution $F$ with a unique oracle order $x^*_{F}$ in the interior of $\mathcal{X}$. Assume that $F$ is differentiable at $x^*_{F}$ with $0<F'(x^*_{F})<\infty$. The outcome-only objective then has second derivative $(h+b)F'(x^*_{F})/(1-\delta)$ at the oracle. Say that the order concentrates on $x^*_{F}$ at rate $r_n$ when
\begin{equation}\label{eqn:concentration}
    x^\pi ~~=~~ x^*_{F}+r_n Z_n\,,
    \qquad r_n\to 0\,,
    \qquad Z_n\Rightarrow Z\,,
\end{equation}
with $\bigl\{|Z_n|^{2}\bigr\}$ uniformly integrable and $Z$ nondegenerate. 
Continuity of demand ensures that the probability of a loss in a shrinking band immediately above the moving atom vanishes as the order approaches the oracle. Appendix~\ref{prf:nv_atom} verifies the corresponding isolation condition and the remaining regularity conditions.

Under these conditions, Theorem~\gref{cor:atom_category} of \cite{BesbesRiegerGeneral} gives $\mathrm{PS}_{\delta,n}(\pi,F)=\Theta(r_n^2)$ and $\mathrm{SR}_{\delta,n}(\pi,F)=\Theta(r_n)$. Thus, expected post-sample regret is quadratic in the scale of order error, while the sample-risk premium is linear and governs the rate of combined-risk regret. These are fixed-distribution conclusions, with constants that can depend on $F$.

More precisely, Proposition~\gref{prop:atom_brackets} of \cite{BesbesRiegerGeneral} gives
\begin{equation}\label{eqn:atom_approx}
    \frac{\mathrm{PS}_{\delta,n}(\pi,F)}{\mathrm{MSE}}
    ~~\underset{n \to \infty}{\longrightarrow}~~
    \frac{h+b}{2(1-\delta)}\,F'(x^*_{F})
\end{equation}
where $\mathrm{MSE}=\mathbb{E}\bigl[(x^\pi-x^*_{F})^2\bigr]$ is the mean squared error of the order about the oracle. The same proposition prices the sample-risk premium by a weighted absolute deviation of the induced cutoff. 
For this loss the cutoff equals $-bx$ near the oracle. At the oracle, the atom has mass $(h+b\delta)/(h+b)$. The atom begins at cumulative mass zero, so the upper weight on the cutoff is $\tau=\delta(h+b)/(h+b\delta)$. Passing from the cutoff to the order multiplies by $b$ and exchanges the two sides. Let $u_\tau$ be a quantile of $x^\pi$ at level $1-\tau$, and write
\begin{equation}\label{eqn:mad_tau}
    \mathrm{MAD}_\tau
    ~~=~~
    \mathbb{E}\bigl[
    (1-\tau)\,(x^\pi-u_\tau)^{+}
    +
    \tau\,(u_\tau-x^\pi)^{+}
    \bigr]
\end{equation}
for that deviation of the order. Then
\begin{equation}\label{eqn:sr_approx}
    \frac{\mathrm{SR}_{\delta,n}(\pi,F)}{\mathrm{MAD}_\tau}
    ~~\underset{n \to \infty}{\longrightarrow}~~
    \frac{b\,(h+b\delta)}{(h+b)\,(1-\delta)} \ .
\end{equation}
Dividing~\eqref{eqn:sr_approx} by~\eqref{eqn:atom_approx} produces
\begin{equation}\label{eqn:sr_over_p}
    \frac{\mathrm{SR}_{\delta,n}(\pi,F)\,/\,\mathrm{PS}_{\delta,n}(\pi,F)}{\mathrm{MAD}_\tau\,/\,\mathrm{MSE}}
    ~~\underset{n \to \infty}{\longrightarrow}~~
    \frac{2b\,(h+b\delta)}{(h+b)^{2}\,F'\!\left(x^*_{F}\right)} \ .
\end{equation}

\section{Regret of SAA} \label{sec:saa}

Let $\hat F_n$ be the empirical distribution placing mass $1/n$ on each observation. SAA plugs $\hat F_n$ into the oracle problem and orders as if that empirical distribution were the truth. With $D\sim\hat F_n$,
\begin{equation}\label{eqn:saa_objective}
    \cvar_\delta[c\left(x, D\right)] ~~=~~ \inf_t  \left\{t + \frac{1}{1-\delta} \cdot \frac{1}{n} \sum_{i=1}^n (c(x, D_i) - t)^+ \right \}
\end{equation}
Write $c_{[1]}(x)\ge\cdots\ge c_{[n]}(x)$ for the losses on the $n$ observations, ordered from largest to smallest. Optimizing out $t$ leaves, for $m=\lceil(1-\delta)n\rceil$,
\[ \cvar_\delta[c\left(x, D\right)] ~~=~~ \frac{1}{1-\delta} \cdot \left[\frac{1}{n} \cdot \sum_{i=1}^{m-1} c_{[i]}\left(x\right) + \left(1 - \delta - \frac{m-1}{n}\right) \cdot c_{[m]}\left(x\right) \right] \ . \]
The SAA order is a measurable solution of
\[
x^{\mathrm{SAA}}
\in
\arg\min_{x\in\mathcal{X}}
\cvar_{\delta,D\sim\hat F_n}[c(x,D)].
\]
We select the smallest empirical minimizer when the problem has several solutions. For this loss the oracle order is the critical quantile, so on $\hat F_n$ the rule returns the corresponding order statistic. We identify this deterministic policy with its selected order. Writing $D_{(1)}\le\cdots\le D_{(n)}$,
\[
\pi^{\mathrm{SAA}}(D_1^n)
~=~
\hat F_n^{-1}(\qstar)
~=~
D_{(k(\qstar))},
\qquad
k(\qstar):=\lceil \qstar n\rceil,
\]
under the left-continuous inverse $F^{-1}(q)=\inf\{x:F(x)\ge q\}$.

\subsection{Post-sample consistency and the sample-risk premium}

\cite[Proposition~\gref{prop:saa_post_sample_consistency}]{BesbesRiegerGeneral} show that sample average approximation is post-sample consistent when the action set is compact, the loss is continuous in the action, and the loss is uniformly bounded. The newsvendor profit is Lipschitz in the order. On a compact order set with bounded demand those conditions hold, so
\[
    \mathrm{PS}_{\delta,n}(\pi^{\mathrm{SAA}},F)\longrightarrow 0.
\]

On a compact order set, for a fixed continuous demand distribution of bounded support whose oracle order is unique, interior, and well separated, and whose demand distribution is differentiable there with a positive finite derivative, \cite[Proposition~\gref{prop:saa_concentration_rate}]{BesbesRiegerGeneral} show that SAA concentrates at rate $n^{-1/2}$. Theorem~\gref{cor:atom_category} of \cite{BesbesRiegerGeneral} then gives expected post-sample regret of order $n^{-1}$ and a sample-risk premium of order $n^{-1/2}$, so the combined regret is of order $n^{-1/2}$. Outside this regime 
convergence can be faster: 
for finite-support demand, the SAA order converges exponentially fast to the oracle when $\qstar$ lies strictly inside a CDF jump. For point-mass demand, SAA has zero regret.

\paragraph{When does the sample-risk premium vanish?}
Under expected-cost regret the sample-risk premium is zero. Under $\cvar$, \cite[Propositions~\gref{prop:endpoints} and~\gref{prop:multiple_optima}]{BesbesRiegerGeneral} 
show that, when the order set is compact and the loss is bounded and continuous in the order, the premium tends to zero when the selected order converges to a constant. It can remain positive when the limiting decision mixes optimal orders whose conditional objectives have no common minimizing threshold. 

Under these conditions, SAA converges to the oracle order when it is unique, and the premium tends to zero. A nonvanishing premium for SAA therefore requires multiple optimal orders.

\begin{theorem}\label{thm:saa_nonvanishing}
    For every fixed $\delta\in(0,1)$, the critical Bernoulli distribution $\mathrm{Ber}(1-\qstar)$ satisfies
    \[
        \RegretCVaR{\mathrm{Ber}(1-\qstar)}{\pi^{\mathrm{SAA}}} ~\to~ \begin{cases}
        \frac{h \cdot b }{2\left(h+b\right)} \quad & \delta \geq \frac{h}{2h+b} \\
        \frac{b \cdot \delta}{2 \cdot \left(1-\delta\right)} \quad & \delta < \frac{h}{2h+b}
        \end{cases} \ .
    \]
    Consequently, for every fixed $\delta\in(0,1)$, SAA is not universally consistent under combined-risk regret: there exists a fixed distribution under which its regret converges to a strictly positive constant.
\end{theorem}
The proof is provided in Appendix~\ref{prf:saa_nonvanishing}. At $\mathrm{Ber}(1-\qstar)$, every order in $[0,1]$ is oracle-optimal, but the empirical quantile selects only the endpoints: the number of zero observations fluctuates around $\qstar n$, so SAA selects $0$ and $1$ with asymptotic probability $1/2$. Ordering $0$ makes the loss identically zero, so its cutoff is $0$; ordering $1$ puts probability $1-\qstar>\delta$ on the loss $-b$, so its cutoff is $-b$. Each pure action therefore has zero outcome-only regret, the expected post-sample regret is zero for every sample size, and the sample-risk premium tends to the positive constant in the theorem. This is the mixing regime of \cite[Proposition~\gref{prop:multiple_optima}]{BesbesRiegerGeneral} for the newsvendor, a property of the cost, the distribution, and the policy together. 
A policy that selects orders in $[0,1]$ and converges to a fixed order has a premium that tends to zero.

In the risk-neutral newsvendor, SAA's worst-case regret over $\mathcal F_T$ vanishes as $n\to\infty$ for fixed $T$. The relative-regret bound of \cite[Theorem~5]{BesbesM23} implies an $O(n^{-1/2})$ absolute-regret bound on this class, since the optimal expected mismatch cost is uniformly bounded. This reflects the stability of expected-value objectives under estimation error.

Theorem~\ref{thm:saa_nonvanishing} shows that under CVaR the same conclusion need not hold. Thus, the failure is not that the empirical distribution fails to converge to $F$. Rather, at the critical Bernoulli instance, convergence of the empirical distribution is compatible with persistent randomness in the empirical $\qstar$-quantile. Under expectation the effect of this randomness is averaged out; under combined CVaR it remains visible in the tail of the unconditional loss distribution.

\begingroup\emergencystretch=3em
Post-sample near-optimality results, such as those in \cite{JiangTYZ26}, concern future-demand risk conditional on the learned policy. Theorem~\ref{thm:saa_nonvanishing} shows why such guarantees need not extend to combined-risk $\cvar$, which also includes sample-induced randomness in the selected order.
\par\endgroup

\subsection{\texorpdfstring{SAA at Fixed Demand Distributions}{SAA at Fixed Demand Distributions}}\label{sec:saa_fixed}

The preceding analysis reveals an important asymmetry in the performance of SAA. Under broad regularity conditions, its expected post-sample regret vanishes, and may already be small at moderate sample sizes. By contrast, the sample-risk premium can be substantial and, as Theorem~\ref{thm:saa_nonvanishing} demonstrates, may even remain nonvanishing asymptotically.

Figure~\ref{fig:saa_decomp_distributions} in Appendix~\ref{apx:additional_figures} illustrates this asymmetry under four demand distributions. Uniform demand on $[0,1]$ lies in $\mathcal{F}_1$. Exponential demand with rate $1$, lognormal demand with $\mu_{\log}=1$ and $\sigma_{\log}=1.8$, and gamma demand with shape $4$ and scale $0.125$ are unbounded and lie outside $\mathcal{F}_1$. Those three distributions illustrate a support-agnostic rule, and the bounded-class worst-case guarantee does not cover them. We decompose the combined-risk regret of SAA into its expected post-sample regret and sample-risk premium. Across these examples, SAA performs relatively well according to the post-sample criterion, whereas a substantial share of its combined-risk regret is generated by the
sample-risk premium.

\medskip

For all four distributions, the expected post-sample regret of SAA is $\Theta(n^{-1})$, whereas its sample-risk premium is $\Theta(n^{-1/2})$. Thus, the sample-risk premium determines the asymptotic rate of combined-risk regret, including for the distributions with unbounded support. Proposition~\ref{prop:saa_fixed_rates} (Appendix~\ref{apx:saa_fixed_rates}) establishes these rates by verifying the conditions of the atomic-regime result of~\cite{BesbesRiegerGeneral}. These comparisons are Monte Carlo estimates of combined-risk regret, and each series records its standard error. Appendix~\ref{apx:fixed_demand_mc} records the computation.

The fixed-distribution failure in Theorem~\ref{thm:saa_nonvanishing} and these numerical examples motivate two questions: how large can SAA's worst-case regret be, and which policies can improve on it? Under combined-risk evaluation, a policy should not be judged solely by its expected post-sample regret. Reducing the sample-risk premium may require accepting some deterioration in post-sample performance. The relevant policy-design problem is therefore to balance the two components of the decomposition, potentially in a different way from SAA. We start by introducing the class of data-driven policies that we will use to study this question.

\section{\texorpdfstring{Worst-Case Analysis and Policy Design}{Worst-Case Analysis and Policy Design}} \label{sec:order_stats}
How much can we improve on SAA within a class that also allows us to average order statistics? We now introduce the focal class of data-driven policies used in the rest of the paper: mixtures of convex combinations of order statistics. This class contains empirical-quantile policies such as SAA, while also allowing convex combinations of several order statistics. 
Once their weights are specified, MCC policies use only the observed demands to select an order. The scaling result below shows that the same optimal weights can be used for every support bound $T>0$.

\begin{definition}
    Consider $n$ observed demands $D_1,\ldots,D_n$, and let $D_{(1)} \leq D_{(2)} \leq \cdots \leq D_{(n)}$ denote the order statistics. A policy $\pi$ in $\Pi$ is said to be a mixture of convex combinations of order statistics if the realized decision it prescribes takes the form
    \[
    x^\pi\left(D_1^n\right) \;=\; \sum_{k=1}^n \alpha_k \cdot D_{(k)} \, ,
    \]
    where the weights $\boldsymbol{\alpha} = (\alpha_1,\ldots,\alpha_n)$ are drawn from a fixed probability distribution on the simplex $\Delta^{n-1} := \{\boldsymbol{\alpha} \in \mathbb{R}^n : \alpha_k \geq 0\;\forall k,\ \sum_{k=1}^n \alpha_k = 1\}$, with this distribution specified ex ante and the draw independent of the sample $D_1^n$.
    \end{definition}
    We refer to policies in this class as MCC policies and use $\pi^{\mcc}$ to denote an arbitrary MCC policy.
    Accordingly, we write $x^\mcc(D_1^n)$ for the action prescribed by $\pi^{\mcc}$.
    With $n$ fixed, let $\Pi^{\mcc}$ denote the subclass of $\Pi$ consisting of MCC policies.
    For any quantile level $q\in(0,1]$, write $k(q):=\lceil qn\rceil$ for the corresponding empirical rank.
    Common policy types are recovered by restricting the distribution of $\boldsymbol{\alpha}$:
    \begin{enumerate}
      \item \emph{Order statistic policies}: the distribution of $\boldsymbol{\alpha}$ puts all mass on a single fixed vertex $\mathbf e_{k(q)}$ of $\Delta^{n-1}$, yielding
      \[
      x^\pi(D_1^n)=D_{(k(q))}=\hat{F}_n^{-1}(q)
      \ . \]
      where $\hat{F}_n$ is the empirical CDF induced by the sample. SAA is recovered by setting $q=\qstar$.
      \item \emph{Mixtures of order statistic policies}: the distribution of $\boldsymbol{\alpha}$ is supported on the vertices $\{\mathbf e_1,\ldots,\mathbf e_n\}$ of $\Delta^{n-1}$, so one first randomizes over an index and then selects the corresponding order statistic.
      \item \emph{Convex combinations of order statistics}: the distribution of $\boldsymbol{\alpha}$ is degenerate at a fixed (not necessarily extreme) point of $\Delta^{n-1}$, so several order statistics are combined with deterministic nonnegative weights summing to one.
      \item \emph{Mixtures of convex combinations of order statistics}: the distribution of $\boldsymbol{\alpha}$ is arbitrary on $\Delta^{n-1}$; this is the full class. 
    \end{enumerate}
    For $n\geq 2$, relative to the counting-policy class in \cite[Definition 4]{BesbesMaMouchtaki2023}, $\Pi^{\mcc}$ is incomparable: MCC contains sample-only convex-combination rules that are not counting policies, while counting policies can prescribe constants outside the sample convex hull. Thus, our structural results apply beyond the counting-policy framework
    ; Appendix~\ref{apx:mcc_counting_incomparable} gives explicit examples. 
    
    With the characterization of our focal class, we are interested in the minimax regret problem over this class:
    \[ \inf_{\pi \in \Pi^{\mcc}} \sup_{F \in \mathcal{F}_T} \RegretCVaR{F}{\pi} \ , \]
    where $\mathcal{F}_T$ is the set of distributions with support contained in $[0,T]$.
   We first reduce worst-case evaluation of MCC policies to Bernoulli demand. We apply this reduction to SAA in Section~\ref{sec:saa_worst_case}, then show that deterministic weights suffice and formulate policy selection as a semi-infinite linear program, which we can approximate to any prescribed accuracy.

%% file: TEX-Computation.tex
\subsection{Worst-Case Analysis}\label{sec:worst_case_analysis}

In this section, we characterize structural properties of worst-case distributions for policies in $\Pi^{\mcc}$. 
Recalling that $\mathcal{F}_T$
denotes the class of demand distributions with demand in $[0,T]$, and observing that policies in $\Pi^{\mcc}$ are homogeneous in the sample, and the newsvendor loss is homogeneous in demand and order quantity, we claim that regret and expected post-sample regret scale linearly with the support bound. Therefore, we can without loss of generality restrict attention to distributions with demand in the unit interval, which are denoted by $\mathcal{F}_1$. We formalize this notion in the following claim.
\begin{claim}\label{clm:scaling}
    Let $T>0$, let $F_T \in \mathcal{F}_T$ and let $\pi^{\mcc}$ be an MCC policy. Then there exists a distribution with CDF $F_1 \in \mathcal{F}_1$ such that
    \begin{align*}
        \RegretCVaR{F_1}{\pi^{\mcc}}
        ~~=~~ \frac{ \RegretCVaR{F_T}{\pi^{\mcc}}}{T}, \quad &
        \mathrm{PS}_{\delta,n}(\pi^{\mcc}, F_1)
        ~~=~~ \frac{ \mathrm{PS}_{\delta,n}(\pi^{\mcc}, F_T)}{T}.
    \end{align*}
    Conversely, for every $F_1\in\mathcal{F}_1$, there exists $F_T\in\mathcal{F}_T$ satisfying the same equalities.
\end{claim}
The proof is presented in Appendix~\ref{prf:scaling}. By Claim~\ref{clm:scaling}, for every support bound $T>0$ and every $\pi^{\mcc}$,
\[
\sup_{F\in\mathcal{F}_T}\RegretCVaR{F}{\pi^{\mcc}}
~=~
T \cdot \sup_{F\in\mathcal{F}_1}\RegretCVaR{F}{\pi^{\mcc}} \ ,
\]
and similarly
\[ \sup_{F\in\mathcal{F}_T}\mathrm{PS}_{\delta,n}(\pi^{\mcc},F) ~=~ T \cdot \sup_{F\in\mathcal{F}_1}\mathrm{PS}_{\delta,n}(\pi^{\mcc},F) \ , \]
so the general-$T$ worst case reduces to $\mathcal{F}_1$. From this point onward, we restrict attention to distributions in $\mathcal{F}_1$.

A priori, the worst-case regret problems are infinite-dimensional, since the adversary optimizes over $\mathcal{F}_1$. Our objective is therefore to identify structure that reduces this search to a tractable subclass without loss of worst-case value.
For related risk-neutral formulations,~\cite{BesbesM23,BesbesMaMouchtaki2023} establish that the worst-case distribution is Bernoulli by exploiting a separable representation of the regret in the form
\[
\sup_{F\in\mathcal{F}_1}\int_0^\infty \Psi\left(F\left(y\right)\right)\,dy.
\]
In such settings, maximizing the integrand pointwise in $F(y)$ for each $y$ yields a global maximizer. In our $\cvar$-regret setting, however, 
the regret depends on the joint ordering
of loss outcomes and is not separable in the values of the CDF. Therefore, we cannot leverage existing approaches. We instead develop a successive support-reduction
argument to show that, for any policy in this broad class, the worst-case value
can be approached arbitrarily closely within a small subclass of $\mathcal{F}_1$. 
Our argument proceeds by support reduction rather than pointwise CDF
optimization.
We first approximate arbitrary bounded distributions by discrete ones. We then show that, for an interior support point, moving its mass to an appropriate neighboring support point cannot decrease regret. The key step is a special
convexity property: along such a one-point support perturbation, the
data-driven combined-CVaR and expected post-sample CVaR terms are convex, while the oracle benchmark is affine.
Iterating this support-merging argument collapses any discrete distribution to
two points while always (weakly) increasing  regret.

\begin{theorem}\label{thrm:worst_case_bernoulli}
    Let $F \in \mathcal{F}_1$ and let $\pi^{\mcc}$ be an MCC policy. For every $\epsilon > 0$, there exist $\mu,\mu'\in[0,1]$ and Bernoulli distributions $F_{\mathrm{Ber}}=\mathrm{Ber}(\mu)$ and $F_{\mathrm{Ber}}'=\mathrm{Ber}(\mu')$, each supported on $\{0,1\}$, such that
    \begin{align*}
        \RegretCVaR{F_{\mathrm{Ber}}}{\pi^{\mcc}}
        ~~\geq~~ \RegretCVaR{F}{\pi^{\mcc}} - \epsilon,
        \quad &
        \mathrm{PS}_{\delta,n}(\pi^{\mcc}, F_{\mathrm{Ber}}')
        ~~\geq~~ \mathrm{PS}_{\delta,n}(\pi^{\mcc}, F) - \epsilon.
    \end{align*}
\end{theorem}
\begin{proof}{Proof of Theorem \ref{thrm:worst_case_bernoulli}.}\label{prf:worst_case_bernoulli}
The proof proceeds in four steps, applied separately to combined-risk regret and to expected post-sample regret. The two criteria may select different merges, so the resulting Bernoulli distributions need not coincide. First, we discretize an arbitrary distribution in $\mathcal{F}_1$ while preserving each criterion up to an arbitrarily small error. Second, we iteratively reduce support size to two points, using convexity in the location of an interior support point (with masses and remaining support fixed). Third, we apply a shifting invariance argument to move the lower support point to zero without changing either criterion. Finally, we apply scaling to reduce each resulting two-point distribution to a Bernoulli distribution on $\{0,1\}$.

\paragraph{Discretization.}
Let $F \in \mathcal{F}_1$, and fix $\epsilon > 0$. We begin by proving that there exists a discrete distribution in $\mathcal{F}_1$, with finite support size, whose combined-risk regret and expected post-sample regret are arbitrarily close to those associated with $F$.

\begin{lemma} \label{lem:discretization}
    Let $F \in \mathcal{F}_1$. There exists a sequence of discrete distributions $\{F_{d,m}\}_{m\geq 1} \subseteq \mathcal{F}_1$ with support of size at most $m$ such that, as $m\to\infty$,
    \[
        \RegretCVaR{F_{d,m}}{\pi^{\mcc}} ~~\to~~ \RegretCVaR{F}{\pi^{\mcc}},
        \qquad
        \mathrm{PS}_{\delta,n}(\pi^{\mcc},F_{d,m}) ~~\to~~ \mathrm{PS}_{\delta,n}(\pi^{\mcc},F).
    \]
\end{lemma}
The proof of the lemma is in Appendix~\ref{prf:discretization_coupling_once}.

\paragraph{Support size reduction.}
Lemma~\ref{lem:discretization} gives a sequence $\{F_{d,m}\}_{m\geq 1}$ of discrete distributions along which both criteria converge to the values associated with $F$. In particular, for every $\epsilon>0$ there exists $m(\epsilon)$ such that if we set $F_d:=F_{d,m(\epsilon)}$, then
\[
    \RegretCVaR{F_d}{\pi^{\mcc}} ~~\geq~~ \RegretCVaR{F}{\pi^{\mcc}}-\epsilon,
    \qquad
    \mathrm{PS}_{\delta,n}(\pi^{\mcc},F_d) ~~\geq~~ \mathrm{PS}_{\delta,n}(\pi^{\mcc},F)-\epsilon .
\]
If $F_d$ is supported on a single point, then both criteria equal zero, and the theorem statement holds for a degenerate Bernoulli distribution.

Henceforth assume that $F_d$ has at least two support points. If it has exactly two, set $F_2:=F_d$ and $F_2':=F_d$, and proceed directly to the shifting step below. Otherwise, write $F_d=\sum_{i=1}^K p_i\,\delta_{a_i}$ with strictly ordered support $0\leq a_1<\cdots<a_K\leq 1$, positive masses $p_i>0$, and set $a_{-\ell}:=(a_1,\ldots,a_{\ell-1},a_{\ell+1},\ldots,a_K)$ and $p=(p_1,\ldots,p_K)$.
  Fix an interior index $\ell\in\{2,\ldots,K-1\}$ and define
\[
F_t^{(\ell)} \;:=\; \sum_{i\neq \ell} p_i\,\delta_{a_i} \;+\; p_\ell\,\delta_{t},
\qquad
t\in[a_{\ell-1},a_{\ell+1}] \ .
\]
For each criterion, we show that we can move mass from an interior support point to an adjacent support point without decreasing its value. The choice of adjacent support point may differ across criteria. To that end, we prove that both criteria are convex in the location of that point. We decompose each criterion into a data-driven term and the oracle term, both parameterized by the perturbation variable $t$.
We first prove that the data-driven terms are convex in $t$.
\begin{lemma}\label{lem:data_driven_cvar_convex}
    Fix $\ell\in\{2,\dots,K-1\}$, fix $a_{-\ell}$ and $p$, and let $F_t^{(\ell)}$ be as defined above.
    Couple the family $\{F_t^{(\ell)}\}$ using fixed atom labels whose probabilities do not depend on $t$, sending label $\ell$ to the moving support point $t$ and every other label $i$ to $a_i$, and write $D_1(t),\ldots,D_n(t),D(t)$ for the resulting demands (independent of any policy randomization). Define
    \[
        C_t \;:=\; c\bigl(x^{\mcc}(D_1^n(t)),D(t)\bigr).
    \]
    Then $t\mapsto \cvar_\delta(C_t)$ and $t\mapsto \mathbb{E}[\Psi_n^{\pi^{\mcc}}]$ are convex on $[a_{\ell-1},a_{\ell+1}]$, where $\Psi_n^{\pi^{\mcc}}$ is the post-sample outcome-only $\cvar_\delta$ in equation~\eqref{eqn:conditional_psi}.
\end{lemma}
The proof is in Appendix~\ref{prf:data_driven_cvar_convex}. At a high level, under this coupling, for almost every realization of the labels and the policy randomization the map $t\mapsto C_t$ is convex on $[a_{\ell-1},a_{\ell+1}]$; convexity and monotonicity of $\cvar_\delta(\cdot)$ then yield both conclusions, the second after conditioning on the sample labels and the realized coefficient vector.

We next prove that the oracle benchmark is affine in the same perturbation variable $t$.
\begin{lemma}\label{lem:off_linear}
    Fix $\ell\in\{2,\dots,K-1\}$, fix $a_{-\ell}$ and $p$, and let $F_t^{(\ell)}$ be as defined above. Then $t\mapsto \opt_\delta(F_t^{(\ell)})$ is affine on $[a_{\ell-1},a_{\ell+1}]$.
\end{lemma}
The proof is in Appendix~\ref{prf:off_linear}. By Lemma~\ref{lem:data_driven_cvar_convex}, both data-driven terms are convex in $t$, and by Lemma~\ref{lem:off_linear}, the oracle term is affine in $t$. Hence $t\mapsto \RegretCVaR{F_t^{(\ell)}}{\pi^{\mcc}}$ and $t\mapsto \mathrm{PS}_{\delta,n}(\pi^{\mcc},F_t^{(\ell)})$ are convex on $[a_{\ell-1},a_{\ell+1}]$. Now fix any interior support point $a_{\ell}$, write
$a_{\ell}=\lambda a_{\ell-1}+(1-\lambda)a_{\ell+1}$ for some $\lambda\in(0,1)$, and apply convexity:
\[
\RegretCVaR{F_{a_\ell}^{(\ell)}}{\pi^{\mcc}}
~~\leq~~
\lambda \RegretCVaR{F_{a_{\ell-1}}^{(\ell)}}{\pi^{\mcc}}
~+~
(1-\lambda)\RegretCVaR{F_{a_{\ell+1}}^{(\ell)}}{\pi^{\mcc}}.
\]
Consequently, the regret at $a_\ell$ is no greater than the regret at at least one endpoint, so we can merge the mass at $a_\ell$ into a neighboring support point without decreasing regret, reducing support size by one. For example, if the maximizing endpoint is $a_{\ell-1}$, replace its mass by $p_{\ell-1}+p_\ell$, delete the zero-mass point at $a_\ell$, and re-index; the case $a_{\ell+1}$ is analogous. 
The resulting distribution is now supported  on $K-1$  points in $[0,1]$; re-index the remaining support points and repeat the same argument on the updated distribution. Repeating this argument inductively yields a 2-point distribution $F_2 \in \mathcal{F}_1$ such that
\[
\RegretCVaR{F_{2}}{\pi^{\mcc}} ~~\geq~~ \RegretCVaR{F_{d}}{\pi^{\mcc}}  ~~\geq~~ \RegretCVaR{F}{\pi^{\mcc}} - \epsilon  \ .
\]
The same endpoint comparison, applied to $\mathrm{PS}_{\delta,n}(\pi^{\mcc},F_t^{(\ell)})$ rather than to combined-risk regret, produces a 2-point distribution $F_2'\in\mathcal{F}_1$ such that
\[
\mathrm{PS}_{\delta,n}(\pi^{\mcc},F_2') ~~\geq~~ \mathrm{PS}_{\delta,n}(\pi^{\mcc},F_d) ~~\geq~~ \mathrm{PS}_{\delta,n}(\pi^{\mcc},F) - \epsilon \ .
\]
The two merges may select different endpoints, so $F_2'$ need not equal $F_2$.
\paragraph{Shifting.}
Let $F_2$ have support $\{a_1,a_2\}\subseteq[0,1]$, with $a_1\leq a_2$, and let $F_2'$ be the two-point distribution obtained above for expected post-sample regret. If a lower support point is positive we shift that distribution, without changing the corresponding criterion. The next claim formalizes this invariance.
\begin{claim}\label{clm:shifting}
    Let $F_2 \in \mathcal{F}_1$ be the CDF of a distribution with support $\{a_1, a_2\}$, and $0<a_1 < a_2$. Then there exists another distribution with CDF $F_0 \in \mathcal{F}_1$, with support $\{0, a_2-a_1\}$ such that
    \[\RegretCVaR{F_{{0}}}{\pi^{\mcc}} ~~=~~ \RegretCVaR{F_2}{\pi^{\mcc}},
    \qquad
    \mathrm{PS}_{\delta,n}(\pi^{\mcc},F_0) ~~=~~ \mathrm{PS}_{\delta,n}(\pi^{\mcc},F_2) \ . \]
\end{claim}
The proof is in Appendix~\ref{prf:shifting}.
Applying Claim~\ref{clm:shifting} to $F_2$ yields $F_0 \in \mathcal{F}_1$ with support $\{0,a\}$ for some $a \in (0,1]$ and $\RegretCVaR{F_0}{\pi^{\mcc}} = \RegretCVaR{F_2}{\pi^{\mcc}}$. Applying the same claim to $F_2'$ yields $F_0'\in\mathcal{F}_1$ with support $\{0,a'\}$ for some $a'\in(0,1]$ and $\mathrm{PS}_{\delta,n}(\pi^{\mcc},F_0')=\mathrm{PS}_{\delta,n}(\pi^{\mcc},F_2')$.
View $F_0$ as an element of $\mathcal{F}_T$ with $T=a$ and apply Claim~\ref{clm:scaling} with this $T$: the rescaled distribution $F_{\mathrm{Ber}}(y):=F_0(ay)$ lies in $\mathcal{F}_1$, is Bernoulli on $\{0,1\}$, and satisfies $\RegretCVaR{F_{\mathrm{Ber}}}{\pi^{\mcc}} = \RegretCVaR{F_0}{\pi^{\mcc}}/a$.
Since $a \in (0,1]$ and regret is nonnegative,
\[
\RegretCVaR{F_{\mathrm{Ber}}}{\pi^{\mcc}} ~~\geq~~ \RegretCVaR{F_0}{\pi^{\mcc}} ~~=~~ \RegretCVaR{F_2}{\pi^{\mcc}} ~~\geq~~ \RegretCVaR{F}{\pi^{\mcc}} - \epsilon \ .
\]
The same scaling applied to $F_0'$, using the identity for $\mathrm{PS}_{\delta,n}$ in Claim~\ref{clm:scaling}, produces a Bernoulli distribution $F_{\mathrm{Ber}}'$ on $\{0,1\}$ with
\[
\mathrm{PS}_{\delta,n}(\pi^{\mcc},F_{\mathrm{Ber}}') ~~\geq~~ \mathrm{PS}_{\delta,n}(\pi^{\mcc},F_0') ~~=~~ \mathrm{PS}_{\delta,n}(\pi^{\mcc},F_2') ~~\geq~~ \mathrm{PS}_{\delta,n}(\pi^{\mcc},F) - \epsilon \ ,
\]
where the first inequality uses $a'\in(0,1]$ and nonnegativity of expected post-sample regret.
\end{proof}

Consequently, for every $\pi^{\mcc}\in\Pi^{\mcc}$,
\[
\sup_{F\in\mathcal{F}_1}\RegretCVaR{F}{\pi^{\mcc}}
~=
\sup_{\mu\in[0,1]}\RegretCVaR{\mathrm{Ber}(\mu)}{\pi^{\mcc}}.
\]

\begingroup\emergencystretch=3em
The theorem establishes that, for an MCC policy $\pi^{\mcc}$, the worst-case $\cvar$-regret over $\mathcal{F}_1$ can be evaluated by restricting attention to Bernoulli distributions. 
Our result complements existing structural reductions. MCC and counting policies are not nested in general (Appendix~\ref{apx:mcc_counting_incomparable}). We establish the Bernoulli reduction under combined-risk CVaR through convexity and successive support reduction, without relying on separability in the demand CDF.
\par\endgroup

Theorem~\ref{thrm:worst_case_bernoulli} is a fundamental building block of
the paper. It replaces an infinite-dimensional adversarial search over
$\mathcal F_1$ by a one-dimensional search over the Bernoulli success
probability $\mu$. 
We first apply this reduction to characterize SAA's asymptotic and finite-sample worst-case regret. The remainder of the paper uses the same
one-parameter representation to compute a near-optimal MCC
policy and to analyze the other policies in the focal class, identifying the
mechanisms that govern their finite-sample and asymptotic performance.

\subsection{\texorpdfstring{Worst-Case Performance of SAA}{Worst-Case Performance of SAA}}\label{sec:saa_worst_case}

Since SAA selects a single order statistic, Theorem~\ref{thrm:worst_case_bernoulli} applies at every sample size. Combining this reduction with the analysis of SAA under Bernoulli demand gives its exact asymptotic worst-case regret. Restoring the support scale $T$, we obtain the following result.

\begin{theorem}\label{thm:saa_exact_worst_case}
For fixed $h,b,T>0$ and $\delta\in(0,1)$, SAA satisfies
\[
\displaystyle
\lim_{n\to\infty}\sup_{F\in\mathcal F_T}
\RegretCVaR{F}{\pi^{\mathrm{SAA}}}
=T\frac{hb\delta}{h+b\delta}.
\]
The limit is attained along a sequence of Bernoulli distributions on $\{0,T\}$ whose probability of zero demand differs from $\qstar$ by $O(n^{-1/2})$.
\end{theorem}

The proof is in Appendix~\ref{prf:saa_exact_worst_case}. Unlike Theorem~\ref{thm:saa_nonvanishing}, the worst-case characterization permits demand to vary with $n$. Along the attaining sequence, expected post-sample regret tends to zero, so the limiting regret is entirely sample-risk premium. At $h=1$, $b=9$, $\delta=.95$, and $T=1$, the asymptotic worst case is approximately $0.8953$, compared with $0.45$ at the fixed critical Bernoulli. The larger finite-sample values below therefore reflect slow convergence to this limit.

\paragraph{Finite-sample worst-case regret.}
\begingroup\emergencystretch=5em
For each sample size, we numerically maximize SAA's combined-risk regret and expected post-sample regret separately over demand distributions supported on $[0,1]$.
By Theorem~\ref{thrm:worst_case_bernoulli}, we optimize each criterion over a Bernoulli parameter. The two criteria need not share a worst case, so the curves are not the two terms of equation~\eqref{eqn:regret_psi_decomp} at one distribution. 
The vertical gap nevertheless provides a lower bound on the sample-risk premium at a worst-case distribution. All suprema below are over $F\in\mathcal F_1$. On a distribution $F^\star$ that attains the worst-case combined-risk regret,
\[
\mathrm{SR}_{\delta,n}(\pi^{\mathrm{SAA}}, F^\star)
=
\RegretCVaR{F^\star}{\pi^{\mathrm{SAA}}}
-
\mathrm{PS}_{\delta,n}(\pi^{\mathrm{SAA}}, F^\star)
\geq
\RegretCVaR{F^\star}{\pi^{\mathrm{SAA}}}
-
\sup_{F} \mathrm{PS}_{\delta,n}(\pi^{\mathrm{SAA}}, F).
\]
The right-hand side is 
the difference between the two worst-case values illustrated in the figure. It lower-bounds the sample-risk premium on a distribution that attains the worst-case combined-risk regret, and the bound is strict when $\mathrm{PS}_{\delta,n}(\pi^{\mathrm{SAA}}, F^\star) < \sup_{F} \mathrm{PS}_{\delta,n}(\pi^{\mathrm{SAA}}, F)$.
\par\endgroup

\begingroup\emergencystretch=3em
The analytical SAA series has been recomputed using exact critical ranks and refined one-dimensional maximization; see Appendix~\ref{apx:saa_numerics}.
\par\endgroup

\medskip
\FIGURE
{%
\begin{minipage}{\linewidth}
\centering
\includegraphics[width=444.865bp,height=165.302bp]{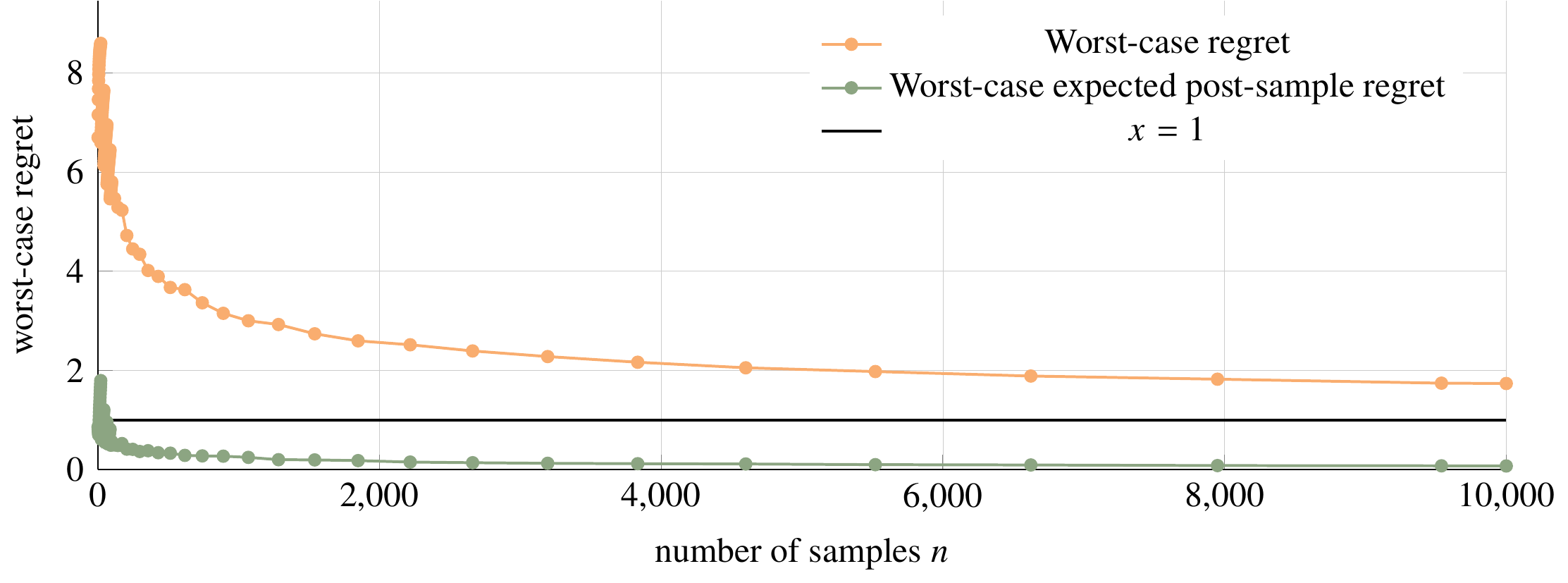}\\
{\small $h = 1$, $b=9$, $\delta=0.95$.}
\end{minipage}
}
{Numerical approximations to SAA worst-case combined-risk and expected post-sample regret.\label{fig:saa_performance_linear}}
{}

\medskip
Two comparisons stand out. 
For $h=1$, $b=9$, and $\delta=0.95$, SAA's worst-case combined-risk regret remains above one throughout the plotted range, up to $n=10{,}000$. The policy that ignores the sample and orders $x=1$ has worst-case regret one over $\mathcal F_1$, so it has lower worst-case regret over this range. The worst-case expected post-sample regret is much smaller and continues to fall. On a distribution that attains the worst-case combined-risk regret, the sample-risk premium is at least the vertical gap between these curves, and that gap accounts for nearly all of the regret at large sample sizes. The poor worst-case performance of SAA is therefore carried by the sample-risk premium, while its expected post-sample regret remains small even in the worst case.
This benchmark uses the known support bound; SAA uses only the observed demands. The comparison highlights the potential cost of using empirical optimization under high risk aversion.

The next figure shows the worst-case performance of SAA for different values of the risk aversion level $\delta$. The figure shows that the performance of SAA is highly sensitive to $\delta$, and that the dependence on $n$ can change substantially across the plotted risk levels. 
For $h=1$, $b=9$, and $\delta\in\{0.9,0.95,0.99\}$, worst-case regret remains above the constant-order benchmark even at $n=10{,}000$.

\begingroup\emergencystretch=3em
The analytical SAA series has been recomputed using exact critical ranks and refined one-dimensional maximization; see Appendix~\ref{apx:saa_numerics}.
\par\endgroup

\medskip
\FIGURE
{%
\begin{minipage}{\linewidth}
\centering
\includegraphics[width=443.565bp,height=164.041bp]{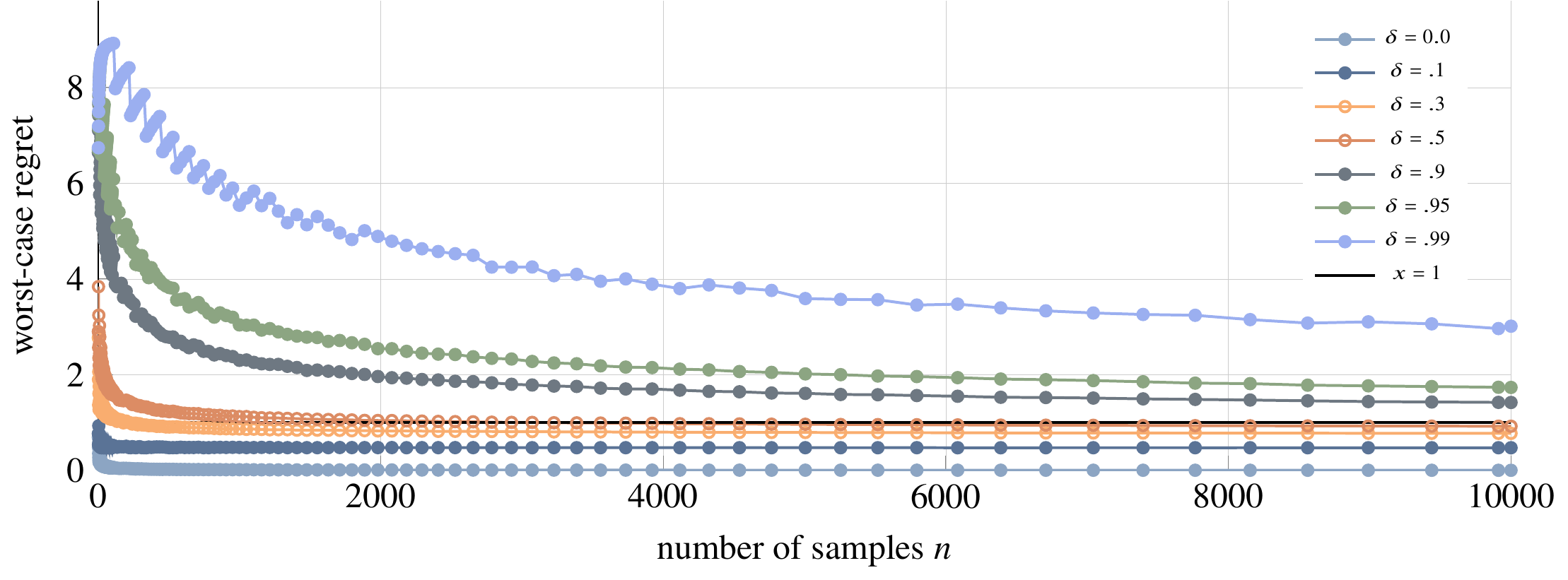}\\
{\small $h = 1$, $b=9$.}
\end{minipage}
}
{Numerical approximations to SAA worst-case regret by $\delta$.\label{fig:saa_performance_by_delta_linear}}
{}

\medskip

\FloatBarrier
These comparisons quantify the scope for improving on SAA. We now use the same Bernoulli reduction to optimize the weights within MCC.

\subsection{Computation of a Near-Optimal MCC Policy}\label{sec:computed_policy}

We now turn from worst-case evaluation to policy optimization. By
Theorem~\ref{thrm:worst_case_bernoulli}, optimizing a deterministic
convex-combination policy over $\mathcal F_1$ is equivalent to optimizing its
worst-case regret over the one-parameter Bernoulli family. 
The search over demand distributions is therefore one-dimensional; policy optimization still chooses the coefficient vector.

We say that an MCC policy is \emph{deterministic} if its convex-combination coefficients are predefined and have no randomness: there exists a fixed vector $(\alpha_1,\ldots,\alpha_n) \in \Delta^{n-1}$, chosen ex ante and independent of the sample and of any data-generating randomness, such that the policy prescribes 
\[ x^{\pi} ~~=~~ \sum_{i=1}^n \alpha_i \cdot D_{(i)} \ , \] 
for every sample $D_1^n$. Equivalently, the weights of the order statistics do not depend on the sample and involve no internal randomization.
For $i=1,\ldots,n$, define the cumulative tail weights
\[
    A_i := \sum_{j=i}^n \alpha_j.
\]
Under Bernoulli demand, $A_i$ is the order selected when the sample contains exactly $i-1$ zeros.

\paragraph{Optimality of the deterministic policy with respect to MCC.} We first show that restricting attention to deterministic MCC policies is without loss, since for any random MCC, there exists a deterministic MCC whose performance is at least as good. Randomization over weights cannot hedge that tail risk: replacing the random weights by their means removes convex-order dispersion and cannot increase the combined-risk $\cvar$. The following claim formalizes this idea.
\begin{claim}\label{clm:optimality_deterministic_mcc}
    Let $\pi$ be an MCC policy. Define the deterministic policy $\bar{\pi}$ by the fixed weights $\bar{\alpha}_k := \mathbb{E}[\alpha_k]$, so that $\bar{x}(D_1^n) := \sum_{k=1}^n \bar{\alpha}_k\, D_{(k)}$ for every sample $D_1^n$. Then $\bar{\pi}$ is deterministic and, for every $F \in \mathcal{F}_1$ (evaluated under the same demand distribution $F$ on both sides),
        \[
            \RegretCVaR{F}{\pi} ~~\geq~~ \RegretCVaR{F}{\bar{\pi}} \, .
        \]
        Consequently, $\sup_{F \in \mathcal{F}_1} \RegretCVaR{F}{\pi} \geq \sup_{F \in \mathcal{F}_1} \RegretCVaR{F}{\bar{\pi}}$.
\end{claim}
The proof of the claim is in Appendix~\ref{prf:optimality_deterministic_mcc}. Consequently, 
\[
\inf_{\pi\in\Pi^{\mcc}}\sup_{F\in\mathcal{F}_1}\RegretCVaR{F}{\pi}
=
\inf_{\alpha\in\Delta^{n-1}}\sup_{F\in\mathcal{F}_1}\RegretCVaR{F}{\pi^\alpha}.
\]

\paragraph{A near-optimal MCC policy.} We now show that we can compute a near-optimal deterministic MCC policy, which by the reduction just recorded is a near-optimal MCC policy. To do so, we show that in fact, the problem of finding the optimal $\alpha \in \Delta^{n-1}$ is a semi-infinite linear program. We begin by showing how to compute the worst case regret of a deterministic policy.

We compute the CVaR of a deterministic convex combination under Bernoulli demand as follows. The key is that the ordering of the realized losses does not depend on the coefficient vector. Once we fix the Bernoulli parameter, the tail-selection weights are fixed and regret is affine in the policy weights.
\begin{enumerate}
    \item Sort the $(D_1^n,D)$ events in descending order of associated costs:
    \[
    \left[
    \begin{aligned}
    &\left(D_{(1)}=1,D=0\right),\ \left(D_{(1)}=0,D_{(2)}=1,D=0\right),\ \ldots,\ \left(D_{(n-1)}=0,D_{(n)}=1,D=0\right), \\
    &\left(D_{(n)}=0\right),\ \left(D_{(n)}=1,D_{(n-1)}=0,D=1\right),\ \ldots,\ \left(D_{(1)}=1,D=1\right)
    \end{aligned}
    \right] \ .
    \]
    Define the events as $O_i := (D_{(i-1)}=0, D_{(i)}=1,D=0)$ for $i=1,\ldots,n$ (with the convention that $D_{(0)}=0$), $U_1 = (D_{(n)}=0)$ and $U_i := (D_{(n-i+2)}=1, D_{(n-i+1)}=0,D=1)$ for $i=2,\ldots,n+1$ (with the convention that $D_{(n+1)}=1$). The event $U_1$ does not condition on the future demand because $D_{(n)}=0$ implies $x^\pi=0$, and hence $c(0,D)=0$ for both $D=0$ and $D=1$.

    \item Define the associated costs of $O_i$ and $U_i$:
    \[
    c\left(O_i\right) := hA_i, \qquad c\left(U_1\right) := 0, \qquad c\left(U_i\right) := -bA_{n-i+2}\quad \text{for } i=2,\ldots,n+1.
    \]
    Since $\alpha_j\geq 0$ for every $j$, these costs satisfy
    \[
    c(O_1)\geq\cdots\geq c(O_n)\geq c(U_1)=0\geq c(U_2)\geq\cdots\geq c(U_{n+1}).
    \]
    In particular, this ordering is independent of the coefficient vector $\alpha$.
\end{enumerate} 
It remains to compute the event probabilities and fill the CVaR tail; the resulting formulas are given in Appendix~\ref{apx:bernoulli_event_enumeration}. 
For fixed $(h,b,\delta)$, the probabilities, tail-selection weights, and oracle term depend on $\mu$ and are independent of $\alpha$. Since the event costs are linear in the cumulative weights $A_i$, the regret is affine in $\alpha$ for each fixed $\mu$.
Combining this affine dependence with Theorem~\ref{thrm:worst_case_bernoulli}, the
minimax problem over deterministic coefficient vectors can be written as
\[
    \inf_{\alpha\in\Delta^{n-1}}
    \sup_{F\in\mathcal F_1}
    \RegretCVaR{F}{\pi^\alpha}
    =
    \inf_{\alpha\in\Delta^{n-1}}
    \sup_{\mu\in[0,1]}
    \RegretCVaR{\mathrm{Ber}(\mu)}{\pi^\alpha}.
\]
The right-hand side is a semi-infinite linear program: the decision variables
are the weights $\alpha$, and the constraints enforce a regret upper bound for
each Bernoulli parameter $\mu\in[0,1]$.
Therefore an optimal deterministic convex combination of order statistics is a solution of the following semi-infinite linear program:
\begin{equation} \label{eqn:semi_infinite_lp}
    \begin{aligned}
        \min_{\alpha \in \mathbb{R}^n, t \in \mathbb{R}} \quad & t \\
        \textrm{s.t.} \quad & \cvar_\delta[c\left(x^{\alpha}, D\right)] - \opt_\delta\left(\mathrm{Ber}\left(\mu\right)\right) \leq t \quad & \forall \mu \in [0,1] \\
        \quad & \sum_{i=1}^n \alpha_i = 1 \\
        \quad & \alpha_i \geq 0 \quad & \forall i \in [n]
    \end{aligned}
\end{equation}
In practice, we want to avoid an infinite number of constraints. One way to approximate the semi-infinite program is to solve a discretized version: for some $m \in \mathbb{N}$, we can solve
\begin{equation*}
    \begin{aligned}
        \min_{\alpha \in \mathbb{R}^n, t \in \mathbb{R}} \quad & t \\
        \textrm{s.t.} \quad & \cvar_\delta[c\left(x^{\alpha}, D\right)] - \opt_\delta\left(\mathrm{Ber}\left(\mu\right)\right) \leq t \quad & \forall \mu \in \left\{0, \frac{1}{m}, \ldots, \frac{m-1}{m}, 1\right\} \\
        \quad & \sum_{i=1}^n \alpha_i = 1 \\
        \quad & \alpha_i \geq 0 \quad & \forall i \in [n]
    \end{aligned}
\end{equation*}
The discretization error can be controlled by an upper bound on the Lipschitz constant of the function $\mu \mapsto \cvar_\delta[c(x^{\pi}, D)] - \opt_\delta(\mathrm{Ber}(\mu))$. The following claim formalizes this.
\begin{claim}\label{clm:lipschitz_constant}
    Fix a deterministic coefficient vector $\alpha\in\Delta^{n-1}$. For each $\mu\in[0,1]$, let $D_1^\mu,\ldots,D_n^\mu,D_\mu$ be i.i.d.\ $\mathrm{Ber}(\mu)$, let $D_{(1)}^\mu\leq \cdots \leq D_{(n)}^\mu$ denote the order statistics of $D_1^\mu,\ldots,D_n^\mu$, and define
    \[
    x^\pi_\mu:=\sum_{i=1}^n \alpha_iD_{(i)}^\mu .
    \]
    Then the function $\mu \mapsto \cvar_\delta[c(x^\pi_\mu, D_\mu)] - \opt_\delta(\mathrm{Ber}(\mu))$ is Lipschitz continuous; in particular, one valid Lipschitz constant is $(h+b) \cdot (n+2)/(1-\delta)$.
\end{claim}
The proof is in Appendix~\ref{prf:lipschitz_constant}. \footnote{When we compute the optimal policy in practice, we use a constraint-generating linear programming scheme to avoid loading all grid constraints at once. We elaborate on the implementation in Appendix~\ref{apx:implementation}.}
Let $L_n:=(h+b)(n+2)/(1-\delta)$, let $V_n^{\mathrm{det}}$ denote the semi-infinite LP value, and let $(\widehat{\alpha}_m,t_m)$ solve its grid version. Then
\[
\sup_{\mu\in[0,1]}\RegretCVaR{\mathrm{Ber}(\mu)}{\pi^{\widehat{\alpha}_m}}
\leq t_m+\frac{L_n}{2m}
\leq V_n^{\mathrm{det}}+\frac{L_n}{2m}.
\]
The upper bound follows by evaluating at a nearest grid point and applying Claim~\ref{clm:lipschitz_constant}.
The grid program relaxes~\eqref{eqn:semi_infinite_lp}, so $t_m\leq V_n^{\mathrm{det}}$. Together with the display,
\[
t_m
~~\leq~~
V_n^{\mathrm{det}}
~~\leq~~
\sup_{\mu\in[0,1]}\RegretCVaR{\mathrm{Ber}(\mu)}{\pi^{\widehat{\alpha}_m}}
~~\leq~~
t_m+\frac{L_n}{2m}.
\]

\subsection{\texorpdfstring{Performance of the Computed MCC Policies}{Performance of the Computed MCC Policies}}
We compare SAA with the computed policy from Section~\ref{sec:computed_policy}, and keep the deterministic benchmark that always orders $x=1$ for reference.
Figure~\ref{fig:finite_sample_opt_vs_saa} reports numerical approximations to worst-case regret for $h=1$, $b=9$, and $\delta=.95$.
At $n=100$ the values are $5.80$ for SAA and $0.84$ for the computed policy; at $n=10{,}000$ they are $1.74$ and $0.41$.
\begingroup\emergencystretch=3em
The analytical SAA series has been recomputed using exact critical ranks and refined one-dimensional maximization; see Appendix~\ref{apx:saa_numerics}.
\par\endgroup
The figure shows that the computed deterministic policy 
substantially improves on SAA and improves on the constant-order benchmark over a broad range of sample sizes, while matching the latter at the smallest plotted sizes.

\medskip
\FIGURE
{%
\begin{minipage}{\linewidth}
\centering
\includegraphics[width=421.463bp,height=184.019bp]{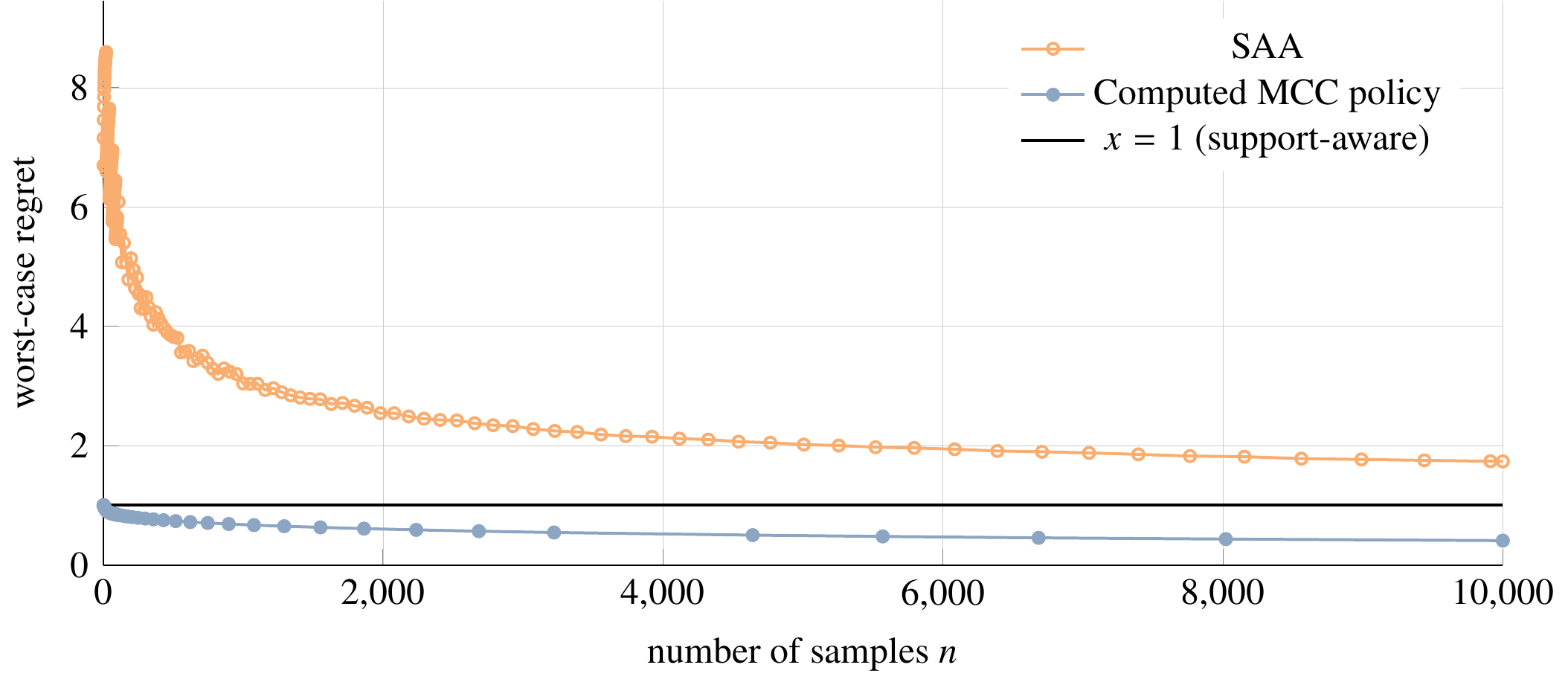}\\
{\small $h = 1$, $b=9$, $\delta = .95$.}
\end{minipage}
}
{Numerical worst-case comparison: refined SAA maxima and grid-evaluated maxima for the computed MCC policy.\label{fig:finite_sample_opt_vs_saa}}
{}

The comparison in Figure~\ref{fig:finite_sample_opt_vs_saa} shows that SAA is not merely a loose benchmark for the focal class; it can be separated from the best computed policies by a substantial margin. This gap suggests that strong worst-case performance may require policies that depart meaningfully from the SAA rule, rather than small perturbations of it. 

We move on to analyze the decomposition of the combined-risk regret into expected post-sample regret and sample-risk premium. 
Figure~\ref{fig:near_opt_psi_decomp} in Appendix~\ref{apx:additional_figures} shows combined-risk regret and expected post-sample regret at the grid distribution that maximizes combined-risk regret. 

Minimizing worst-case combined-risk regret can leave substantial expected post-sample regret. The decomposition suggests that we may improve this component while allowing a small increase in the worst-case combined-risk objective.
Figure~\ref{fig:benign_near_opt} in Appendix~\ref{apx:additional_figures} shows these weights on uniform demand and on exponential demand with rate $1$, the two distributions shared with the later comparison. On both, the near-optimal policy has higher combined-risk regret than SAA over nearly the whole plotted range.
\medskip

We therefore propose a two-stage linear programming procedure. We first minimize worst-case combined-risk regret. We then minimize worst-case expected post-sample regret subject to a prescribed slack in the first objective.
 In the numerical experiments, at the grid distribution that maximizes combined-risk regret, this procedure substantially reduces expected post-sample regret, with an optimization slack of $0.01$\footnote{Appendix~\ref{apx:two_stage_procedure} gives the formulation and implementation. Theorem~\ref{thrm:worst_case_bernoulli} reduces both worst-case criteria to Bernoulli demand.}. The resulting policy's worst-case combined-risk regret exceeds the MCC minimax value by at most the prescribed slack plus numerical approximation errors.

The Bernoulli parameter that binds for one need not bind for the other.
Figure~\ref{fig:psi_decomp} in Appendix~\ref{apx:additional_figures} shows the two-stage policy decomposition at a grid distribution maximizing combined-risk regret, 
with optimization slack $0.01$. This slack is separate from the first-stage approximation error, solver tolerances, and discretization error.

The comparisons below are Monte Carlo estimates of combined-risk regret for the balancing policy and for SAA on shared draws. SAA places unit weight on the order statistic of rank $\lceil nb(1-\delta)/(h+b)\rceil$. Each series records the standard error of each estimate, and the standard error of the paired difference. Appendix~\ref{apx:fixed_demand_mc} records the computation.

Finally, the balancing policy improves combined-risk regret relative to SAA over a substantial range of sample sizes.
The four distributions are again those of Figure~\ref{fig:saa_decomp_distributions}, each with a unique oracle order. 
The improvement is substantial across a broad range of the scenarios considered, although it is not uniform.

\medskip
\FIGURE
{%
\begin{minipage}{\linewidth}
\centering
    \begin{minipage}{0.48\linewidth}
        \centering
        \includegraphics[width=225.339bp,height=118.526bp]{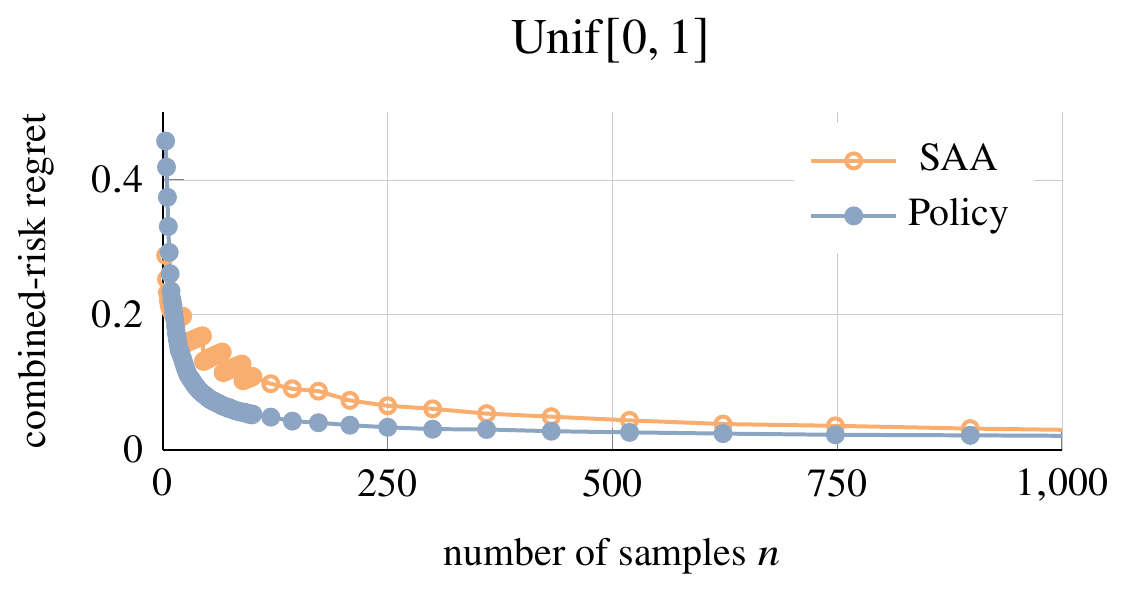}
    \end{minipage}
    \hfill
    \begin{minipage}{0.48\linewidth}
        \centering
        \includegraphics[width=225.339bp,height=118.724bp]{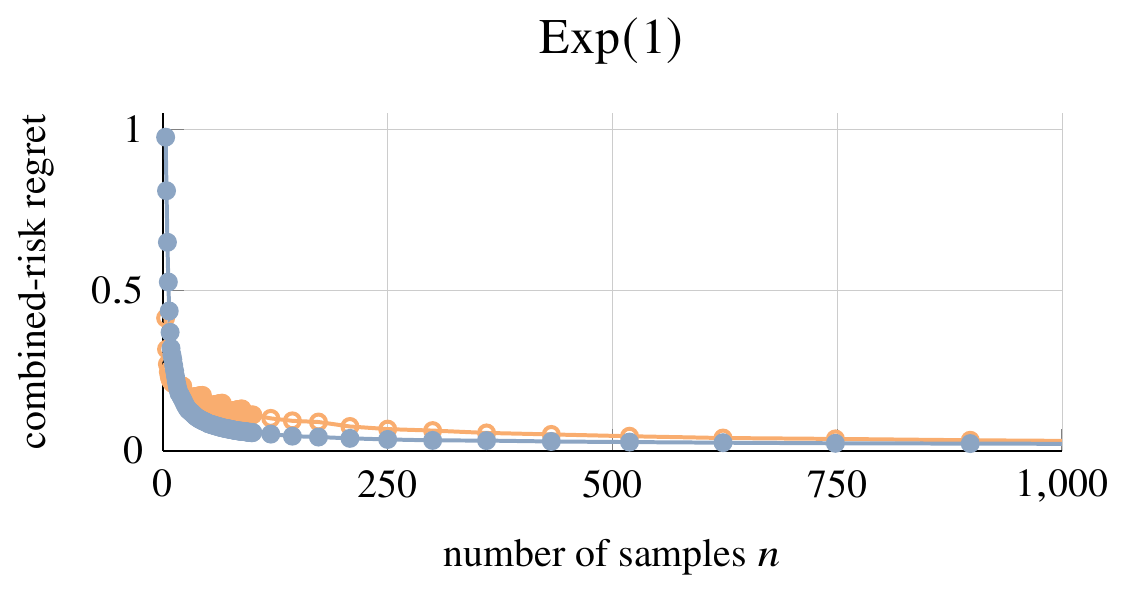}
    \end{minipage}\\[0.8em]
    \begin{minipage}{0.48\linewidth}
        \centering
        \includegraphics[width=218.796bp,height=120.201bp]{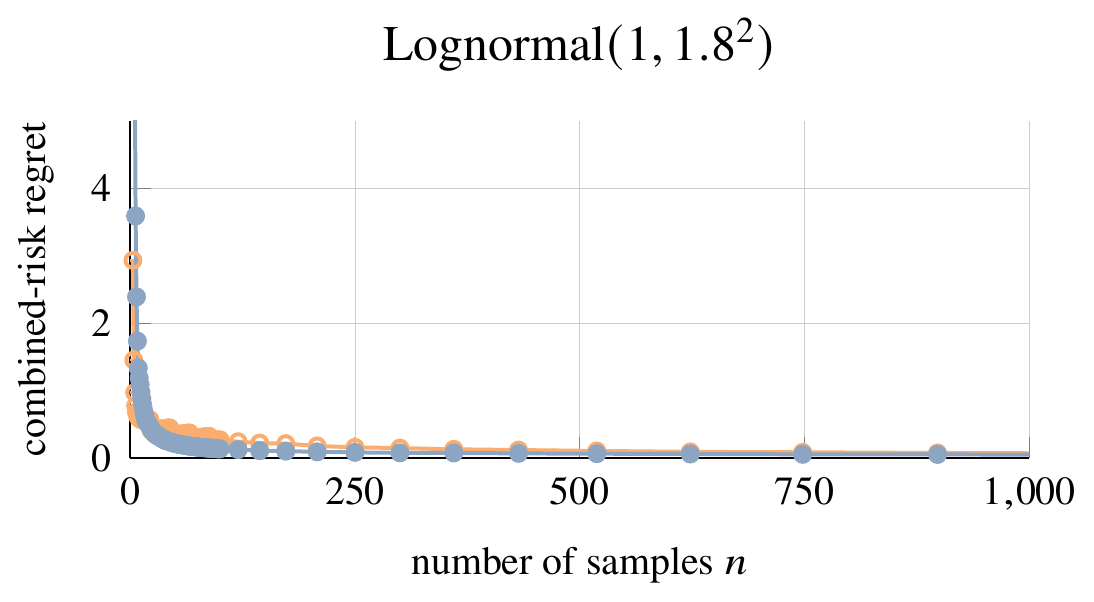}
    \end{minipage}
    \hfill
    \begin{minipage}{0.48\linewidth}
        \centering
        \includegraphics[width=225.339bp,height=118.523bp]{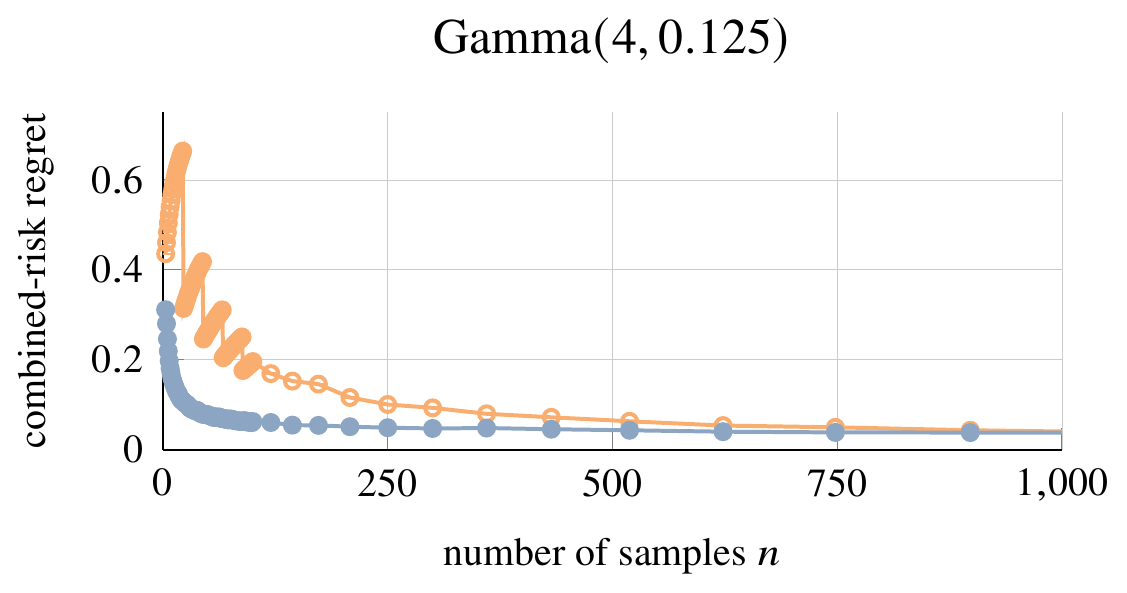}
    \end{minipage}\\
{\small $h = 1$, $b=9$, $\delta = .95$.}
\end{minipage}
}
{Monte Carlo estimates of combined-risk regret for SAA and the balancing policy under uniform demand on $[0,1]$, exponential demand with rate $1$, lognormal demand with $\mu_{\log}=1$ and $\sigma_{\log}=1.8$, and gamma demand with shape $4$ and scale $0.125$. Over $n\le 1000$, the Monte Carlo standard errors are at most $0.0012$ under uniform and exponential demand, $0.0021$ under gamma demand, and $0.035$ under lognormal demand. The lognormal axis ends at $5$; the policy's combined-risk regret at $n=3$, $4$, and $5$ is approximately $17.6$, $11.6$, and $6.1$. Appendix~\ref{apx:fixed_demand_mc} records the computation.\label{fig:saa_vs_policy_distributions}}
{}

\medskip

%% file: TEX-SAA.tex
\section{\texorpdfstring{Minimax Learning Rates and Stabilized Policies}{Minimax Learning Rates and Stabilized Policies}}\label{sec:mcc_asymptotic_success}

The decomposition in Section~\ref{sec:decomposition} identifies sample-induced switching in the learned action as the mechanism behind SAA's combined-risk failure (Theorem~\ref{thm:saa_nonvanishing}). At the critical Bernoulli instance, the $O(\sqrt{n})$ binomial fluctuation in the sample count causes SAA's single empirical rank to switch between the two corner actions with nonvanishing probability. This switching leaves the post-sample regret at zero but generates a persistent sample-risk premium.

We therefore average order statistics over a window centered at the critical rank $n\qstar$ and wider than the $\sqrt{n}$ scale of the binomial fluctuations.

Let $f(n)$ satisfy $f(n)=\omega(\sqrt{n})$ and $f(n)=o(n)$. For each $n$, define
\[
I_n:=\left\{i\in\{1,\ldots,n\}: |i-n\qstar|\leq f(n)\right\}.
\]
A simple choice assigns equal weight to every rank in $I_n$: $\alpha_i^{(n)}=\mathbf{1}_{\{i\in I_n\}}/|I_n|$. We also allow nonuniform weights, provided no individual rank receives too much weight.
More generally, let $\boldsymbol{\alpha}^{(n)}\in\Delta^{n-1}$ satisfy $\alpha_i^{(n)}=0$ for $i\notin I_n$ and $\sqrt{n}\|\boldsymbol{\alpha}^{(n)}\|_\infty\to0$, and let $\pi^{\mathrm{win}}$ be the deterministic MCC policy
\[
x^{\mathrm{win}}(D_1^n)=\sum_{i=1}^n \alpha_i^{(n)}D_{(i)}.
\]
Under Bernoulli demand, a change of $k$ in the number of sample zeros changes the order by at most $k\|\boldsymbol{\alpha}^{(n)}\|_\infty$. The maximum-weight condition therefore makes the effect of fluctuations on the scale $\sqrt n$ vanish.

The following theorem shows that this windowed policy achieves vanishing worst-case regret.

\begin{theorem}\label{thm:mcc_vanishing}
Let $\pi^{\mathrm{win}}$ be the window policy defined above. If
$f(n)=\omega(\sqrt{n})$, $f(n)=o(n)$, and
$\sqrt{n}\|\boldsymbol{\alpha}^{(n)}\|_\infty\to0$, then
\[
\lim_{n\to\infty}
\sup_{F\in\mathcal{F}_1}
\RegretCVaR{F}{\pi^{\mathrm{win}}}=0.
\]

If we further suppose that the weights are window-flat, meaning that no coefficient exceeds a constant times $1/f(n)$ for all sufficiently large $n$:
\[
\|\boldsymbol{\alpha}^{(n)}\|_\infty=O\!\left(\frac{1}{f(n)}\right),
\]
then
\[
    \Omega\!\left(\frac{f(n)}{n}\right)
    \leq
    \sup_{F\in\mathcal{F}_1}\RegretCVaR{F}{\pi^{\mathrm{win}}}
    \leq
    O\!\left(
        \frac{f(n)}{n}
        +\frac{\sqrt n}{f(n)}
    \right).
    \]
    Consequently, $f(n)=n^{3/4}$ yields
    \[
    \sup_{F\in\mathcal{F}_1}
    \RegretCVaR{F}{\pi^{\mathrm{win}}}
    =
    \Theta(n^{-1/4}).
    \]
\end{theorem}

The proof is provided in Appendix~\ref{prf:mcc_vanishing}. The theorem establishes that SAA's asymptotic failure is not a limitation of the MCC class. 
Uniform weights with $f(n)=n^{3/4}$ yield worst-case regret $\Theta(n^{-1/4})$, as do the more general window-flat weights in the theorem. The following theorem shows that their rate is optimal over all data-driven policies. 

\begin{samepage}
\begin{theorem}\label{thm:global_minimax_rate}
Fix $h,b,T>0$. At each risk level $\delta$, define
\[
\displaystyle
\begin{aligned}
V_n^{\mathrm{all}}(T)&:=\inf_{\pi\in\Pi}\sup_{F\in\mathcal F_T}\RegretCVaR{F}{\pi},\\
V_n^{\mcc}(T)&:=\inf_{\pi\in\Pi^{\mcc}}\sup_{F\in\mathcal F_T}\RegretCVaR{F}{\pi}.
\end{aligned}
\]
For every fixed $\delta\in(0,1)$, there are constants $0<c_0\leq C_0<\infty$, depending only on $(h,b,\delta)$, such that, for all sufficiently large $n$,
\[
\displaystyle
c_0Tn^{-1/4}\leq V_n^{\mathrm{all}}(T)
\leq V_n^{\mcc}(T)\leq C_0Tn^{-1/4}.
\]
Uniform window policies with $f(n)=n^{3/4}$ attain this order. In contrast, at $\delta=0$, both minimax values have order $\Theta(Tn^{-1/2})$.
\end{theorem}
\end{samepage}

The proof is in Appendix~\ref{prf:global_minimax_rate}. The lower bound holds for all sufficiently large $n$ and all policies, even when $T$ is known. Small regret requires the mean order to change across nearby Bernoulli distributions, while the sample-risk premium requires the order to concentrate at each distribution. Samples from nearby parameters have overlapping distributions, which limits how quickly a concentrated order can change. This tension yields the $n^{-1/4}$ lower bound. Together with the window upper bound, it shows that positive risk aversion changes the best achievable learning rate, in addition to making SAA inconsistent.

The binomial count fluctuates on the scale $\sqrt{n}$. Under $f(n)=\omega(\sqrt{n})$ and $\sqrt{n}\|\boldsymbol{\alpha}^{(n)}\|_\infty\to 0$, including a width only slightly larger than $\sqrt{n}$, the learned order concentrates around a deterministic, possibly $n$-dependent proxy. At the critical Bernoulli 
every order in $[0,1]$ is optimal, so post-sample regret is identically zero and the sample-risk premium vanishes as the learned order concentrates. 
For uniform weights on $I_n$, the proxy converges to $1/2$. 
At finite sample sizes, truncation at ranks $1$ and $n$, as well as integer rounding, can move the proxy away from $1/2$. Concentration does not require the balanced width $f(n)=n^{3/4}$.

The upper bound separates the suboptimality of the deterministic proxy, controlled by $f(n)/n$, from the effect of sampling variation around it, controlled by $\sqrt n/f(n)$ under window-flat weights. A wider window increases the first bound and reduces the second. Balancing these bounds gives $f(n)=n^{3/4}$. The requirement $f(n)=o(n)$ makes the bound on the proxy's suboptimality vanish.

Under the regularity and concentration conditions accompanying equation~\eqref{eqn:concentration}, a policy with concentration rate $r_n=n^{-1/2}$ has combined-risk regret of order $n^{-1/2}$ at the fixed demand distribution. The comparison is between this window's uniform rate and that pointwise rate. It does not show that these windows are slower at that distribution.
A window whose width stays a fixed fraction of $n$, and whose weights are spread over that width, 
can retain positive worst-case regret asymptotically. A wide declared window can still put its mass on a smaller valid subwindow.

The window policy also certifies that the MCC minimax value vanishes as the sample size grows. Consequently, computed near-minimax policies achieve vanishing worst-case regret if their approximation error tends to zero. For the two-stage policies of Section~\ref{sec:order_stats}, which allow a fixed slack of $0.01$ above the minimax value, vanishing approximation error instead guarantees that their worst-case regret is asymptotically at most $0.01$.

\section{Conclusion}

We show that the move from expected cost to combined-risk CVaR changes the learning guarantees of the data-driven newsvendor. For every fixed positive risk level, SAA can retain regret bounded away from zero at a fixed demand distribution, even as the sample grows. The regret decomposition identifies sample-induced variability in the learned order as the source of this failure: in our construction, expected post-sample regret is zero for every sample size, while the sample-risk premium remains positive asymptotically. Under our regularity and concentration conditions, both components vanish, with the sample-risk premium governing the rate of combined-risk regret.

We also construct explicit window policies that stabilize the learned order by averaging nearby order statistics. With uniform weights and window half-width $f(n)=n^{3/4}$, these policies achieve worst-case regret $\Theta(n^{-1/4})$. We prove a matching lower bound over all data-driven policies. Thus, these windows attain the minimax rate, which changes from $\Theta(n^{-1/2})$ under expectation to $\Theta(n^{-1/4})$ at every fixed positive risk level. Additional data can therefore support vanishing regret under combined-risk CVaR, even when empirical optimization fails to realize that improvement. Our results identify control of sample-induced variation in the learned decision as a central principle for realizing the value of data.

\begingroup\emergencystretch=3em
Within the MCC class, we reduce worst-case analysis over $\mathcal F_T$ to Bernoulli demand on $\{0,T\}$. We use this reduction to compute minimax policies in the class to any prescribed accuracy through linear programming. The decomposition also guides a two-stage procedure: we first minimize worst-case combined-risk regret, then minimize worst-case expected post-sample regret subject to a prescribed slack in the first objective. In our numerical experiments, the resulting policies substantially improve on SAA across a broad range of demand distributions and sample sizes, although the improvement is not uniform.
\par\endgroup

We leave open whether policies outside MCC can improve finite-sample performance.

%% file: TEX-Appendix.tex
\section{CVaR auxiliary claims and proofs}\label{apx:cvar_auxiliary}

\subsection{Proof of Lemma~\ref{lem:nv_atom}}\label{prf:nv_atom}

\begin{proof}{Proof}
The identity $c(x,D)=-bx+(h+b)(x-D)^+$ gives $c(x,D)=-bx$ on $\{D\geq x\}$ and $c(x,D)=-bx+(h+b)(x-D)>-bx$ on $\{D<x\}$. The distribution of the loss therefore 
places probability $1-F(x-)$ at $-bx$, and no mass below $-bx$.

Since $\qstar=b(1-\delta)/(h+b)$ and $h>0$, one has $\qstar<1-\delta$. At $x^*_{F}=F^{-1}(\qstar)$, every distribution function satisfies $F(x^*_{F}-)\leq\qstar$, so
\[
    \mathbb{P}\bigl(c(x^*_{F},D)=-bx^*_{F}\bigr) ~~=~~ 1-F(x^*_{F}-) ~~\geq~~ 1-\qstar ~~>~~ \delta,
    \qquad
    \mathbb{P}\bigl(c(x^*_{F},D)<-bx^*_{F}\bigr) ~~=~~ 0.
\]
The smallest threshold at which the loss distribution has accumulated mass $\delta$ is therefore $-bx^*_{F}$. That is $\mathrm{VaR}_\delta[c(x^*_{F},D)]$, and the same two displays are the strict straddling condition.
\Halmos
\end{proof}

If in addition $F$ is continuous, $x^*_{F}$ lies in the interior of $\mathcal{X}$ and $0<F'(x^*_{F})<\infty$, then 
$F(x^*_{F})=\qstar$ and the atom mass equals $1-\qstar$. Continuity at $x^*_{F}$ and $\qstar<1-\delta$ give $F(x)<1-\delta$ throughout a neighborhood of $x^*_{F}$, and on that neighborhood $\mathrm{VaR}_\delta[c(x,D)]=-bx$. By \cite[Theorem 10]{RockafellarU02} the variational formula is attained at this threshold. Using $c(x,D)+bx=(h+b)(x-D)^+$,
\[
    \cvar_\delta[c(x,D)] ~~=~~ -bx + \frac{h+b}{1-\delta}\,\mathbb{E}\bigl[(x-D)^+\bigr].
\]
The derivative of $\mathbb{E}[(x-D)^+]$ is $F(x)$, so the derivative of the display is $-b+\frac{h+b}{1-\delta}F(x)$. It vanishes at $x^*_{F}$ because $F(x^*_{F})=\qstar$. Differentiating once more produces $(h+b)F'(x^*_{F})/(1-\delta)\in(0,\infty)$.

We next verify the probability limits required by the isolation condition. Write $\ell(x)=-bx$. The identity $c(x,D)=\ell(x)+(h+b)(x-D)^+$ puts every loss at or above $\ell(x)$, so for every $u>0$,
\[
    \mathbb{P}\bigl(c(x,D)\ge\ell(x)-u\bigr)
    ~~=~~
    1
    ~~=~~
    \mathbb{P}\bigl(c(x^*_{F},D)\ge -b x^*_{F}\bigr).
\]
A loss strictly above $\ell(x)+u$ is a demand that falls short of the order by more than $u/(h+b)$,
\[
    \mathbb{P}\bigl(c(x,D)>\ell(x)+u\bigr)
    ~~=~~
    \mathbb{P}\bigl(D<x-u/(h+b)\bigr).
\]
As $(x,u)\to(x^*_{F},0^+)$ the threshold tends to $x^*_{F}$. Continuity of $F$ gives $\mathbb{P}(D<y)\to F(x^*_{F})$ as $y\to x^*_{F}$, and the same continuity gives $\mathbb{P}(c(x^*_{F},D)>-bx^*_{F})=F(x^*_{F})$. These are the isolation limits in Theorem~\gref{cor:atom_category} of \cite{BesbesRiegerGeneral}. The loss is differentiable in the order except on $\{D=x\}$. That event has probability zero, and on its complement the derivative equals $h$ or $-b$.

The loss is also Lipschitz in the order with constant $\max(b,h)$, independent of demand, and the atom level $\ell(x)=-bx$ is linear in the order. On a compact order set with bounded demand, the loss is uniformly bounded. Together with the almost-everywhere derivative and the Lipschitz constant, which has every moment, these properties verify the conditions on the loss in Proposition~\gref{prop:saa_concentration_rate} of \cite{BesbesRiegerGeneral}.

In the proof of Theorem~\gref{cor:atom_category} of \cite{BesbesRiegerGeneral}, the one-sided slopes of the variational objective are the quantities in equation~\gref{eqn:atom_slopes}. For this loss they equal
\[
    s_+ ~~=~~ \frac{1-\qstar-\delta}{1-\delta} ~~=~~ \frac{h}{h+b},
    \qquad
    s_- ~~=~~ \frac{\delta}{1-\delta}.
\]
Their sum equals $(1-\qstar)/(1-\delta)$. The cutoff equals $-bx$, so transferring the weighted deviation to the order multiplies by $b$ and exchanges the sides. The resulting constant is the limit in~\eqref{eqn:sr_approx},
\[
    \frac{b\,(h+b\delta)}{(h+b)\,(1-\delta)} \ .
\]
Dividing by the limit in~\eqref{eqn:atom_approx} produces the limit in~\eqref{eqn:sr_over_p}.

\subsection{Differences between Conditional Values at Risk}
\label{prf:bound_cvar_diff}

\begin{lemma}\label{lem:bound_cvar_diff}
For any random variables $Z$ and $Z'$ with finite first moments,
\[
|\cvar_\delta\left(Z\right) - \cvar_\delta\left(Z'\right)| \;\leq\; \frac{1}{1-\delta}\,\mathbb{E}|Z - Z'|.
\]
\end{lemma}

\begin{proof}{Proof}
Using the following representation:
\[\cvar_\delta\left(Z\right) ~~=~~ \inf_{t\in\mathbb{R}}\left\{t + \left(1-\delta\right)^{-1}\mathbb{E}\left[\left(Z-t\right)^+\right]\right\} \ , \]
for each $t$, the map $z\mapsto(z-t)^+$ is $1$-Lipschitz, so $|(Z-t)^+-(Z'-t)^+|\leq|Z-Z'|$ and thus
\[
\left|t + \frac{1}{1-\delta}\,\mathbb{E}\left[\left(Z-t\right)^+\right] - \left(t + \frac{1}{1-\delta}\,\mathbb{E}\left[\left(Z'-t\right)^+\right]\right)\right|
\;\leq\; \frac{1}{1-\delta}\,\mathbb{E}|Z-Z'|.
\]
Let $\varepsilon>0$ and let $t_\varepsilon$ be such that the variational objective of $Z$ at $t_\varepsilon$ is at most $\cvar_\delta(Z)+\varepsilon$. The same objective for $Z'$ at $t_\varepsilon$ is at most that value plus $(1-\delta)^{-1}\mathbb{E}|Z-Z'|$, and it is at least $\cvar_\delta(Z')$. The symmetric choice, starting from $Z'$, gives the other direction. Letting $\varepsilon\downarrow 0$ yields
$|\cvar_\delta(Z)-\cvar_\delta(Z')|\leq (1-\delta)^{-1}\mathbb{E}|Z-Z'|$.
\Halmos
\end{proof}

\subsection{Mismatch Cost-CVaR and Profit-CVaR}\label{apx:cost_profit_cvar}
Recall the cost $c(x,D)$ from equation~\eqref{eqn:neg_profit}.
We write the \emph{demand-supply mismatch cost}
\begin{equation}\label{eqn:mismatch_cost}
    \mcost\left(x,D\right) ~~=~~ b \cdot \left(D-x\right)^+ + h \cdot \left(x-D\right)^+ \ ,
\end{equation}
so that $c(x,D) = \mcost(x,D) - bD$.
For a fixed demand distribution $F$, two natural formulations arise under $\cvar$ (with $D \sim F$ throughout):
\begin{enumerate}
    \item \emph{CVaR of demand-supply mismatch cost:} minimize $\cvar_\delta\bigl[\mcost(x,D)\bigr]$.
    \item \emph{Profit-CVaR:} minimize $\cvar_\delta\bigl[c(x,D)\bigr]$, where $c(x,D)$ is the cost defined in equation~\eqref{eqn:neg_profit}.
\end{enumerate}
Minimizing $\cvar_\delta[c]$ is equivalent to protecting against the worst $(1-\delta)$-tail of profits: the loss $c$ attains its largest values when profit $\Phi(x,D)$ attains its smallest, so the $\cvar$ tail of the loss corresponds to the worst tail of profits.

These two formulations coincide under expectation (when $\delta=0$) but generally differ under $\cvar$: under an expectation objective, the random term $-bD$ is constant in $x$, so although it shifts the objective value, it does not change the minimizer; under $\cvar$, by contrast, it \emph{does} shift which realizations populate the tail event defining $\cvar$, so the optimizer may differ even though the added term is decision-independent.
\cite{GotohT07} derive closed-form solutions for both problems when the demand is known.
For the $\cvar$ objective of demand-supply mismatch cost, $\cvar_\delta\bigl[\mcost(x,D)\bigr]$, an optimal solution is
\begin{equation}\label{eqn:gotoh_takano_mismatch_apx}
    x^*_{F} ~~=~~ \frac{h}{h+b} \cdot F^{-1} \left(\qstar \right) + \frac{b}{h+b} \cdot F^{-1} \left(\frac{h \cdot \delta + b}{h+b} \right) \ ,
\end{equation}
where $\qstar = b \cdot (1-\delta)/(h+b)$.
For the $\cvar$ objective of the \emph{profit} problem $\cvar_\delta\bigl[c(x,D)\bigr]$, an optimal solution is
\begin{equation}\label{eqn:gotoh_takano_profit_apx}
    x^*_{F} ~~=~~  F^{-1} \left(\qstar \right) \ .
\end{equation}
The contrast between equations~\eqref{eqn:gotoh_takano_mismatch_apx} and~\eqref{eqn:gotoh_takano_profit} illustrates that the mismatch-cost and profit formulations can yield distinct optimizers under $\cvar$.

\section{Incomparability with Counting Policies}\label{apx:mcc_counting_incomparable}
We now show that, for $n\geq 2$, the class $\Pi^{\mcc}$ is incomparable with the counting-policy class of \cite[Definition 4]{BesbesMaMouchtaki2023}.
We use their threshold-event characterization: for each threshold $y$, the event $\{x(D_1^n)\leq y\}$ is determined by the vector of indicators of which observations fall at or below $y$. The same function of these indicators must apply at every threshold.

\paragraph{$\Pi^{\mcc}$ is not contained in counting policies.}
Consider $n=2$ and the policy
\[
x\left(D_1,D_2\right) = \frac{D_{(1)} + D_{(2)}}{2},
\]
which is a deterministic convex combination of order statistics, hence belongs to $\Pi^{\mcc}$. For threshold $y=0.55$, compare two samples with demands in $[0,1]$:
\[
\left(0,1\right) \quad \text{and} \quad \left(0.3,1\right).
\]
For both samples the first observation is at most $0.55$ and the second is above it, so any counting-based representation uses the same vector of threshold indicators. Yet
\[
x(0,1) = 0.5,
\qquad
x(0.3,1) = 0.65.
\]
Thus the event $\{x \leq 0.55\}$ holds for $(0,1)$ but not for $(0.3,1)$, so this policy is not a counting policy.

\paragraph{Counting policies are not contained in $\Pi^{\mcc}$.}
Consider
\[
x\left(D_1,D_2\right) = \min\{D_{(1)},0\}.
\]
Since demands lie in $[0,1]$, this is the constant policy $x\equiv 0$, so $\{x\leq y\}$ holds for every $y\geq 0$ and is empty for $y<0$; the induced threshold events are counting-based. However, this policy cannot be represented as a convex combination of observed demands, because when $D_{(1)}>0$ it outputs $0$, which need not be in the convex hull of $\{D_1,D_2\}$. Therefore this policy is a counting policy but it is not in $\Pi^{\mcc}$.

Together, these two examples establish incomparability.
For any $n>2$, the first example extends by taking $x=(D_{(1)}+D_{(n)})/2$ and appending $n-2$ observations equal to one to each sample. The constant-zero policy provides the second example for every $n$. At $n=1$, however, MCC contains only the counting policy $x=D_1$, so the classes are not incomparable.

\section{Auxiliary Lemma Proofs for Theorem~\ref{thrm:worst_case_bernoulli}}

\subsection{Proof of Lemma~\ref{lem:discretization}} \label{prf:discretization_coupling_once}
Throughout, $F,F_m\in\mathcal{F}_1$, so losses are uniformly bounded on $[0,1]$ and Lemma~\ref{lem:bound_cvar_diff} transfers almost-sure loss convergence to regret convergence.
\paragraph{Coupling lemma.}
We begin by stating conditions that ensure that regret converges along a sequence of approximating distributions in $\mathcal{F}_1$.
\begin{lemma}\label{lem:coupling_regret_convergence}
Let $\pi^{\mcc}$ be an MCC policy.
Let $F \in \mathcal{F}_1$ and $(F_m)_{m\geq 1} \subseteq \mathcal{F}_1$.
Fix a probability space and, on it, a coupling of demand vectors: let $(D_1,\dots,D_n,D)\sim F$ and, for each $m\geq 1$, let $(D_1^{(m)},\dots,D_n^{(m)},D^{(m)})\sim F_m$.
Let $W$ be the exogenous randomness used by the policy to generate coefficients on the same space, and define
$\boldsymbol{\alpha} := \boldsymbol{\alpha}(W)$.
Assume:
\begin{enumerate}[label=\textup{(\roman*)}]
    \item for each $i$, $D_i^{(m)}\to D_i$ and $D^{(m)}\to D$ almost surely as $m\to\infty$; \label{itm:demand_convergence_coupling}
    \item $x^*_{F_m}\to x^*_F$ as $m\to\infty$; \label{itm:offline_covergence_coupling}
\end{enumerate}
Then, as $m\to\infty$,
\[
    \RegretCVaR{F_m}{\pi^{\mcc}} ~~\to~~ \RegretCVaR{F}{\pi^{\mcc}},
    \qquad
    \mathrm{PS}_{\delta,n}(\pi^{\mcc},F_m) ~~\to~~ \mathrm{PS}_{\delta,n}(\pi^{\mcc},F).
\]
\end{lemma}
The proof is in Appendix~\ref{prf:coupling_regret_convergence}.

\paragraph{Mesh lemma.}
Using Lemma~\ref{lem:coupling_regret_convergence}, we show that for any $F\in\mathcal{F}_1$, we can construct a discrete distribution with finite support such that, under a pointwise coupling, the sampled demands and oracle decision converge and the corresponding regret is appropriately close to that associated with~$F$.

\begin{lemma}\label{lem:mesh_for_coupling_once}
Let $F \in \mathcal{F}_1$.
Fix a probability space and $(D_1,\ldots,D_n,D)\sim F$.
For every $M\geq 1$, there exists $F^{(M)} \in \mathcal{F}_1$ with $|\supp(F^{(M)})|\leq M$ and a deterministic map $\psi_M\colon[0,1]\to[0,1]$ such that, under the pointwise coupling
\[
D_i^{(M)}:=\psi_M\left(D_i\right),\qquad D^{(M)}:=\psi_M\left(D\right),
\]
and $F^{(M)}$ is the distribution of $\psi_M(D)$, we have $D_i^{(M)}\to D_i$ and $D^{(M)}\to D$ almost surely as $M\to\infty$, and $x^*_{F^{(M)}}\to x^*_{F}$.
\end{lemma}
The proof is in Appendix~\ref{prf:mesh_for_coupling_once}.
Fix $F\in\mathcal{F}_1$ as in Lemma~\ref{lem:discretization}. For each $m\geq 1$, apply Lemma~\ref{lem:mesh_for_coupling_once} with $M=m$ and let $F_{d,m}:=F^{(m)} \in \mathcal{F}_1$ denote the corresponding discretization of $F$. Applying Lemma~\ref{lem:coupling_regret_convergence} with $F_m=F_{d,m}$ yields, as $m\to\infty$,
\[
    \RegretCVaR{F_{d,m}}{\pi^{\mcc}}
    \;\to\;
    \RegretCVaR{F}{\pi^{\mcc}},
    \qquad
    \mathrm{PS}_{\delta,n}(\pi^{\mcc},F_{d,m})
    \;\to\;
    \mathrm{PS}_{\delta,n}(\pi^{\mcc},F).
\]
Finally, by Lemma~\ref{lem:mesh_for_coupling_once}, each $F_{d,m}$ belongs to $\mathcal{F}_1$, is discrete, and has support of size at most $m$.

\subsection{Proof of Lemma~\ref{lem:coupling_regret_convergence}} \label{prf:coupling_regret_convergence}
Fix $m$ and write $F':=F_m$, $D_i':=D_i^{(m)}$, $D':=D^{(m)}$.
\paragraph{Convergence of data-driven decisions.}
Let $\boldsymbol{\alpha}=\boldsymbol{\alpha}(W)$ denote the exogenous coefficient vector, and write
\[
x ~~=~~ \sum_{k=1}^n \alpha_k D_{(k)},
\qquad
x' ~~=~~ \sum_{k=1}^n \alpha_k D'_{(k)},
\]
where $D_{(k)}$ and $D'_{(k)}$ are the tie-broken order statistics of $(D_1,\ldots,D_n)$ and $(D_1',\ldots,D_n')$.
Since $D_i'\to D_i$ almost surely for each $i$, we have
\[
\max_{1\leq k\leq n}|D'_{(k)}-D_{(k)}|
~\leq~
\max_{1\leq i\leq n}|D_i'-D_i|
~\to~
0
\quad\text{almost surely},
\]
and therefore $x'\to x$ almost surely.

\paragraph{Convergence of data-driven CVaR.}
Let $X' := c(x',D')$ and $X := c(x,D)$.
Using the Lipschitz bound $|c(x,d)-c(x',d')|\leq (h+b)(|x-x'|+|d-d'|)$ yields
$|X'-X|\to 0$ almost surely.
Since $F,F'\in\mathcal{F}_1$, each demand lies in $[0,1]$ almost surely, so $x$ and $x'$ are convex combinations in $[0,1]$ and $|X'|,|X|\leq 2(h+b)$.
Therefore bounded convergence gives $\mathbb{E}|X'-X|\to 0$, and Lemma~\ref{lem:bound_cvar_diff} implies
\[
    |\cvar_\delta\left(X'\right)-\cvar_\delta\left(X\right)|
    \;\leq\; \frac{1}{1-\delta}\,\mathbb{E}|X'-X|
    \;\to\; 0.
\]

\paragraph{Convergence of the oracle term.} Define $Y' := c(x^*_{F'},D')$ and $Y := c(x^*_F,D)$.
Assumptions~\ref{itm:offline_covergence_coupling} and~\ref{itm:demand_convergence_coupling} give $Y'\to Y$ almost surely.
Again, Lipschitz continuity and $|Y'|,|Y|\leq 2(h+b)$ (since $D',D\in[0,1]$ and $x^*_{F'},x^*_F\in[0,1]$) yield $\mathbb{E}|Y'-Y|\to 0$ by bounded convergence, and hence
$\opt_\delta(F')=\cvar_\delta(Y')\to \cvar_\delta(Y)=\opt_\delta(F)$.
Finally, regret is the difference of the two components, so
$\RegretCVaR{F_m}{\pi^{\mcc}} - \RegretCVaR{F}{\pi^{\mcc}}\to 0$ as $m\to\infty$.

The same coupling yields convergence of expected post-sample regret. Write $\Psi$ and $\Psi'$ for the post-sample outcome-only $\cvar_\delta$ in~\eqref{eqn:conditional_psi} under $F$ and $F'$, evaluated at the coupled orders $x$ and $x'$. For every fixed order $z\in[0,1]$, Lemma~\ref{lem:bound_cvar_diff} applied to the coupled out-of-sample demands gives
\[
\bigl|\cvar_\delta[c(z,D')]-\cvar_\delta[c(z,D)]\bigr|
\;\leq\;
\frac{h+b}{1-\delta}\,\mathbb{E}|D'-D|,
\]
which does not depend on $z$ and tends to $0$. For a common demand distribution, the same lemma and the Lipschitz bound on $c$ give $\bigl|\cvar_\delta[c(x',D)]-\cvar_\delta[c(x,D)]\bigr|\leq (h+b)/(1-\delta)\,|x'-x|$. Therefore $\mathbb{E}|\Psi'-\Psi|\to 0$. Since the oracle term converges,
$\mathrm{PS}_{\delta,n}(\pi^{\mcc},F_m)-\mathrm{PS}_{\delta,n}(\pi^{\mcc},F)\to 0$.

\subsection{Proof of Lemma~\ref{lem:mesh_for_coupling_once}} \label{prf:mesh_for_coupling_once}
We begin by stating the following Lemma:
\begin{lemma}\label{lem:mesh_construction}
    Let $F \in \mathcal{F}_1$, and let $M\geq 1$ be an integer.
    There exists a distribution $F^{(M)} \in \mathcal{F}_1$ with finite support of size at most $M$, together with a deterministic bin map $\phi_M\colon[0,1]\to[0,1]$ such that $F^{(M)}$ is the distribution of $\phi_M(D)$ for $D\sim F$, and
    \[
        |\phi_M\left(d\right)-d|\leq 1/M \quad\text{for every }d\in[0,1],
        \qquad
        |x^*_{F^{(M)}}-x^*_{F}|\leq 1/M \ .
    \]
    \end{lemma}
    The proof is in Appendix~\ref{prf:mesh_construction}. We apply Lemma~\ref{lem:mesh_construction} to $F$ and take $\psi_M:=\phi_M$.
By Lemma~\ref{lem:mesh_construction},
$|D_i^{(M)}-D_i|\leq 1/M$ and $|D^{(M)}-D|\leq 1/M$ under the pointwise coupling
$D_i^{(M)}:=\psi_M(D_i)$, $D^{(M)}:=\psi_M(D)$, hence $D_i^{(M)}\to D_i$ and $D^{(M)}\to D$ almost surely.
The same lemma gives $|x^*_{F^{(M)}}-x^*_{F}|\leq 1/M\to 0$.

\subsection{Proof of Lemma~\ref{lem:mesh_construction} } \label{prf:mesh_construction}

Fix an integer $M \geq 1$.
Recall that
\[
    x^*_{F} \;=\; F^{-1}\!\left(\qstar\right) \;=\; \min\left\{x \geq 0 : F\left(x\right) \geq \qstar\right\}
\]
is the oracle's optimal decision under the profit-CVaR formulation~\eqref{eqn:gotoh_takano_profit_apx}.
Since $x^*_{F} = F^{-1}(\qstar)$, we have $F(x^*_{F})\geq\qstar$ and $F(x)<\qstar$ for every $x<x^*_{F}$.

Partition $[0,1]$ uniformly into $M$ subintervals of width $1/M$:
\[
    B_j := \left[\left(j{-}1\right)\frac{1}{M},\; j\frac{1}{M}\right),
    \quad j=1,\dots,M{-}1,
    \qquad
    B_M := \left[\left(M{-}1\right)\frac{1}{M},\; 1\right].
\]
Write $r_j:=(j{-}1)/M$ for the left endpoint of $B_j$, and define the \emph{bin map} $\phi_M\colon [0,1]\to [0,1]$ by $\phi_M(d):=r_j$ whenever $d\in B_j$.
The \emph{discretized distribution} $F^{(M)}$ is the distribution of $\phi_M(D)$ when $D\sim F$.
Since $\phi_M$ maps $[0,1]$ into $\{0,1/M,\ldots,(M{-}1)/M\}$, we have $F^{(M)}\in\mathcal{F}_1$ and $|\supp(F^{(M)})|\leq M$.

\paragraph{Bin approximation.}
For every $d\in[0,1]$, $d$ and $\phi_M(d)$ lie in the same bin $B_j$, which has width at most $1/M$, so $|\phi_M(d)-d|\leq 1/M$.

\paragraph{Oracle approximation.}
Let $j^*$ be the unique index with $x^*_{F}\in B_{j^*}$, and write $\ell_{j^*}:=r_{j^*}=(j^*{-}1)/M$.
Then $|\ell_{j^*}-x^*_{F}|\leq 1/M$.
Since $B_1\cup\cdots\cup B_{j^*}$ covers $[0,x^*_{F}]$, every $d\leq x^*_{F}$ satisfies $\phi_M(d)\leq r_{j^*}$, hence
\[
    F^{(M)}\left(r_{j^*}\right)
    = P\left(\phi_M\left(D\right)\leq r_{j^*}\right)
    \;\geq\; P\left(D\leq x^*_{F}\right)
    = F\left(x^*_{F}\right)
    \;\geq\; \qstar.
\]
If $j^*\geq 2$, then $B_1\cup\cdots\cup B_{j^*-1}$ covers $[0,\ell_{j^*})$, so
\[
    F^{(M)}\left(r_{j^*-1}\right)
    = P\left(\phi_M\left(D\right)\leq r_{j^*-1}\right)
    = P\left(D<\ell_{j^*}\right)
    = F\left(\ell_{j^*}{}^-\right).
\]
When $\ell_{j^*}<x^*_{F}$, we have $F(\ell_{j^*}{}^-)\leq F(\ell_{j^*})<\qstar$, so $F^{(M)}(r_{j^*-1})<\qstar$ and $x^*_{F^{(M)}}=r_{j^*}$.
When $\ell_{j^*}=x^*_{F}$, we have $F(\ell_{j^*}{}^-)=F(x^*_{F}{}^-)\leq\qstar$; if $F(x^*_{F}{}^-)<\qstar$, then again $F^{(M)}(r_{j^*-1})<\qstar$ and $x^*_{F^{(M)}}=r_{j^*}$, while if $F(x^*_{F}{}^-)=\qstar$, then $F^{(M)}(r_{j^*-1})=\qstar$ and $x^*_{F^{(M)}}=r_{j^*-1}$.
If $j^*=1$, then $x^*_{F^{(M)}}=r_1=0$.
In every case, $x^*_{F^{(M)}}$ equals $r_{j^*}$ or $r_{j^*-1}$, so
\[
    |x^*_{F^{(M)}}-x^*_{F}| \;\leq\; \frac{1}{M}.
\]

\subsection{Proof of Lemma~\ref{lem:data_driven_cvar_convex}}\label{prf:data_driven_cvar_convex}

Fix $\ell\in\{2,\dots,m-1\}$. Let $F_d=\sum_{i=1}^m p_i\,\delta_{a_i}$ with strictly ordered support $0\leq a_1<\cdots<a_m\leq 1$, positive masses $p_i>0$, $a_{-\ell}:=(a_1,\ldots,a_{\ell-1},a_{\ell+1},\ldots,a_m)$, and $p=(p_1,\ldots,p_m)$.
Consider the perturbation family
\[
F_t^{(\ell)} \;:=\; \sum_{i\neq \ell} p_i\,\delta_{a_i} \;+\; p_\ell\,\delta_{t},
\qquad
t\in[a_{\ell-1},a_{\ell+1}] \ .
\]
As in the lemma statement, couple the family $\{F_t^{(\ell)}\}$ using fixed atom labels whose probabilities do not depend on $t$, sending label $\ell$ to the moving support point $t$ and every other label $i$ to $a_i$, and write $D_1(t),\ldots,D_n(t),D(t)$ for the resulting demands (independent of any policy randomization). Set
\[
    C_t \;:=\; c\bigl(x^{\mcc}(D_1^n(t)),D(t)\bigr).
\]
We show that $t\mapsto \cvar_\delta(C_t)$ is convex on $[a_{\ell-1},a_{\ell+1}]$.

Since the label probabilities are fixed as $t$ varies on $[a_{\ell-1},a_{\ell+1}]$, the joint distribution of the labels does not depend on $t$. Fix an outcome $\omega$, i.e., fix the atom labels and the policy randomization. For $t\in(a_{\ell-1},a_{\ell+1})$, because $a_1<\cdots<a_m$ and $t$ stays strictly between $a_{\ell-1}$ and $a_{\ell+1}$, there are no crossings between the moving value $t$ and the fixed support points $\{a_k:k\neq \ell\}$. Therefore, the set of order-statistic ranks occupied by sample coordinates labeled $\ell$ is constant. Its total coefficient $\gamma$ is consequently constant, while all remaining order-statistic values are fixed, and hence
\[
    x^{\mcc}\bigl(D_1^n(t)\bigr)=x_0+\gamma t
\]
for some constants $x_0\in\mathbb{R}$ and $\gamma\in[0,1]$ (depending on $\omega$), where $\gamma=\sum_{k:\,D_{(k)}(t)=t}\alpha_k$ is the total weight on order-statistic ranks occupied by the value $t$.

Recall that the newsvendor cost can be written as the maximum of two affine functions:
\[
    c\left(x,d\right)=\max\{-bx,\;hx-\left(h+b\right)d\}.
\]
If the out-of-sample label is not $\ell$, then $D(t)$ equals some fixed $a_i$ and is constant in $t$. Hence both
\[
    f_1\left(t\right)=-b\,\left(x_0+\gamma t\right),
    \qquad
    f_2\left(t\right)=h\,\left(x_0+\gamma t\right)-\left(h+b\right)D(t)
\]
are affine in $t$, so $t\mapsto C_t(\omega)$ is convex as the pointwise maximum of two affine functions.
If the out-of-sample label is $\ell$, then $D(t)=t$, so
\[
    f_1\left(t\right)=-b\,\left(x_0+\gamma t\right),
    \qquad
    f_2\left(t\right)=h\,\left(x_0+\gamma t\right)-\left(h+b\right)t
\]
are again affine in $t$, and therefore $t\mapsto C_t(\omega)=\max\{f_1(t),f_2(t)\}$ is convex.
Thus for almost every outcome $\omega$, the map $t\mapsto C_t(\omega)$ is convex on $(a_{\ell-1},a_{\ell+1})$.

We now use the Rockafellar--Uryasev representation:
\[
    \cvar_\delta\left(C_t\right)
    =
    \inf_{u\in\mathbb{R}}
    \left\{
        u + \frac{1}{1-\delta}\,\mathbb{E}\left[\left(C_t-u\right)^+\right]
    \right\}
    ~=:~ \inf_{u\in\mathbb{R}}\phi\left(t,u\right).
\]
We show that $\phi$ is \emph{jointly} convex in $(t,u)$ on $(a_{\ell-1},a_{\ell+1})\times\mathbb{R}$.
Fix an outcome $\omega$ for which $t\mapsto C_t(\omega)$ is convex; as shown above, this holds for almost every $\omega$.
Then $(t,u)\mapsto C_t(\omega)-u$ is jointly convex, being convex in $t$ and affine in $u$, and composing with the convex nondecreasing map $x\mapsto x^+$ preserves joint convexity, so $(t,u)\mapsto (C_t(\omega)-u)^+$ is jointly convex.
Taking expectations (a nonnegative average of jointly convex maps) and adding the affine term $u$ shows that $\phi$ is jointly convex in $(t,u)$.
We now invoke the standard partial-minimization theorem for convex functions: if $\phi(t,u)$ is jointly convex, then $t\mapsto \inf_u \phi(t,u)$ is convex on its effective domain. Since $\cvar_\delta(C_t)=\inf_{u\in\mathbb{R}}\phi(t,u)$ is finite because $C_t$ is bounded, $t\mapsto\cvar_\delta(C_t)$ is convex on $(a_{\ell-1},a_{\ell+1})$.
It remains to extend convexity to the closed interval $[a_{\ell-1},a_{\ell+1}]$. Let $t_0\in\{a_{\ell-1},a_{\ell+1}\}$. Under the fixed-label coupling, the moving atom merges with a neighboring atom as $t\to t_0$; tie-invariance of the policy and continuity of $c$ imply that $C_t\to C_{t_0}$ almost surely. The losses are uniformly bounded, so dominated convergence and Lemma~\ref{lem:bound_cvar_diff} give $\cvar_\delta(C_t)\to \cvar_\delta(C_{t_0})$. Thus, $t\mapsto\cvar_\delta(C_t)$ is continuous on the closed interval, and its convexity on the open interval extends to the endpoints.

The same conclusion holds for the expected post-sample term. Conditional on the sample labels and the realized coefficient vector, $\Psi_n^{\pi^{\mcc}}$ in equation~\eqref{eqn:conditional_psi} is the outcome-only $\cvar_\delta$ of an independent out-of-sample demand. For almost every realization of these labels and coefficients, the loss is convex in $t$ for each out-of-sample label, so the same partial-minimization step under that conditional distribution yields convexity of $t\mapsto\Psi_n^{\pi^{\mcc}}$ on $(a_{\ell-1},a_{\ell+1})$. Taking expectations over the sample labels and coefficients preserves convexity. At an endpoint $t_0$, the same almost-sure convergence of the loss and Lemma~\ref{lem:bound_cvar_diff}, applied conditionally on the sample labels and coefficients, give $\Psi_n^{\pi^{\mcc}}(t)\to\Psi_n^{\pi^{\mcc}}(t_0)$ almost surely. The conditional tails are uniformly bounded, so $\mathbb{E}[\Psi_n^{\pi^{\mcc}}(t)]\to\mathbb{E}[\Psi_n^{\pi^{\mcc}}(t_0)]$, and convexity extends to the closed interval.

\subsection{Proof of Lemma~\ref{lem:off_linear}}\label{prf:off_linear}

Fix $\ell\in\{2,\dots,m-1\}$. Let $F_d=\sum_{i=1}^m p_i\,\delta_{a_i}$ with strictly ordered support $0\leq a_1<\cdots<a_m\leq 1$, positive masses $p_i>0$, $a_{-\ell}:=(a_1,\ldots,a_{\ell-1},a_{\ell+1},\ldots,a_m)$, and $p=(p_1,\ldots,p_m)$.
Consider the perturbation family
\[
F_t^{(\ell)} \;:=\; \sum_{i\neq \ell} p_i\,\delta_{a_i} \;+\; p_\ell\,\delta_{t},
\qquad
t\in[a_{\ell-1},a_{\ell+1}] \ .
\]
Let $t\in[a_{\ell-1},a_{\ell+1}]$. Set
\[
g\left(t\right):=\opt_\delta\left(F_t^{(\ell)}\right),\qquad x_t^\star:=\left(F_t^{(\ell)}\right)^{-1}\left(q^\star\right),
\]
so $g(t)=\cvar_\delta(c(x_t^\star,D))$ with $D\sim F_t^{(\ell)}$.

We prove that $g$ is affine in $t$ by four steps:
(i) identify the form of $x_t^\star$,
(ii) write finite-support CVaR as an affine function of sorted losses,
(iii) show the CVaR coefficients are constant in $t$, and
(iv) show the selected losses are affine in $t$.

\textbf{Step 1 (form of $x_t^\star$).}
Only one atom location moves (from $a_\ell$ to $t$). For $t\in(a_{\ell-1},a_{\ell+1})$, support order and atom probabilities are unchanged, so the quantile atom containing $q^\star$ is fixed in index. Consequently, on this open interval, either:
\[
x_t^\star=a_k \text{ for some fixed } k\neq \ell,
\]
or
\[
x_t^\star=t.
\]

\textbf{Step 2 (finite-support CVaR representation).}
Fix $t$ and a decision value $x$. Let support atoms of $D$ be $d_1(t),\dots,d_m(t)$, and define losses
\[
L_i\left(x,t\right):=c\left(x,d_i\left(t\right)\right).
\]
Sort them:
\[
L_{(1)}\left(x,t\right)\geq\cdots\geq L_{(m)}\left(x,t\right),
\]
with corresponding probabilities $p_{(1)},\dots,p_{(m)}$. Let
\[
r=\min\Big\{j:\sum_{i=1}^j p_{(i)}\geq 1-\delta\Big\},\qquad
\theta=\left(1-\delta\right)-\sum_{i=1}^{r-1}p_{(i)}.
\]
Then
\[
\cvar_\delta\left(c\left(x,D\right)\right)
=\frac{1}{1-\delta}\left(\sum_{i=1}^{r-1}p_{(i)}L_{(i)}\left(x,t\right)+\theta L_{(r)}\left(x,t\right)\right).
\]
So CVaR is affine in the selected sorted losses.

\textbf{Step 3 (coefficients do not change with $t$).}
Uniformly over every decision $x$, the map $d\mapsto c(x,d)$ is weakly decreasing in $d$ (all demand atoms above $x$ tie at loss $-bx$). Thus higher losses correspond to smaller-demand atoms. Since demand-atom order and probabilities are unchanged with $t$, the tail is selected by demand order and the CVaR value is invariant to tie-breaking among tied losses; hence $r$, $\theta$, and the coefficients in the CVaR combination are constant on $(a_{\ell-1},a_{\ell+1})$.

\textbf{Step 4 (selected losses are affine in $t$).}
If $x_t^\star=a_k$, then
\[
g\left(t\right)=\cvar_\delta\left(c\left(a_k,D\right)\right),
\]
a fixed-coefficient combination of $\{c(a_k,a_i)\}_{i\neq\ell}$ and $c(a_k,t)$. For fixed $a_k$,
\[
c\left(a_k,t\right)=b\left(t-a_k\right)^+ + h\left(a_k-t\right)^+ - bt,
\]
which is affine on $[a_{\ell-1},a_{\ell+1}]$: if $k<\ell$, then $a_k\leq t$ throughout and $c(a_k,t)=-ba_k$; if $k>\ell$, then $a_k\geq t$ throughout and $c(a_k,t)=ha_k-(h+b)t$. Hence $g$ is affine.

If $x_t^\star=t$, then
\[
g\left(t\right)=\cvar_\delta\left(c\left(t,D\right)\right),
\]
a fixed-coefficient combination of $\{c(t,a_i)\}_{i\neq\ell}$ and $c(t,t)$. For each $i\neq\ell$, $t\mapsto c(t,a_i)$ is affine on $[a_{\ell-1},a_{\ell+1}]$, and $c(t,t)=-bt$ is affine. Hence $g$ is affine on $(a_{\ell-1},a_{\ell+1})$.

It remains to extend affinity to the endpoints. Let $t_0\in\{a_{\ell-1},a_{\ell+1}\}$, and couple $D_t\sim F_t^{(\ell)}$ and $D_{t_0}\sim F_{t_0}^{(\ell)}$ using the same atom label. Then $\mathbb{E}|D_t-D_{t_0}|=p_\ell|t-t_0|$. For every fixed $x$, Lemma~\ref{lem:bound_cvar_diff} and the Lipschitz continuity of $c$ in demand give
\[
\left|\cvar_\delta[c(x,D_t)]-\cvar_\delta[c(x,D_{t_0})]\right|
\leq
\frac{h+b}{1-\delta}\,p_\ell|t-t_0|.
\]
This bound is uniform in $x$, so taking minima over $x$ yields $g(t)\to g(t_0)$ as $t\to t_0$. Thus $g$ is continuous on $[a_{\ell-1},a_{\ell+1}]$, and its affinity on the open interval extends to the endpoints.
In both cases, $t\mapsto \opt_\delta(F_t^{(\ell)})$ is affine on $[a_{\ell-1},a_{\ell+1}]$.

\section{Additional Proofs from Section~\ref{sec:order_stats}}

\subsection{Proof of Claim~\ref{clm:shifting}} \label{prf:shifting}
We define $F_0$ as a shifting of $F_2$ by $a_1$:
\[F_0\left(x\right) ~~=~~ F_2\left(x+a_1\right) \ . \]
Consequently, we have 
\[F_0\left(0\right) ~~=~~ F_2\left(a_1\right) ~~>~~ 0 \ . \]
Moreover, for $x < 0$ we have
\[F_0\left(x\right) ~~=~~ F_2\left(x+a_1\right) ~~=~~ 0 \ . \]
Here, the second equality holds since $x+a_1 < a_1$, and $a_1$ is the smallest supported value of $F_2$. 
If $F_2$ has support $\{a_1,a_2\}\subseteq [0,1]$, then $F_0$ has support $\{0,a_2-a_1\}\subseteq [0,1]$, and therefore $F_0 \in \mathcal{F}_1$.
It remains to show that the regret stays the same. Recalling that for each acquired quantity $x$ and demand realization $D$ the associated cost is
\[c\left(x,D\right) ~~=~~ \begin{cases}
    hx - \left(h+b\right)D \quad & D \leq x \\
    -bx \quad & D > x
\end{cases} \ , \]
we observe that 
\[c\left(x+a_1,D+a_1\right) ~~=~~ \begin{cases}
    h\left(x+a_1\right) - \left(h+b\right)\left(D+a_1\right) \quad & D \leq x \\
    -b\left(x+a_1\right) \quad & D > x
\end{cases} ~~=~~ c\left(x,D\right) -ba_1  \ . \]
Let $D_1,\ldots,D_n,D$ be independent with $D_i\sim F_2$ and $D\sim F_2$, and define
$\tilde D_i:=D_i-a_1$ and $\tilde D:=D-a_1$.
Then $\tilde D_1,\ldots,\tilde D_n,\tilde D$ are independent with $\tilde D_i\sim F_0$ and $\tilde D\sim F_0$.

Translation preserves order, and $x^\mcc(D_1^n)=\sum_{k=1}^n \alpha_k D_{(k)}$, so
\[
x^\mcc\left(\tilde D_1^n\right) ~~=~~ x^\mcc\left(D_1^n\right)-a_1.
\]
Therefore, almost surely,
\[
c\!\left(x^\mcc\left(D_1^n\right),D\right)+ba_1
~~=~~ c\!\left(x^\mcc\left(\tilde D_1^n\right),\tilde D\right).
\]
By translation invariance of $\cvar_\delta$,
\begin{eqnarray*}
    \cvar_\delta[c\left(x^\mcc(D_1^n),D\right)]
    &=& \cvar_\delta[c\left(x^\mcc(\tilde D_1^n),\tilde D\right)] - ba_1 \ ,
\end{eqnarray*}
where the first $\cvar$ is under $(D_1^n,D)$ with $D_i,D\sim F_2$,
and the second under $(\tilde D_1^n,\tilde D)$ with $\tilde D_i,\tilde D\sim F_0$.

Similarly, $x^*_{F_2}-a_1=x^*_{F_0}$ and $c(x^*_{F_2},D)+ba_1=c(x^*_{F_0},\tilde D)$ almost surely, so
\begin{eqnarray*}
    \opt_\delta\left(F_2\right) &=& \opt_\delta\left(F_0\right) - ba_1 \ .
\end{eqnarray*}
Therefore
\begin{eqnarray*}
    \RegretCVaR{F_2}{\pi^{\mcc}}
    &=& \cvar_\delta[c\left(x^\mcc(D_1^n),D\right)] - \opt_\delta\left(F_2\right) \\
    &=& \left(\cvar_\delta[c\left(x^\mcc(\tilde D_1^n),\tilde D\right)] - ba_1\right) - \left(\opt_\delta\left(F_0\right) - ba_1\right) \\
    &=& \RegretCVaR{F_0}{\pi^{\mcc}} \ .
\end{eqnarray*}
The same translation applies to expected post-sample regret. Conditional on the sample and the realized coefficient vector, $\Psi_n^{\pi^{\mcc}}$ is the outcome-only $\cvar_\delta$ of $c(x^{\mcc},D)$. The identity $c(x+a_1,D+a_1)=c(x,D)-ba_1$ and translation invariance of $\cvar_\delta$ shift $\Psi_n^{\pi^{\mcc}}$ by $-ba_1$ almost surely, so its expectation shifts by $-ba_1$. Subtracting the oracle, which shifts by the same amount, leaves $\mathrm{PS}_{\delta,n}$ unchanged.

\subsection{Proof of Claim~\ref{clm:scaling}} \label{prf:scaling}

Fix $T>0$ and $F_T\in\mathcal{F}_T$, and define
\[
F_1\left(y\right) ~:=~ F_T\left(Ty\right), \qquad y\in\mathbb{R}.
\]
Because $\support{F_T}\subseteq[0,T]$, we have $\support{F_1}\subseteq[0,1]$, so $F_1\in\mathcal{F}_1$.

Recalling that for each acquired quantity $x$ and demand realization $D$ the associated cost is
\[c\left(x,D\right) ~~=~~ \begin{cases}
    hx - \left(h+b\right)D \quad & D \leq x \\
    -bx \quad & D > x
\end{cases} \ , \]
we have $c(x/T,D/T)=c(x,D)/T$ for every $x,D\geq 0$.

Let $D_1,\ldots,D_n,D$ be independent with $D_i\sim F_T$ and $D\sim F_T$, and define
$\tilde D_i:=D_i/T$ and $\tilde D:=D/T$.
Then $\tilde D_1,\ldots,\tilde D_n,\tilde D$ are independent with $\tilde D_i\sim F_1$ and $\tilde D\sim F_1$.

Since $\pi^{\mcc}$ prescribes $x^\mcc(D_1^n)=\sum_{k=1}^n \alpha_k D_{(k)}$ and scaling preserves order,
\[
x^\mcc\left(\tilde D_1^n\right) ~~=~~ \frac{1}{T}\,x^\mcc\left(D_1^n\right).
\]
Therefore, almost surely,
\[
c\!\left(x^\mcc\left(D_1^n\right),D\right) ~~=~~ T\cdot c\!\left(x^\mcc\left(\tilde D_1^n\right),\tilde D\right).
\]
By positive homogeneity of $\cvar_\delta$,
\begin{eqnarray*}
    \cvar_\delta[c\left(x^\mcc,D\right)]
    &=& T\cdot \cvar_\delta[c\left(x^\mcc,\tilde D\right)] \ ,
\end{eqnarray*}
where the first $\cvar$ is computed under $(D_1^n,D)$ with $D_i,D\sim F_T$ as above,
and the second under $(\tilde D_1^n,\tilde D)$ with $\tilde D_i,\tilde D\sim F_1$.

Similarly, $x^*_{F_T}/T=x^*_{F_1}$ because $x^*_F$ depends on $F$ only through quantiles.
With the same coupling, $c(x^*_{F_T},D)=T\cdot c(x^*_{F_1},\tilde D)$ almost surely, hence
\begin{eqnarray*}
    \opt_\delta\left(F_T\right) &=& T\cdot \opt_\delta\left(F_1\right) \ .
\end{eqnarray*}
Therefore
\[
\RegretCVaR{F_1}{\pi^{\mcc}}
~~=~~ \cvar_\delta[c\left(x^\mcc,\tilde D\right)] - \opt_\delta\left(F_1\right)
~~=~~ \frac{\cvar_\delta[c\left(x^\mcc,D\right)] - \opt_\delta\left(F_T\right)}{T}
~~=~~ \frac{\RegretCVaR{F_T}{\pi^{\mcc}}}{T} \ .
\]
The same scaling applies to expected post-sample regret. Conditional on the realized order, $\Psi_n^{\pi^{\mcc}}$ is the outcome-only $\cvar_\delta$ in~\eqref{eqn:conditional_psi}. The loss identity and the substitution $t=Ts$ multiply that infimum by $T$, so $\mathbb{E}[\Psi_n^{\pi^{\mcc}}]$ scales by $T$. Subtracting $\opt_\delta(F_T)=T\cdot\opt_\delta(F_1)$ gives
\[
\mathrm{PS}_{\delta,n}(\pi^{\mcc},F_1)
~~=~~
\frac{\mathrm{PS}_{\delta,n}(\pi^{\mcc},F_T)}{T}.
\]
Conversely, given any $F_1\in\mathcal{F}_1$, define $F_T(y):=F_1(y/T)$. Then $F_T\in\mathcal{F}_T$, and the same scaling calculation gives both equalities in the reverse direction.

\subsection{Proof of Claim~\ref{clm:optimality_deterministic_mcc}}\label{prf:optimality_deterministic_mcc}
Fix $F \in \mathcal{F}_1$, let $D_1^n$ be an i.i.d.\ sample from $F$, and let $D \sim F$ be independent of the sample and of the policy randomization.
Let $L_\pi := c(x^\pi, D)$ and $L_{\bar{\pi}} := c(\bar{x}, D)$ denote the induced loss random variables under the respective policies. 
The regret inequality $\RegretCVaR{F}{\pi} \geq \RegretCVaR{F}{\bar{\pi}}$ is equivalent to $\cvar_\delta[L_\pi] \geq \cvar_\delta[L_{\bar{\pi}}]$, since $\opt_\delta(F)$ is policy-independent. We therefore prove $\cvar_\delta[L_\pi] \geq \cvar_\delta[L_{\bar{\pi}}]$.

\paragraph{Step 1: Determinism of $\bar{\pi}$.}
By construction, $\bar{\alpha}_k = \mathbb{E}[\alpha_k]$ for each $k$, so the weights of $\bar{\pi}$ are fixed ex ante and do not depend on the sample or on any randomness. Hence $\bar{\pi}$ is deterministic.

\paragraph{Step 2: Conditional mean of the decision.}
Since the randomization of the weight vector $(\alpha_1,\ldots,\alpha_n)$ is done ex-ante of observing the sample, it is independent of the sample $D_1^n = (D_1,\ldots,D_n)$, so $\mathbb{E}[\alpha_k \mid D_1,\ldots,D_n] = \mathbb{E}[\alpha_k] = \bar{\alpha}_k$. The policy's prescribed quantity is $x^\pi = \sum_{k=1}^n \alpha_k\, D_{(k)}$. Thus
\[
\mathbb{E}[x^\pi \mid D_1,\ldots,D_n] \;=\; \sum_{k=1}^n \mathbb{E}[\alpha_k \mid D_1,\ldots,D_n]\, D_{(k)} \;=\; \sum_{k=1}^n \bar{\alpha}_k\, D_{(k)} \;=\; \bar{x}\,.
\]
The future demand $D$ is independent of the sample and of the weight vector, so $\mathbb{E}[x^\pi \mid D_1^n, D] = \bar{x}$ as well.

\paragraph{Step 3: Conditional Jensen for the cost.}
The newsvendor loss $c(x,D)$ is convex in $x$ for each fixed $D$: using $c(x,D)=-bx+(h+b)(x-D)^+$, it can be written as $\max\{-bx,\,hx-(h+b)D\}$. Conditional on $(D_1^n, D)$, the quantity $\bar{x} = \bar{x}(D_1^n)$ is deterministic while $x^\pi$ is random with $\mathbb{E}[x^\pi \mid D_1^n, D] = \bar{x}$. By Jensen's inequality,
\[
\mathbb{E}\left[c\left(x^\pi, D\right) \,\middle|\, D_1^n, D\right] \;\geq\; c\!\left(\mathbb{E}[x^\pi \mid D_1^n, D],\, D\right) \;=\; c\left(\bar{x}, D\right) \;=\; L_{\bar{\pi}}\,.
\]
Thus, almost surely,
\[
\mathbb{E}[L_\pi \mid D_1^n, D] \;\geq\; L_{\bar{\pi}}\,.
\]

\paragraph{Step 4: CVaR comparison.}
From Step~3, $\mathbb{E}[L_\pi \mid D_1^n, D] \geq L_{\bar{\pi}}$ almost surely. For each $t \in \mathbb{R}$, the map $\phi(u) := (u-t)^+$ is convex, so conditional Jensen yields
\[
\mathbb{E}\left[\left(L_\pi - t\right)^+ \mid D_1^n, D\right]
\;\geq\;
\left(\mathbb{E}[L_\pi \mid D_1^n, D] - t\right)^+ \,.
\]
Since $\phi$ is also nondecreasing and $\mathbb{E}[L_\pi \mid D_1^n, D] \geq L_{\bar{\pi}}$, we obtain
\[
\mathbb{E}\left[\left(L_\pi - t\right)^+ \mid D_1^n, D\right] \;\geq\; \left(L_{\bar{\pi}} - t\right)^+ \,.
\]
Taking expectations and using the Rockafellar--Uryasev representation $\cvar_\delta[X] = \inf_{t \in \mathbb{R}} \bigl\{ t + (1-\delta)^{-1} \mathbb{E}[(X-t)^+] \bigr\}$ gives
\[
t + \frac{1}{1-\delta}\,\mathbb{E}\left[\left(L_\pi - t\right)^+\right] \;\geq\; t + \frac{1}{1-\delta}\,\mathbb{E}\left[\left(L_{\bar{\pi}} - t\right)^+\right] \qquad \text{for all } t \in \mathbb{R}\,,
\]
and taking the infimum over $t$ yields $\cvar_\delta[L_\pi] \geq \cvar_\delta[L_{\bar{\pi}}]$.

Hence $\cvar_\delta[L_\pi] \geq \cvar_\delta[L_{\bar{\pi}}]$, which yields $\RegretCVaR{F}{\pi} \geq \RegretCVaR{F}{\bar{\pi}}$ for every $F \in \mathcal{F}_1$. The consequence $\sup_{F \in \mathcal{F}_1} \RegretCVaR{F}{\pi} \geq \sup_{F \in \mathcal{F}_1} \RegretCVaR{F}{\bar{\pi}}$ follows immediately.

\subsection{Bernoulli event probabilities and CVaR tail filling}
\label{apx:bernoulli_event_enumeration}
Using the events $O_i,U_i$ and associated costs defined in Section~\ref{sec:computed_policy}, the associated probabilities under $\mathrm{Ber}(\mu)$ are
\[p\left(O_i\right) := \underbrace{\binom{n}{i-1} \cdot \mu^{n-i+1} \cdot \left(1-\mu\right)^{i-1}}_{P(D_{(i-1)} = 0, D_{(i)} = 1)} \cdot \underbrace{\left(1-\mu\right)}_{P(D=0)} ~~=~~ \binom{n}{i-1} \cdot \mu^{n-i+1} \cdot \left(1-\mu\right)^{i} \ . \]
For $i >1$:
\[p\left(U_i\right) := \underbrace{\binom{n}{i-1} \cdot \mu^{i-1} \cdot \left(1-\mu\right)^{n-i+1}}_{P(D_{(n-i+2)} = 1, D_{(n-i+1)} = 0)} \cdot \underbrace{\mu}_{P(D=1)} ~~=~~ \binom{n}{i-1} \cdot \mu^{i} \cdot \left(1-\mu\right)^{n-i+1} \ , \]
and $P(U_1) = (1-\mu)^n$. The resulting data-driven CVaR is
\begin{eqnarray*}
    \cvar_\delta[c\left(x^{\pi}, D\right)] &=& \frac{1}{1-\delta} \cdot \left[\sum_{i=1}^n c\left(O_i\right) \cdot \min \left\{p\left(O_i\right), \left[1-\delta - \sum_{j=1}^{i-1} p\left(O_j\right)\right]^+\right\} \right . \\
    && \left. + \sum_{i=1}^{n+1} c\left(U_i\right) \cdot \min \left\{p\left(U_i\right), \left[1-\delta - \sum_{j=1}^n p\left(O_j\right) - \sum_{j=1}^{i-1} p\left(U_j\right)\right]^+\right\} \right].
\end{eqnarray*}
For fixed $\mu$, the event probabilities and tail-selection weights are independent of $A_i$. Since the event costs are linear in the cumulative tail weights, the data-driven CVaR is linear in $(A_1,\ldots,A_n)$, and hence in $\boldsymbol{\alpha}$.

\subsection{Proof of Claim~\ref{clm:lipschitz_constant}}\label{prf:lipschitz_constant}
Let
\[
R\left(\mu\right):=\cvar_\delta[c\left(x^\pi_\mu,D_\mu\right)]-\opt_\delta\left(\mathrm{Ber}\left(\mu\right)\right),
\]
with the notation from Claim~\ref{clm:lipschitz_constant}. We show that for every
$\mu,\nu\in[0,1]$,
\[
|R\left(\mu\right)-R\left(\nu\right)|\leq \frac{\left(h+b\right)\left(n+2\right)}{1-\delta}|\mu-\nu|.
\]

Couple the two Bernoulli systems by common uniforms. Let
$U_1,\ldots,U_n,U_0$ be i.i.d.\ uniform random variables on $[0,1]$, and set
$D_i^\mu=\mathbbm{1}\{U_i\leq\mu\}$, $D_i^\nu=\mathbbm{1}\{U_i\leq\nu\}$ for
$i=1,\ldots,n$, and similarly $D_\mu=\mathbbm{1}\{U_0\leq\mu\}$,
$D_\nu=\mathbbm{1}\{U_0\leq\nu\}$. Under this coupling, $D_i^\mu$ and
$D_i^\nu$ differ only when $U_i$ falls between $\mu$ and $\nu$, an event of
probability $|\mu-\nu|$. Thus
\[
\mathbb{E}|D_\mu-D_\nu|=|\mu-\nu|,\qquad
\Pr\left(D_i^\mu\neq D_i^\nu\right)=|\mu-\nu| \quad \text{for each } i.
\]
The two sample vectors differ only if at least one coordinate differs, so the
union bound gives
\[
\Pr\{\left(D_1^\mu,\ldots,D_n^\mu\right)\neq\left(D_1^\nu,\ldots,D_n^\nu\right)\}
\leq n|\mu-\nu|.
\]
Because the policy action is a convex combination of Bernoulli order
statistics, $x^\pi_\mu,x^\pi_\nu\in[0,1]$. Moreover, if the two samples are
identical then their order statistics, and hence the two policy actions, are
identical. Therefore
\[
\mathbb{E}|x^\pi_\mu-x^\pi_\nu|
\leq
\Pr\{\left(D_1^\mu,\ldots,D_n^\mu\right)\neq\left(D_1^\nu,\ldots,D_n^\nu\right)\}
\leq n|\mu-\nu|.
\]

The newsvendor cost satisfies
\[
|c\left(x,d\right)-c\left(x',d'\right)|\leq \left(h+b\right)\left(|x-x'|+|d-d'|\right)
\]
for all feasible decisions $x,x'$ and all $d,d'\in[0,1]$. Hence
\[
\mathbb{E}\left|c\left(x^\pi_\mu,D_\mu\right)-c\left(x^\pi_\nu,D_\nu\right)\right|
\leq \left(h+b\right)\left(n+1\right)|\mu-\nu|.
\]
Applying Lemma~\ref{lem:bound_cvar_diff} gives
\[
\left|\cvar_\delta[c\left(x^\pi_\mu,D_\mu\right)]
-\cvar_\delta[c\left(x^\pi_\nu,D_\nu\right)]\right|
\leq
\frac{\left(h+b\right)\left(n+1\right)}{1-\delta}|\mu-\nu|.
\]

It remains to bound the oracle component. For any fixed feasible decision $x$,
the same coupling gives
\[
\mathbb{E}|c\left(x,D_\mu\right)-c\left(x,D_\nu\right)|\leq \left(h+b\right)|\mu-\nu|,
\]
and therefore, by Lemma~\ref{lem:bound_cvar_diff},
\[
\left|\cvar_\delta[c\left(x,D_\mu\right)]-\cvar_\delta[c\left(x,D_\nu\right)]\right|
\leq \frac{h+b}{1-\delta}|\mu-\nu|.
\]
Since this bound is uniform over feasible decisions, taking the infimum over
$x$ preserves the same Lipschitz bound for the value function, so
\[
|\opt_\delta\left(\mathrm{Ber}\left(\mu\right)\right)-\opt_\delta\left(\mathrm{Ber}\left(\nu\right)\right)|
\leq \frac{h+b}{1-\delta}|\mu-\nu|.
\]
Combining the two bounds by the triangle inequality yields
\[
|R\left(\mu\right)-R\left(\nu\right)|
\leq
\frac{\left(h+b\right)\left(n+2\right)}{1-\delta}|\mu-\nu|.
\]

\subsection{Constraint generating linear programming implementation of the minimax problem}\label{apx:implementation}
When we compute the linear programming policy, we avoid loading all grid constraints at once. Fix the grid $\mathcal{G}_m:=\{0,1/m,\ldots,1\}$ and use the following constraint-generating procedure:
\begin{enumerate}
    \item Start with an initial set of $\mu$ values $M \subseteq \mathcal{G}_m$.
    \item Solve the following linear program:
    \begin{equation*}
        \begin{aligned}
            \min_{\alpha \in \mathbb{R}^n, t \in \mathbb{R}} \quad & t \\
            \textrm{s.t.} \quad & \cvar_\delta[c\left(x^{\pi}, D\right)] - \opt_\delta\left(\mathrm{Ber}\left(\mu\right)\right) \leq t \quad & \forall \mu \in M \\
            \quad & \sum_{i=1}^n \alpha_i = 1 \\
            \quad & \alpha_i \geq 0 \quad & \forall i \in [n]
        \end{aligned}
    \end{equation*}

    \item Let $(\alpha^*, t^*)$ be an optimal solution of the above program. Evaluate the regret at every point of $\mathcal{G}_m$ and let $\mu^*$ be a maximizer.
    
    \item If $\RegretCVaRDelta{\mathrm{Ber}(\mu^*)}{x^{\alpha^*}}{\delta} \leq t^* + \epsilon$, stop and return $\alpha^*$. Otherwise, add $\mu^*$ to $M$ and repeat from step 2.
\end{enumerate}
At termination, the returned policy is feasible for the full grid LP up to $\epsilon$.
We implement the constraint-generation procedure above to obtain the coefficient vectors used in the numerical comparisons below.

\subsection{Constraint generating two-stage procedure for combined-risk regret and expected post-sample regret multiobjective optimization}\label{apx:two_stage_procedure}
We minimize worst-case expected post-sample regret subject to a slack of $\epsilon$ above the minimum worst-case combined-risk regret. Formally, if we solved the semi-infinite linear program from equation~\eqref{eqn:semi_infinite_lp}, and used $t^*$ to denote its optimal value, we would wish to solve the following semi-infinite linear program:
\begin{equation} \label{eqn:semi_infinite_lp_stage_2}
    \begin{aligned}
        \min_{\alpha \in \mathbb{R}^n, s \in \mathbb{R}} \quad & s \\
        \textrm{s.t.} \quad & \cvar_\delta[c\left(x^{\alpha}, D\right)] - \opt_\delta\left(\mathrm{Ber}\left(\mu\right)\right) \leq t^* + \epsilon \quad & \forall \mu \in [0,1] \\
        \quad & \mathbb{E}_{D_1^n, x^\alpha}[\Psi_n^\alpha] - \opt_\delta\left(\mathrm{Ber}\left(\mu\right)\right) \leq s \quad & \forall \mu \in [0,1] \\
        \quad & \sum_{i=1}^n \alpha_i = 1 \\
        \quad & \alpha_i \geq 0 \quad & \forall i \in [n]
    \end{aligned}
\end{equation}
To claim that this is indeed a linear program, we must show that $\mathbb{E}_{D_1^n, x^\alpha}[\Psi_n^\alpha]$ is linear in $\alpha$ (the other components were already shown to be affine in the linear program from equation~\eqref{eqn:semi_infinite_lp}). Under a Bernoulli distribution with parameter $\mu$, $\mathbb{E}_{D_1^n, x^\alpha}[\Psi_n^\alpha]$ is computed as follows: 
index sample events by their number of zeros, $j-1$, for $j=1,\ldots,n+1$. On event $j$, the order is $x=x_j=\sum_{i=j}^n\alpha_i$, with $x_{n+1}=0$.
\begin{enumerate}
    \item The relevant events are $[D=0, D=1]$.
    \item The associated costs are $c(x,0) = hx$ and $c(x,1) = -bx$.
    \item The probability of $c(x,0)$ is $\Pr(D=0) = 1-\mu$ and the probability of $c(x,1)$ is $\Pr(D=1) = \mu$.
    \item The vector that defines $\Psi_n^\alpha$ has components $\frac{hx_j \cdot \min\{1-\mu, 1-\delta\} - bx_j [\mu -\delta]^+}{1-\delta}$, where $x_j = \sum_{i=j}^n \alpha_i$. This term is linear in $\alpha$.
    \item 
    Event $j$ has probability $\binom{n}{j-1}(1-\mu)^{j-1}\mu^{n-j+1}$, independent of $\alpha$. We keep these sample events distinct even when their prescribed orders coincide.
    \item The expectation is the inner product of the probabilities vector and the CVaR vector. Therefore it is linear in $\alpha$.
\end{enumerate}
Having established the linearity of the second stage semi-infinite linear program, we move on to the discretization and constraint generation procedure. 
Let $t^*$ now denote the computed first-stage value, and let $\epsilon_1$ be a valid bound on its distance from the semi-infinite optimum $V_n^{\mathrm{det}}$. The discretized version of equation~\eqref{eqn:semi_infinite_lp_stage_2} is:
\begin{equation*} 
    \begin{aligned}
        \min_{\alpha \in \mathbb{R}^n, s \in \mathbb{R}} \quad & s \\
        \textrm{s.t.} \quad & \cvar_\delta[c\left(x^{\alpha}, D\right)] - \opt_\delta\left(\mathrm{Ber}\left(\mu\right)\right) \leq t^* + \epsilon \quad & \forall \mu \in M \\
        \quad & \mathbb{E}_{D_1^n, x^\alpha}[\Psi_n^\alpha] - \opt_\delta\left(\mathrm{Ber}\left(\mu\right)\right) \leq s \quad & \forall \mu \in M \\
        \quad & \sum_{i=1}^n \alpha_i = 1 \\
        \quad & \alpha_i \geq 0 \quad & \forall i \in [n]
    \end{aligned}
\end{equation*}
Initialize $M$ to the active set returned by the first-stage procedure; that set is only the initial constraint index. After each solve, evaluate both the combined-risk regret and the expected post-sample regret at every point of $\mathcal{G}_m$. If either exceeds its right-hand side by more than the separation tolerance, add a violating $\mu$ to $M$ and resolve. At termination both families hold on $\mathcal{G}_m$ up to the solver and separation tolerances, and the returned $s$ is the maximum expected post-sample regret on that grid up to those tolerances. The slack $\epsilon$ is the optimization budget in the combined-risk cap.

\section{Proof of Theorem~\ref{thm:saa_nonvanishing}}\label{prf:saa_nonvanishing}

We analyze the asymptotic SAA policy regret. The key observation is that if $\mu = 1-\qstar$, then (in the asymptotic regime) SAA facing a $\mathrm{Ber}(\mu)$ distribution picks a quantity of zero with probability $1/2$ and one with probability $1/2$. 
\begin{lemma} \label{lem:erm_asymptotic_error}
    Let $D_1,\ldots,D_n \overset{i.i.d.}{\sim} \mathrm{Ber}(\mu)$ for some $\mu\in(0,1)$, and let $k_n:=\lceil (1-\mu)n\rceil$.
    Write $X_{(k_n)}$ for the $k_n$-th order statistic of $(D_1,\ldots,D_n)$ on $\{0,1\}$.
    Then
    \[
    \Pr\left(X_{(k_n)}=0\right) \underset{n \to \infty}{\to} \frac{1}{2}
    \qquad\text{and}\qquad
    \Pr\left(X_{(k_n)}=1\right) \underset{n \to \infty}{\to} \frac{1}{2}.
    \]
\end{lemma}
The proof of the lemma appears in Appendix~\ref{prf:erm_asymptotic_error}.
We specialize to SAA at the critical Bernoulli parameter $\mu^*:=1-\qstar$, so $\Pr(D=0)=\qstar$ and $x^{\mathrm{SAA}}=D_{(\lceil \qstar n\rceil)}$.
Lemma~\ref{lem:erm_asymptotic_error} at $\mu^*=1-\qstar$ gives
\[
\Pr\left(x^{\mathrm{SAA}}=0\right)\to \frac12,
\qquad
\Pr\left(x^{\mathrm{SAA}}=1\right)\to \frac12 .
\]
At this parameter, $\Pr(D=0)=\qstar$; substituting $\qstar=b(1-\delta)/(h+b)$ gives $h\qstar-b(1-\delta-\qstar)=0$. For an order $x\in[0,1]$ the loss equals $hx$ when demand is zero and $-bx$ when demand is one, so the same cancellation gives $(1-\delta)\cdot\cvar_\delta[c(x,D)]=0$. Every order in $[0,1]$ is therefore optimal. Hence the expected post-sample regret of SAA is zero for every sample size. At $x=0$ the loss is identically zero, so the cutoff is zero. At $x=1$ the loss equals $-b$ on the event that demand equals one, and that event has mass $1-\qstar>\delta$, so the cutoff sits at $-b$. The two cutoffs disagree. The loss is continuous in the order and bounded, and $x^{\mathrm{SAA}}$ converges in distribution to the equal mixture of these two optimal orders. The proposition on multiple optimal actions in \cite{BesbesRiegerGeneral} therefore 
implies that the sample-risk premium tends to a positive constant. The calculation below evaluates that constant; because the post-sample term is zero, combined regret converges to it.
Future demand $D$ is independent of the sample and of $x^{\mathrm{SAA}}$. We claim that
\[
\left(1-\delta\right)\cdot \cvar_\delta[c\left(x^{\mathrm{SAA}},D\right)]
= h\cdot \Pr\left(D=0 \land x^{\mathrm{SAA}}=1\right)
- b\cdot \left[1-\delta-\Pr\left(D=0 \land x^{\mathrm{SAA}}=1\right)-\Pr\left(x^{\mathrm{SAA}}=0\right)\right]^+ .
\]
Since $\mathbb P(D=0)=\qstar$ and $x^{\mathrm{SAA}}\in\{0,1\}$, we have $\mathbb P(D=0\land x^{\mathrm{SAA}}=1)\leq \qstar<1-\delta$; hence the entire highest-loss atom is included in the CVaR tail. The positive-part term determines whether the remaining tail mass reaches the $-b$ atom, yielding the two cases below.
The idea is that the costs are accumulated from highest to lowest, until a mass of $1-\delta$ is reached. The highest cost is $h$, when we order $1$ and the demand is zero. The next highest cost is zero, when $x^{\mathrm{SAA}}=0$. The lowest cost is $-b$, when we order one and the demand is one. 
By independence, $\Pr(D=0 \land x^{\mathrm{SAA}}=1)=\qstar/2$ and $\Pr(x^{\mathrm{SAA}}=0)=1/2$ in the limit; continuity of the right-hand side in these probabilities yields
\begin{eqnarray*}
\left(1-\delta\right) \cdot \cvar_\delta[c\left(x^{\mathrm{SAA}},D\right)]
&\to& h \cdot \Pr\left(D=0\right)\cdot \frac12
-b\cdot \max\left\{0,1-\delta-\Pr\left(D=0\right)\cdot\frac12-\frac12\right\} \\
&=& \frac{h\qstar}{2}
-b\cdot \max\left\{0,1-\delta-\frac{\qstar}{2}-\frac12\right\}.
\end{eqnarray*}
Substituting $\qstar=b(1-\delta)/(h+b)$ gives
\[
\left(1-\delta\right) \cdot \RegretCVaR{\mathrm{Ber}\left(\mu^*\right)}{\pi^{\mathrm{SAA}}}
~\to~
\begin{cases}
\frac{h \cdot b \cdot \left(1-\delta\right)}{2\left(h+b\right)} \quad & \delta \geq \frac{h}{2h+b} \\
\frac{b}{2}\cdot \delta \quad & \delta < \frac{h}{2h+b}
\end{cases} .
\label{eqn:saa_bernoulli_regret_limit}
\]
Dividing the preceding display by $(1-\delta)$ establishes the stated critical-instance convergence. 
The limit is strictly positive for every fixed $\delta\in(0,1)$, so SAA is not universally consistent under combined-risk regret.

\subsection{Proof of Lemma~\ref{lem:erm_asymptotic_error}}\label{prf:erm_asymptotic_error}
On $\{0,1\}$ we have $\Pr(X_{(k_n)}<0)=\Pr(X_{(k_n)}>1)=0$.
Define $N_n(0):=\sum_{i=1}^n \mathbf{1}_{\{D_i=0\}}$. Then
\[
\{X_{(k_n)}\leq 0\}=\{N_n\left(0\right)\geq k_n\},
\qquad
N_n\left(0\right)\sim \mathrm{Bin}\left(n,1-\mu\right).
\]
Since $k_n=\lceil (1-\mu)n\rceil$ and $N_n(0)\sim \mathrm{Bin}(n,1-\mu)$, the CLT gives
\[
\frac{N_n\left(0\right)-n(1-\mu)}{\sqrt{n(1-\mu)\mu}} ~\Rightarrow~ \mathcal{N}\left(0,1\right).
\]
Therefore
\[
\Pr\left(X_{(k_n)}\leq 0\right)=\Pr\left(N_n\left(0\right)\geq k_n\right)
=
\Pr\!\left(
\frac{N_n\left(0\right)-n(1-\mu)}{\sqrt{n(1-\mu)\mu}}
\geq
\frac{k_n-n(1-\mu)}{\sqrt{n(1-\mu)\mu}}
\right).
\]
Since $0\leq k_n-n(1-\mu)<1$, we have $(k_n-n(1-\mu))/\sqrt{n(1-\mu)\mu}\to 0$, hence
\[
\Pr\left(N_n\left(0\right)\geq k_n\right)\to\frac12.
\]
The threshold mass $\Pr(N_n(0)=k_n)$ is $O(n^{-1/2})$, so the use of a weak inequality does not affect this limit.
On $\{0,1\}$, $\Pr(X_{(k_n)}=0)=\Pr(X_{(k_n)}\leq 0)$ and
$\Pr(X_{(k_n)}=1)=1-\Pr(X_{(k_n)}=0)\to\frac12$.

\section{Proof of Theorem~\ref{thm:mcc_vanishing}}\label{prf:mcc_vanishing}
We now prove Theorem~\ref{thm:mcc_vanishing}. By Theorem~\ref{thrm:worst_case_bernoulli}, it is enough to prove the result for Bernoulli demand. Write
\[
p:=\Pr(D=0)=1-\mu,
\]
and let $N_n(0):=\sum_{i=1}^n \mathbf{1}_{\{D_i=0\}}$. Since Bernoulli observations take values in $\{0,1\}$, the increasing order statistics form a deterministic two-block pattern: exactly $N_n(0)$ zeros followed by $n-N_n(0)$ ones. Equivalently, for each realized sample,
\[
D_{(i)}=0 \text{ for } i\leq N_n(0),
\qquad
D_{(i)}=1 \text{ for } i>N_n(0).
\]
(For a fixed rank $i$, the event $\{i\leq N_n(0)\}$ is random, so $D_{(i)}$ is random; what is exact is the block structure once $N_n(0)$ is fixed.) Define
\[
H_n(z):=\sum_{i=z+1}^n \alpha_i^{(n)},
\qquad z\in\{0,\ldots,n\}.
\]
Then $x^{\mathrm{win}}=H_n(N_n(0))$, and for all $z,z'\in\{0,\ldots,n\}$,
\[
|H_n(z)-H_n(z')|
\leq
\|\boldsymbol{\alpha}^{(n)}\|_\infty |z-z'|.
\]
Let $z_p\in\{0,\ldots,n\}$ be an integer nearest $np$, and set
$\bar x_p:=H_n(z_p)$. Since $\operatorname{Var}(N_n(0))=np(1-p)\leq n/4$,
\[
\mathbb{E}|x^{\mathrm{win}}-\bar x_p|
\leq
\|\boldsymbol{\alpha}^{(n)}\|_\infty
\left(\frac{\sqrt n}{2}+\frac12\right).
\]
The loss is $(h+b)$-Lipschitz in the action, so Lemma~\ref{lem:bound_cvar_diff} gives
\[
\left|
\cvar_\delta[c(x^{\mathrm{win}},D)]
-
\cvar_\delta[c(\bar x_p,D)]
\right|
\leq
\frac{h+b}{1-\delta}\,
\|\boldsymbol{\alpha}^{(n)}\|_\infty
\left(\frac{\sqrt n}{2}+\frac12\right).
\]

\paragraph{Near-critical regime.} We first establish a deterministic proxy bound for Bernoulli demand. Fix a generic deterministic order quantity $x\in[0,1]$; this is a dummy variable and should not be confused with the random policy action $x^{\mathrm{win}}$. Under $\mathrm{Ber}(1-p)$, the loss takes the two values $c(x,0)=hx$ and $c(x,1)=-bx$. Since $hx\geq -bx$, the worst $(1-\delta)$-tail gives
\begin{eqnarray*}
    \cvar_\delta[c(x,D)] &=&
    \begin{cases}
    \frac{p hx+(1-\delta-p)(-bx)}{1-\delta}, & p<1-\delta,\\[1.2ex]
    hx, & p\geq 1-\delta .
    \end{cases} \\
    &\underset{\qstar = \frac{b (1-\delta)}{h+b}}{=}& \begin{cases}
    \frac{h+b}{1-\delta}x(p-\qstar), & p<1-\delta,\\[1.2ex]
    hx, & p\geq 1-\delta .
    \end{cases}
\end{eqnarray*}
The oracle value is $\opt_\delta(\mathrm{Ber}(1-p))=\min_{x\in[0,1]}\cvar_\delta[c(x,D)]$. Since $\qstar<1-\delta$, minimization yields
\[
\opt_\delta(\mathrm{Ber}(1-p))=
\begin{cases}
\frac{h+b}{1-\delta}(p-\qstar), & p<\qstar,\\[1.2ex]
0, & p\geq \qstar .
\end{cases}
\]
For any fixed $x\in[0,1]$, we verify the suboptimality bound
\[
\cvar_\delta[c(x,D)]-\opt_\delta(\mathrm{Ber}(1-p))
\leq \frac{h+b}{1-\delta}|p-\qstar|
\]
in three cases. If $p<\qstar$, then also $p<1-\delta$, so
\[
\cvar_\delta[c(x,D)]-\opt_\delta(\mathrm{Ber}(1-p))
=
\frac{h+b}{1-\delta}(p-\qstar)(x-1)
=
\frac{h+b}{1-\delta}(\qstar-p)(1-x)
\leq
\frac{h+b}{1-\delta}|p-\qstar|.
\]
If $\qstar\leq p<1-\delta$, then $\opt_\delta(\mathrm{Ber}(1-p))=0$ and
\[
\cvar_\delta[c(x,D)]-\opt_\delta(\mathrm{Ber}(1-p))
=
\frac{h+b}{1-\delta}x(p-\qstar)
\leq
\frac{h+b}{1-\delta}(p-\qstar)
=
\frac{h+b}{1-\delta}|p-\qstar|.
\]
If $p\geq 1-\delta$, then $\opt_\delta(\mathrm{Ber}(1-p))=0$ and $\cvar_\delta[c(x,D)]=hx$. Since
\[
p-\qstar\geq (1-\delta)-\qstar=\frac{h(1-\delta)}{h+b},
\]
we have $(h+b)(p-\qstar)/(1-\delta)\geq h$, and therefore
\[
\cvar_\delta[c(x,D)]-\opt_\delta(\mathrm{Ber}(1-p))
=
hx
\leq
h
\leq
\frac{h+b}{1-\delta}(p-\qstar)
=
\frac{h+b}{1-\delta}|p-\qstar|.
\]
Applying the bound with $x=\bar x_p$ and using $|p-\qstar|\leq 2f(n)/n$ in the near-critical regime gives
\[
\cvar_\delta[c(\bar x_p,D)]-\opt_\delta(\mathrm{Ber}(1-p))
\leq \frac{2(h+b)}{1-\delta}\frac{f(n)}{n}.
\]
Combining this with the comparison between $x^{\mathrm{win}}$ and $\bar x_p$, we obtain
\begin{eqnarray*}
    \RegretCVaR{\mathrm{Ber}(1-p)}{\pi^{\mathrm{win}}} &=&
    \cvar_\delta[c(x^{\mathrm{win}},D)]-\opt_\delta(\mathrm{Ber}(1-p)) \\
    &\leq& \left|\cvar_\delta[c(x^{\mathrm{win}},D)]-\cvar_\delta[c(\bar x_p,D)]\right| \\
    &&+ \left[\cvar_\delta[c(\bar x_p,D)]-\opt_\delta(\mathrm{Ber}(1-p))\right] \\
    &\leq&  \frac{2(h+b)}{1-\delta}\frac{f(n)}{n}
    +\frac{h+b}{1-\delta}
    \|\boldsymbol{\alpha}^{(n)}\|_\infty
    \left(\frac{\sqrt n}{2}+\frac12\right)
    ~~\to~~ 0
\end{eqnarray*}
where the first term controls the suboptimality of the deterministic proxy $\bar x_p$, and the second term controls the gap between that proxy and the random policy action.

\paragraph{Far-from-critical regime.} Suppose first that $p>\qstar+2f(n)/n$. Since $z_p\geq np-\frac12$, for all sufficiently large $n$,
\[
z_p>n\qstar+f(n)\geq\max I_n,
\]
so $\bar x_p=H_n(z_p)=0$, the oracle action. If $p<\qstar-2f(n)/n$, then $z_p\leq np+\frac12$, and for all sufficiently large $n$,
\[
z_p<n\qstar-f(n)\leq\min I_n,
\]
so $\bar x_p=H_n(z_p)=1$, again the oracle action. In either case, the CVaR comparison above gives
\[
\RegretCVaR{\mathrm{Ber}(1-p)}{\pi^{\mathrm{win}}}
=
O\!\left(\sqrt n\,\|\boldsymbol{\alpha}^{(n)}\|_\infty\right).
\]

Combining the near- and far-critical regimes, there exists a constant $C$ that depends only on $(h,b,\delta)$ such that
\[
\sup_{\mu\in[0,1]}\RegretCVaR{\mathrm{Ber}(\mu)}{\pi^{\mathrm{win}}}
\leq
C\left(
\frac{f(n)}{n}
+\sqrt n\,\|\boldsymbol{\alpha}^{(n)}\|_\infty
+\|\boldsymbol{\alpha}^{(n)}\|_\infty
\right).
\]
The right-hand side converges to zero because $f(n)=o(n)$ and $\sqrt n\|\boldsymbol{\alpha}^{(n)}\|_\infty\to0$.

\paragraph{Worst-case learning rate (Theorem~\ref{thm:mcc_vanishing}).}
Under the window-flat hypothesis $\|\boldsymbol{\alpha}^{(n)}\|_\infty=O\!\bigl(1/f(n)\bigr)$, the preceding bound gives
\[
\sup_{\mu\in[0,1]}\RegretCVaR{\mathrm{Ber}(\mu)}{\pi^{\mathrm{win}}}
=
O\!\left(\frac{f(n)}{n}+\frac{\sqrt n}{f(n)}\right).
\]
Theorem~\ref{thrm:worst_case_bernoulli} extends this bound to $\mathcal{F}_1$.

The lower bounds for window policies below reduce, after a conditional-Jensen step, to deterministic actions against Bernoulli demand. The next lemma records the resulting closed-form regret lower bound in terms of the action error.

\begin{lemma}[Bernoulli action-error lower bound]\label{lem:bernoulli_action_error_lb}
Let $D\sim\mathrm{Ber}(1-p)$ on $\{0,1\}$ with $p=\Pr(D=0)\in(0,1-\delta)$, and let $x\in[0,1]$ be deterministic. Then
\[
\cvar_\delta[c(x,D)]-\opt_\delta(\mathrm{Ber}(1-p))
~\geq~
\begin{cases}
\dfrac{h+b}{1-\delta}(p-\qstar)\,x,&p\geq\qstar,\\[1ex]
\dfrac{h+b}{1-\delta}(\qstar-p)\,(1-x),&p<\qstar.
\end{cases}
\]
\end{lemma}
\begin{proof}{Proof}
The closed forms derived in the near-critical regime above give, for $p\in(0,1-\delta)$,
\[
\cvar_\delta[c(x,D)]
=\frac{h+b}{1-\delta}x(p-\qstar)
\]
and
\[
\opt_\delta(\mathrm{Ber}(1-p))
=
\begin{cases}
\dfrac{h+b}{1-\delta}(p-\qstar),&p<\qstar,\\[1ex]
0,&p\geq\qstar.
\end{cases}
\]
Subtracting yields the claim.
\end{proof}

\paragraph{Lower bounds.}
Fix $K<\infty$ such that $\|\boldsymbol{\alpha}^{(n)}\|_\infty\leq K/f(n)$ eventually, and set $a:=1/(8K)$. The total weight on indices satisfying $|i-n\qstar|<af(n)$ is at most
\[
\bigl(2af(n)+1\bigr)\frac{K}{f(n)}=\frac14+o(1).
\]
Consequently, along every $n$, at least one of the two sets
\[
I_n^+:=\{i:i-n\qstar\geq af(n)\},
\qquad
I_n^-:=\{i:n\qstar-i\geq af(n)\}
\]
has total weight at least $3/8-o(1)$.

The certificate will be one of the following two explicit Bernoulli distributions. Set
\[
\epsilon_n:=\frac{af(n)}{2n},
\qquad
F_n^+:=\mathrm{Ber}(1-\qstar-\epsilon_n),
\qquad
F_n^-:=\mathrm{Ber}(1-\qstar+\epsilon_n).
\]
Since $f(n)=o(n)$, $\epsilon_n\to0$; as $0<\qstar<1-\delta$, both distributions are valid and their zero probabilities lie in $(0,1-\delta)$ eventually.
Thus, writing $p=\Pr(D=0)$, we have $p=\qstar+\epsilon_n$ under $F_n^+$ and $p=\qstar-\epsilon_n$ under $F_n^-$. For either distribution, write $X=x^{\mathrm{win}}$ and $\bar x=\mathbb{E}[X]$.
Because $X$ is independent of the out-of-sample demand and the Bernoulli endpoint losses are linear in $X$,
\[
\mathbb{E}[c(X,D)\mid D]=c(\bar x,D).
\]
CVaR decreases under conditional expectation: indeed, conditional Jensen gives
$\mathbb{E}[(\mathbb{E}[L\mid D]-t)^+]\leq\mathbb{E}[(L-t)^+]$ for every $t$, and the claim follows by minimizing the variational formula for CVaR over $t$. Therefore
\[
\cvar_\delta[c(X,D)]\geq\cvar_\delta[c(\bar x,D)].
\]
Applying Lemma~\ref{lem:bernoulli_action_error_lb} at $x=\bar x$ therefore gives
\[
\RegretCVaR{\mathrm{Ber}(1-p)}{\pi^{\mathrm{win}}}
\geq
\begin{cases}
\dfrac{h+b}{1-\delta}(p-\qstar)\mathbb{E}[X],&p\geq\qstar,\\[1ex]
\dfrac{h+b}{1-\delta}(\qstar-p)\mathbb{E}[1-X],&p<\qstar.
\end{cases}
\]

If $I_n^+$ has weight at least $3/8-o(1)$, use $F_n^+$, so $p_n=\qstar+\epsilon_n$. For every $i\in I_n^+$, we have $i-np_n\geq af(n)/2$; hence Chebyshev's inequality gives $\Pr(D_{(i)}=1)=1-O(n/f(n)^2)=1-o(1)$, uniformly over these indices. Thus $\mathbb{E}[X]\geq3/8-o(1)$. If instead $I_n^-$ has weight at least $3/8-o(1)$, use $F_n^-$, so $p_n=\qstar-\epsilon_n$. The same argument gives $\Pr(D_{(i)}=0)=1-o(1)$ uniformly over $i\in I_n^-$, and hence $\mathbb{E}[1-X]\geq3/8-o(1)$. Therefore, for each sufficiently large $n$, choose $F_n^+$ or $F_n^-$ according to the applicable case; then
\[
\sup_{F\in\mathcal{F}_1}\RegretCVaR{F}{\pi^{\mathrm{win}}}
=\Omega\!\left(\frac{f(n)}{n}\right).
\]
For $f(n)=n^{3/4}$, this lower bound and the preceding upper bound yield
\[
\sup_{F\in\mathcal{F}_1}\RegretCVaR{F}{\pi^{\mathrm{win}}}
=\Theta(n^{-1/4}).
\]

\subsection{Square integrability of the sample quantile}\label{prf:quantile_ui}

The uniform, exponential, lognormal, and gamma panels rely on uniform integrability of the squared sample-average order. Each of those distributions is supported on $[0,1]$ or on $[0,\infty)$, the critical quantile lies in the interior of the support, the density is continuous and positive there, and the fourth moment is finite. For each of them the sample-average order is the order statistic in the following statement, so the conclusion is uniform integrability of $n(x^{\mathrm{SAA}}-x^*_{F})^{2}$.
\begin{lemma}[Square integrability of the sample quantile]\label{lem:quantile_ui}
Let $q\in(0,1)$ and let $D$ have distribution $F$, supported on $[0,\infty)$ or on $[0,1]$. Suppose that $x^*=F^{-1}(q)$ lies in the interior of the support, that $F'$ is continuous on a neighborhood of $x^*$ inside the support with $F'(x^*)>0$, and that $\mathbb{E}[D^4]<\infty$. Let $\xi_n=D_{(\lceil q n\rceil)}$ be the corresponding order statistic of an i.i.d.\ sample from $F$. Then $\bigl\{n(\xi_n-x^*)^2\bigr\}$ is uniformly integrable.
\end{lemma}

\begin{proof}{Proof}
Write $Z_n=\sqrt{n}\,(\xi_n-x^*)$ and $k_n=\lceil qn\rceil$. Then $qn\le k_n<qn+1$, so $k_n\ge qn$ and $k_n-1\le qn$. Set $\gamma=F'(x^*)/2$. Since $F'$ is continuous at $x^*$ and $F'(x^*)>0$, shrink $\eta>0$ until $[x^*-\eta,x^*+\eta]$ lies in the interior of the support and $F'\ge \gamma$ throughout that interval. The mean value theorem gives, for every $t\in(0,\eta]$,
\[
F(x^*+t)-q ~~\geq~~ \gamma t,
\qquad
q-F(x^*-t) ~~\geq~~ \gamma t.
\]
The fixed gap $\beta=\min\{F(x^*+\eta)-q,\,q-F(x^*-\eta)\}$ then satisfies $\beta\ge \gamma\eta>0$. Every demand is nonnegative, and $0\le\xi_n\le\max_i D_i$, so
\[
\mathbb{E}\bigl[(\xi_n-x^*)^4\bigr]
~~\leq~~
8n\,\mathbb{E}[D^4]+8(x^*)^4
~~<~~
\infty.
\]
Thus $\mathbb{E}[Z_n^2]<\infty$ for every $n$.

For $0<t\le\eta$ let $S$ be binomial with parameters $n$ and $F(x^*+t)$, and let $S'$ be binomial with parameters $n$ and $F(x^*-t)$. Then $\xi_n>x^*+t$ if and only if $S\le k_n-1$, and $\xi_n\le x^*-t$ if and only if $S'\ge k_n$. The mean of $S$ is at least $n(q+\gamma t)$ and $k_n-1\le nq$, so $\{\,S\le k_n-1\,\}$ is contained in the event that $S$ lies at least $n\gamma t$ below its mean. The mean of $S'$ is at most $n(q-\gamma t)$ and $k_n\ge nq$, so $\{\,S'\ge k_n\,\}$ is contained in the event that $S'$ lies at least $n\gamma t$ above its mean. Hoeffding's inequality \citep{Hoeffding1963} therefore yields
\[
\mathbb{P}\bigl(|Z_n|>u\bigr)
~~\leq~~
2\exp(-2\gamma^2 u^2),
\qquad
0<u\le\eta\sqrt{n}.
\]
At the fixed level $\eta$ the same comparison, with gap $\beta$ in place of $\gamma t$, yields
\[
\mathbb{P}\bigl(|\xi_n-x^*|>\eta\bigr)
~~\leq~~
2\exp(-2n\beta^2).
\]
Cauchy--Schwarz and the fourth-moment bound give
\begin{align*}
\mathbb{E}\bigl[Z_n^2\mathbf{1}_{\{|\xi_n-x^*|>\eta\}}\bigr]
&~~=~~
n\,\mathbb{E}\bigl[(\xi_n-x^*)^2\mathbf{1}_{\{|\xi_n-x^*|>\eta\}}\bigr] \\
&~~\leq~~
n\bigl(8n\,\mathbb{E}[D^4]+8(x^*)^4\bigr)^{1/2}
\bigl(2\exp(-2n\beta^2)\bigr)^{1/2} \\
&~~\leq~~
\bigl(4\sqrt{\mathbb{E}[D^4]}\,n^{3/2}+4(x^*)^2 n\bigr)\exp(-n\beta^2).
\end{align*}
The right-hand side tends to zero as $n\to\infty$.

For $A>0$, the tail-integration identity gives
\[
\begin{aligned}
\mathbb{E}\bigl[Z_n^2\mathbf{1}_{\{|Z_n|>A\}}\bigr]
={}&A^2\mathbb{P}(|Z_n|>A)\\
&+\int_A^\infty 2u\,\mathbb{P}(|Z_n|>u)\,du.
\end{aligned}
\]
For $A\leq\eta\sqrt{n}$, the Hoeffding bound implies
\[
\begin{aligned}
&A^2\mathbb{P}(|Z_n|>A)
+\int_A^{\eta\sqrt n}2u\,\mathbb{P}(|Z_n|>u)\,du\\
&\qquad\leq
2A^2e^{-2\gamma^2A^2}
+\int_A^\infty4u e^{-2\gamma^2u^2}\,du\\
&\qquad=
(2A^2+\gamma^{-2})e^{-2\gamma^2A^2}.
\end{aligned}
\]
Moreover,
\[
\int_{\eta\sqrt n}^\infty
2u\,\mathbb{P}(|Z_n|>u)\,du
\leq
\mathbb{E}\bigl[
Z_n^2\mathbf{1}_{\{|\xi_n-x^*|>\eta\}}
\bigr],
\]
which tends to zero as $n\to\infty$ by the preceding bound.

Consequently, for $A\leq\eta\sqrt n$,
\[
\begin{aligned}
\mathbb{E}\bigl[Z_n^2\mathbf{1}_{\{|Z_n|>A\}}\bigr]
\leq{}&
(2A^2+\gamma^{-2})e^{-2\gamma^2A^2}\\
&+\mathbb{E}\bigl[
Z_n^2\mathbf{1}_{\{|\xi_n-x^*|>\eta\}}
\bigr].
\end{aligned}
\]
For $A>\eta\sqrt n$, the left-hand side is bounded directly by the same far-event contribution.

The first term vanishes as $A\to\infty$, while the second vanishes as $n\to\infty$. Since $Z_n^2$ is integrable for every fixed $n$, the finitely many remaining sample sizes can be handled by taking $A$ sufficiently large. Therefore $\{Z_n^2\}$ is uniformly integrable.
\Halmos
\end{proof}

\subsection{Fixed-distribution asymptotics of SAA}
\label{apx:saa_fixed_rates}

\begin{proposition}[SAA at fixed demand distributions]
\label{prop:saa_fixed_rates}
Fix $h,b>0$ and $\delta\in(0,1)$. Let $F$ be any
of the four distributions considered in
Section~\ref{sec:saa_fixed}: uniform, exponential,
lognormal, or gamma demand, with the parameters
specified there. Then
\[
\mathrm{PS}_{\delta,n}(\pi^{\mathrm{SAA}},F)
=\Theta(n^{-1}),
\qquad
\mathrm{SR}_{\delta,n}(\pi^{\mathrm{SAA}},F)
=\Theta(n^{-1/2}).
\]
Consequently,
\[
\RegretCVaR{F}{\pi^{\mathrm{SAA}}}
=\Theta(n^{-1/2}).
\]
\end{proposition}

\begin{proof}{Proof}
Each of the four distributions has a unique interior
oracle order $x_F^*$ at which the density is continuous,
positive, and finite. SAA selects the order statistic
$D_{(\lceil n\qstar\rceil)}$. By
\cite[Corollary~21.5]{vanDerVaart1998},
\[
\sqrt{n}(x^{\mathrm{SAA}}-x_F^*)
\Rightarrow
N\left(
0,\frac{\qstar(1-\qstar)}{F'(x_F^*)^2}
\right),
\]
with a nondegenerate limit.

All four distributions have finite fourth moments
and are supported on $[0,1]$ or $[0,\infty)$.
Lemma~\ref{lem:quantile_ui} therefore establishes
uniform integrability of
$\{n(x^{\mathrm{SAA}}-x_F^*)^2\}$.
Together with the quantile CLT, this verifies
the concentration conditions in
equation~\eqref{eqn:concentration} with
$r_n=n^{-1/2}$.

By Lemma~\ref{lem:nv_atom}, the loss distribution
at the oracle has an atom at $-bx_F^*$ that strictly
straddles the CVaR cutoff. The required isolation
conditions are verified in
Appendix~\ref{prf:nv_atom}. Moreover, the second
derivative of the outcome-only CVaR objective at
the oracle is
\[
\frac{(h+b)F'(x_F^*)}{1-\delta},
\]
which is positive and finite.

Thus, the conditions of the atomic-regime result
of~\cite{BesbesRiegerGeneral} are satisfied.
Applying that result with $r_n=n^{-1/2}$ gives
\[
\mathrm{PS}_{\delta,n}(\pi^{\mathrm{SAA}},F)
=\Theta(n^{-1}),
\qquad
\mathrm{SR}_{\delta,n}(\pi^{\mathrm{SAA}},F)
=\Theta(n^{-1/2}).
\]
Since combined-risk regret is the sum of these
two nonnegative components, it follows that
$\RegretCVaR{F}{\pi^{\mathrm{SAA}}}
=\Theta(n^{-1/2})$.
\end{proof}

\section{\texorpdfstring{{Finite-sample computation for the fixed-demand comparisons}}{Finite-sample computation for the fixed-demand comparisons}}\label{apx:fixed_demand_mc}

Figures~\ref{fig:saa_decomp_distributions} and~\ref{fig:saa_vs_policy_distributions} use the four demand distributions of Section~\ref{sec:saa_fixed}: uniform on $[0,1]$, exponential with rate $1$, lognormal with $\mu_{\log}=1$ and $\sigma_{\log}=1.8$, and gamma with shape $4$ and scale $0.125$. Each distribution is untruncated. The oracle order is $x_F^*=F^{-1}(\qstar)$ from equation~\eqref{eqn:gotoh_takano_profit}, and $\opt_\delta(F)=\cvar_\delta[c(x_F^*,D)]$. Proposition~\ref{prop:saa_fixed_rates} verifies the fixed-distribution rate conclusions for these four distributions.

For each sample size $n$ stored with the stage-two weights, draw $N_n$ replications $(D_{r,1},\ldots,D_{r,n},D_{r,n+1})$. One uniform sample is mapped to all four distributions by the inverse cdf. The balancing policy applies the stage-two weights to the ascending order statistics. SAA places unit weight on $D_{r,(k_n)}$ with $k_n=\lceil n\qstar\rceil$.

The oracle value and the conditional tail are computed from the demand law, with the order held fixed. For a realized order $x$, the conditional value is the outcome-only CVaR in equation~\eqref{eqn:conditional_psi},
\[
\Psi(x)
=
\inf_{t\in\mathbb{R}}
\left\{
t+\frac{1}{1-\delta}\,\mathbb{E}\bigl[(c(x,D)-t)^+\bigr]
\right\}.
\]
On replication $r$, write $Z_r=c(x_r,D_{r,n+1})$ and let $q$ be a $\delta$-quantile of the empirical distribution of $Z_1,\ldots,Z_{N_n}$. The two plotted estimates are
\begin{align*}
\widehat R
&=
q+\frac{1}{N_n(1-\delta)}\sum_{r=1}^{N_n}(Z_r-q)_+
-\opt_\delta(F),
\\
\widehat{\mathrm{PS}}
&=
\frac{1}{N_n}\sum_{r=1}^{N_n}\Psi(x_r)
-\opt_\delta(F).
\end{align*}
The sample-risk premium estimate is $\widehat R-\widehat{\mathrm{PS}}$. If more than $2048$ distinct orders occur, $\Psi$ is evaluated on $512$ empirical quantiles of those orders and extended by linear interpolation.

The standard errors account for variation across the $N_n$ replications. Let $e_r=(Z_r-q)_+$ and let $\bar e$ be their average. Then
\[
\operatorname{se}(\widehat R)
=
\frac{1}{(1-\delta)\sqrt{N_n}}
\Biggl(
\frac{1}{N_n-1}\sum_{r=1}^{N_n}(e_r-\bar e)^2
\Biggr)^{1/2}.
\]
These standard errors measure variation across replications only. The oracle value and each conditional tail are retained only when the quadrature error estimate is at most $10^{-5}$, which is at most $0.0002$ on the CVaR scale at $\delta=.95$. On a fresh sample, linear interpolation of $\Psi$ changes its average by at most $1.7\times 10^{-4}$ in absolute value.
The paired standard error of the balancing policy minus SAA uses the same expression on the paired excesses. For each $n$, $N_n$ begins at $10{,}000$ and increases until the Bernoulli, stability, and paired-error checks hold, or until $N_n=10^7$.

\section{Additional figures}\label{apx:additional_figures}
The decomposition panels and the first-stage comparison on typical demand are recorded here. Uniform demand on $[0,1]$ lies in $\mathcal{F}_1$. Exponential demand with rate $1$ is unbounded and lies outside $\mathcal{F}_1$. The decomposition panels also include lognormal demand with $\mu_{\log}=1$ and $\sigma_{\log}=1.8$ and gamma demand with shape $4$ and scale $0.125$, which are likewise unbounded and outside $\mathcal{F}_1$. The first-stage comparison uses uniform demand and exponential demand with rate $1$.

\FIGURE
{%
\begin{minipage}{\linewidth}
\centering
    \begin{minipage}{0.48\linewidth}
        \centering
        \includegraphics[width=224.454bp,height=118.525bp]{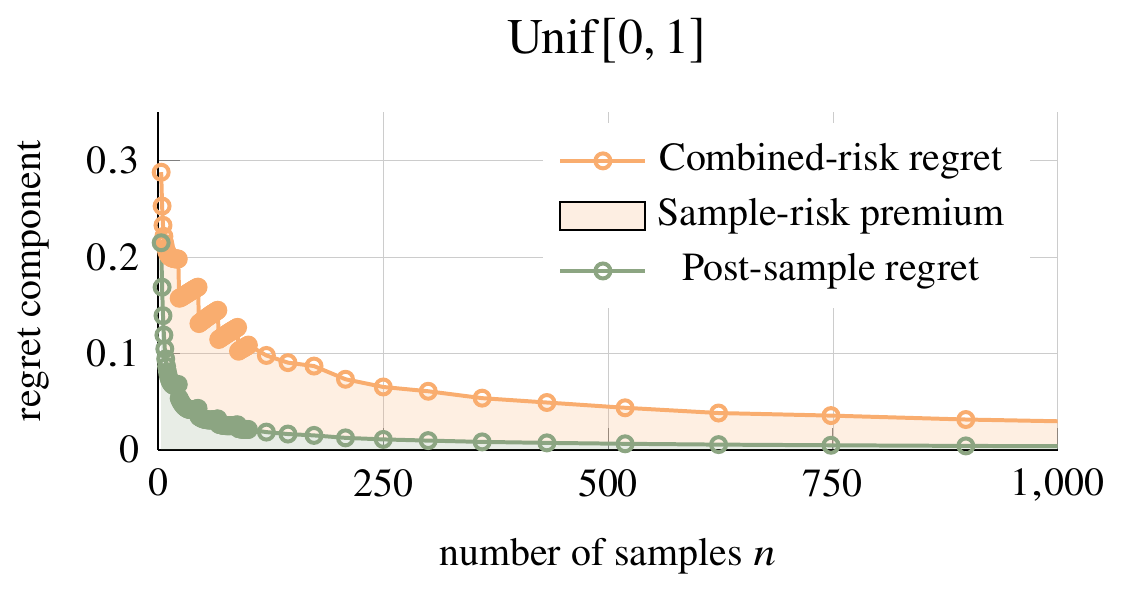}
    \end{minipage}
    \hfill
    \begin{minipage}{0.48\linewidth}
        \centering
        \includegraphics[width=224.454bp,height=118.722bp]{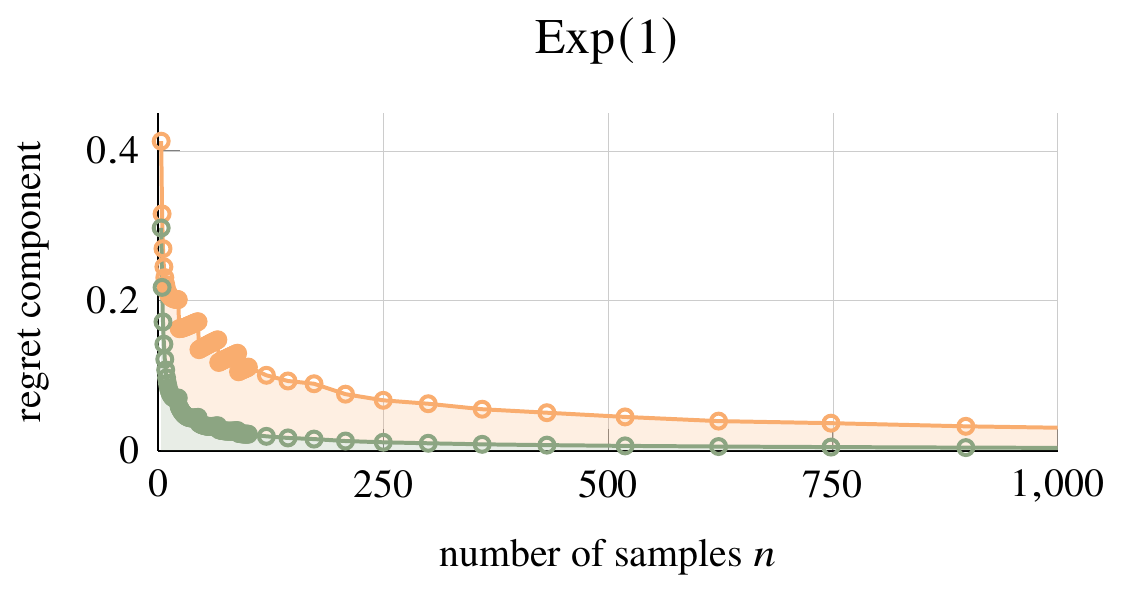}
    \end{minipage}\\[0.8em]
    \begin{minipage}{0.48\linewidth}
        \centering
        \includegraphics[width=217.911bp,height=120.201bp]{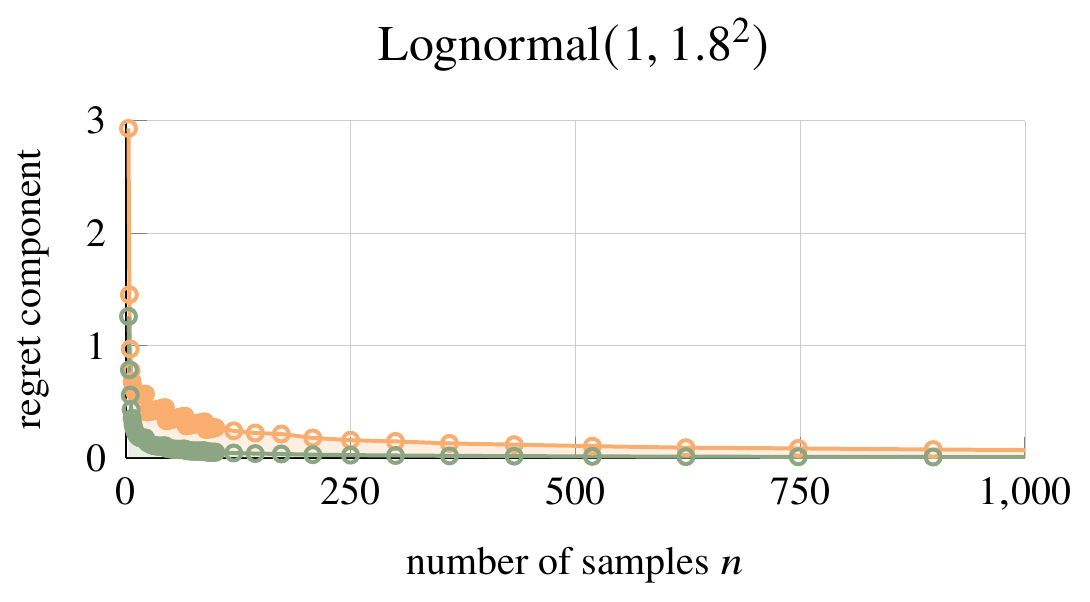}
    \end{minipage}
    \hfill
    \begin{minipage}{0.48\linewidth}
        \centering
        \includegraphics[width=224.454bp,height=118.523bp]{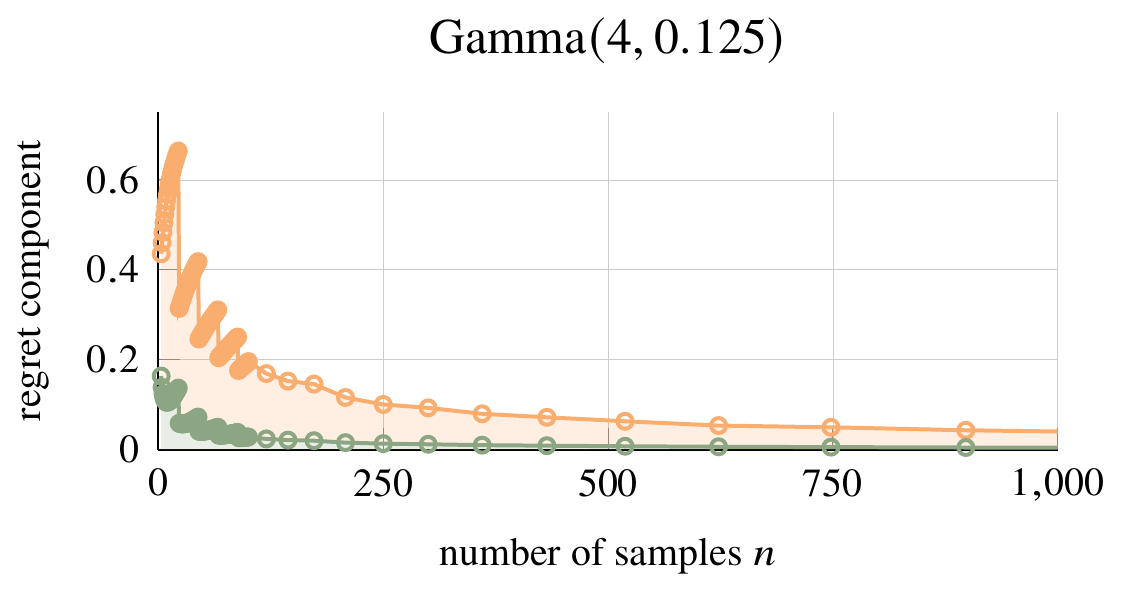}
    \end{minipage}\\
{\small $h = 1$, $b=9$, $\delta = .95$.}
\end{minipage}
}
{Monte Carlo estimates of combined-risk regret and expected post-sample regret for SAA  under uniform demand on $[0,1]$, exponential demand with rate $1$, lognormal demand with $\mu_{\log}=1$ and $\sigma_{\log}=1.8$, and gamma demand with shape $4$ and scale $0.125$. Each sample size uses between $617{,}956$ and $10^{7}$ replications, and each replication consists of the $n$ training demands and one out-of-sample demand. The Monte Carlo standard errors are at most $0.0009$ under uniform and exponential demand, $0.0020$ under gamma demand, and $0.0073$ under lognormal demand. Appendix~\ref{apx:fixed_demand_mc} records the computation.\label{fig:saa_decomp_distributions}}
{}

\FIGURE
{%
\begin{minipage}{\linewidth}
\centering
\includegraphics[width=429.264bp,height=184.025bp]{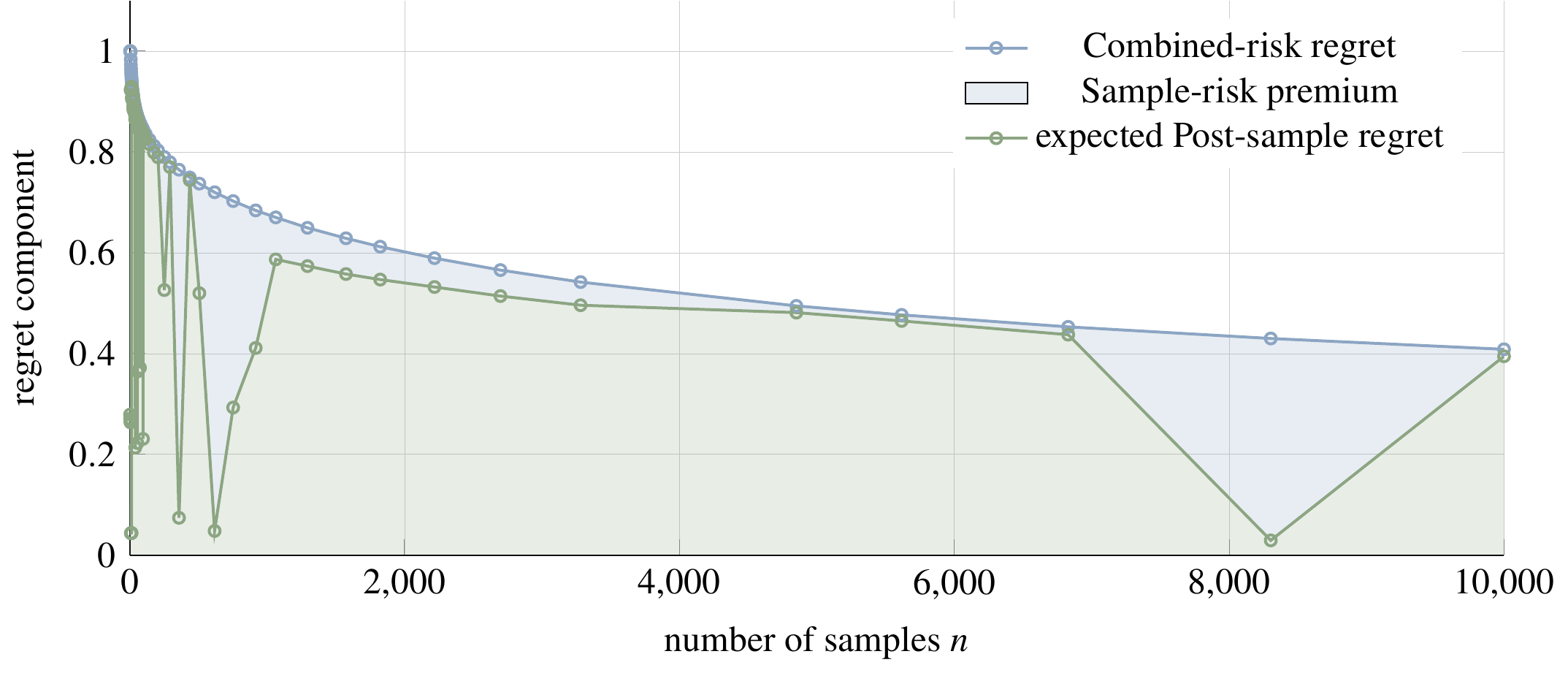}\\
{\small $h = 1$, $b=9$, $\delta = .95$.}
\end{minipage}
}
{Computed MCC policy: combined-risk regret and its decomposition at the grid distribution that maximizes combined-risk regret. \label{fig:near_opt_psi_decomp}}
{}

\FIGURE
{%
\begin{minipage}{\linewidth}
\centering
    \begin{minipage}{0.48\linewidth}
        \centering
        \includegraphics[width=225.338bp,height=174.679bp]{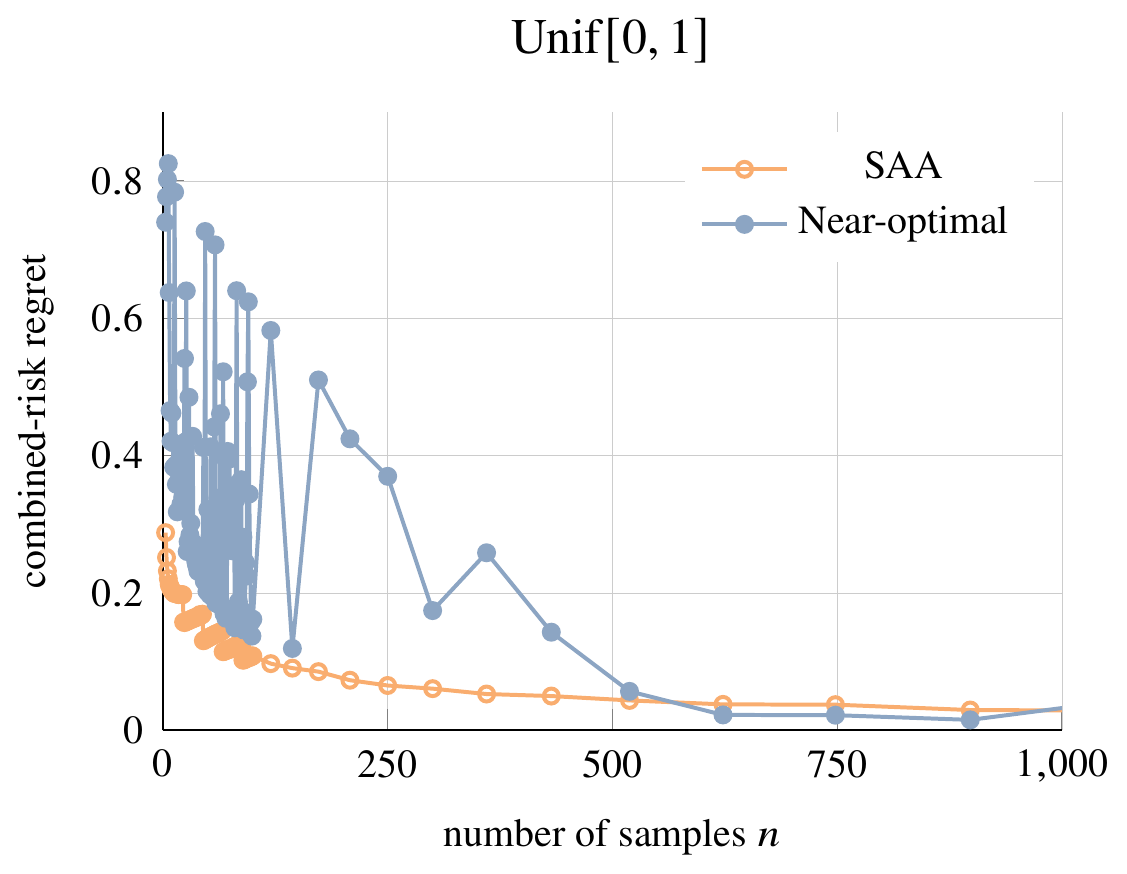}
    \end{minipage}
    \hfill
    \begin{minipage}{0.48\linewidth}
        \centering
        \includegraphics[width=218.795bp,height=174.883bp]{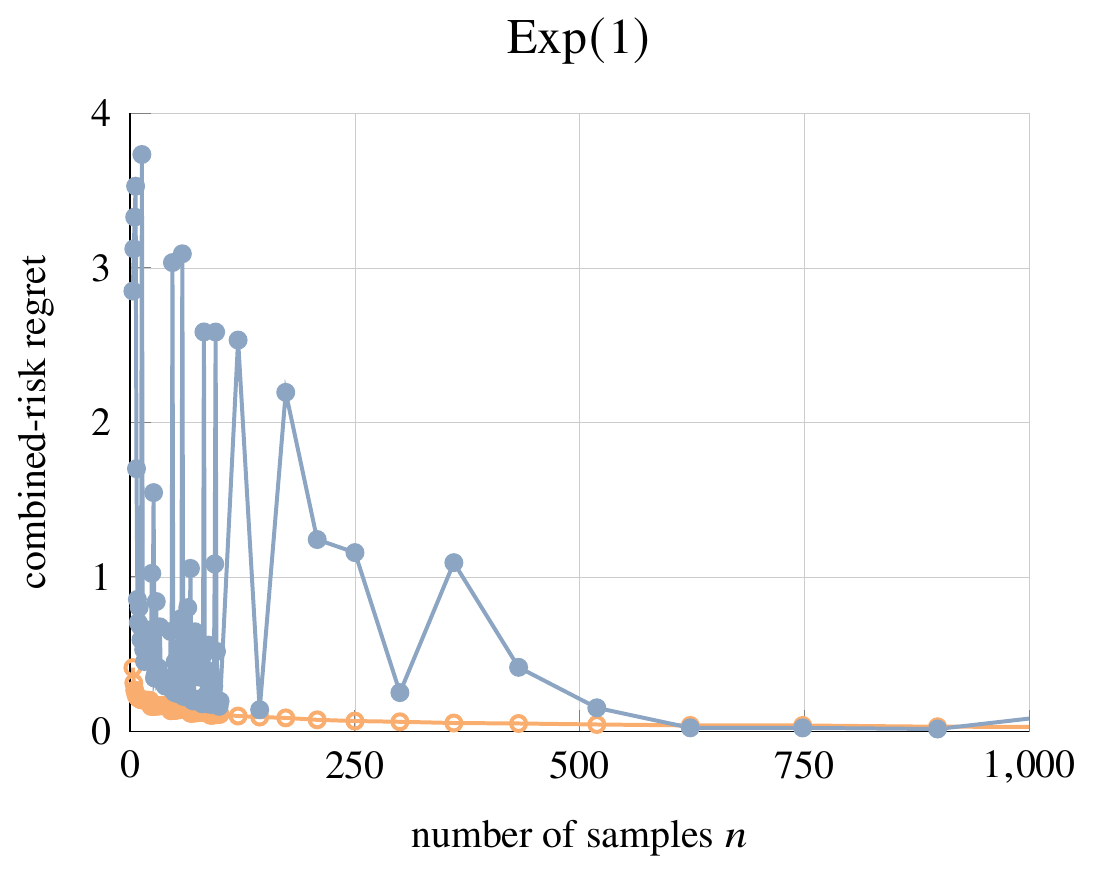}
    \end{minipage}\\
{\small $h = 1$, $b=9$, $\delta = .95$.}
\end{minipage}
}
{Monte Carlo estimates of combined-risk regret for SAA and the first-stage computed policy under uniform demand on $[0,1]$ and exponential demand with rate $1$. On both, the near-optimal policy lies above SAA over nearly the whole plotted range.\label{fig:benign_near_opt}}
{}

\FIGURE
{%
\begin{minipage}{\linewidth}
\centering
\includegraphics[width=430.094bp,height=184.024bp]{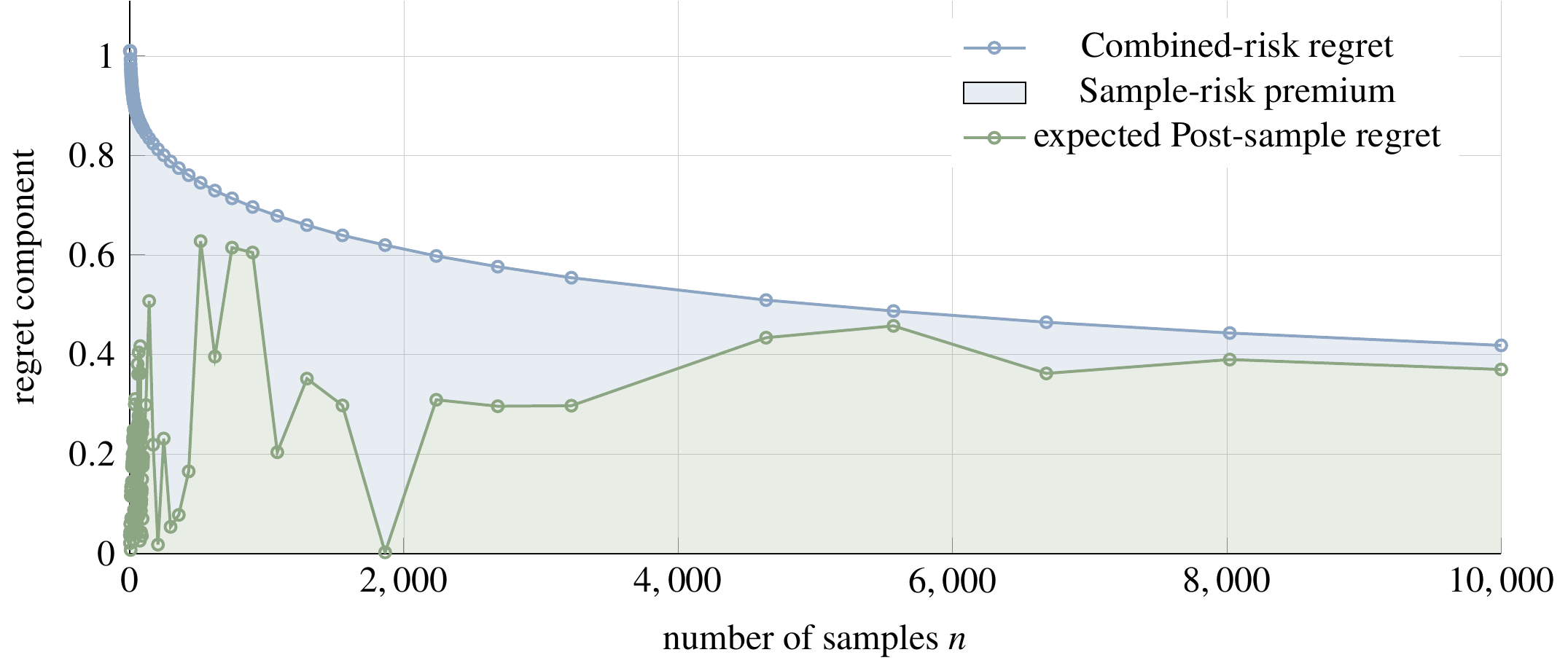}\\
{\small $h = 1$, $b=9$, $\delta = .95$.}
\end{minipage}
}
{Two-stage policy: combined-risk regret and its decomposition at the grid distribution that maximizes combined-risk regret. \label{fig:psi_decomp}}
{}

\input{TEX-Strengthened-Proofs}

%% file: TEX-Strengthened-Proofs.tex
\section{Proof of Theorem~\ref{thm:saa_exact_worst_case}}\label{prf:saa_exact_worst_case}

\begin{proof}{Proof of Theorem~\ref{thm:saa_exact_worst_case}}
By scaling, it suffices to take $T=1$. SAA is an MCC policy, and Theorem~\ref{thrm:worst_case_bernoulli} reduces its worst-case regret to
\[
\sup_{\mu\in[0,1]}\RegretCVaR{\mathrm{Ber}(\mu)}{\pi^{\mathrm{SAA}}}.
\]
Under $\mathrm{Ber}(\mu)$ one has $\Pr(D=0)=1-\mu$ and $x^{\mathrm{SAA}}=D_{(\lceil \qstar n\rceil)}$. This order equals $1$ when fewer than $\lceil \qstar n\rceil$ sample points are zero, and equals $0$ otherwise.

Future demand is independent of the sample. The combined loss equals $h$ when the order is $1$ and demand is zero, equals $0$ when the order is $0$, and equals $-b$ when the order is $1$ and demand is one. These events have probabilities $(1-\mu)\Pr(x^{\mathrm{SAA}}=1)$, $\Pr(x^{\mathrm{SAA}}=0)$, and $\mu\Pr(x^{\mathrm{SAA}}=1)$. Accumulating cost from highest to lowest until mass $1-\delta$ is reached, as in the proof of Theorem~\ref{thm:saa_nonvanishing}, gives
\[
\begin{aligned}
&(1-\delta)\cdot\cvar_\delta[c(x^{\mathrm{SAA}},D)] \\
&\qquad=
h\min\{(1-\mu)\Pr(x^{\mathrm{SAA}}=1),\,1-\delta\}
-b\big(\mu\Pr(x^{\mathrm{SAA}}=1)-\delta\big)^+.
\end{aligned}
\]
If $1-\mu\geq\qstar$, the oracle order is $0$ and the oracle value is $0$. If $1-\mu<\qstar$, the oracle order is $1$. Since $\qstar<1-\delta$, the atom of mass $1-\mu$ at loss $h$ lies entirely in the upper tail, and the oracle value equals $\frac{h+b}{1-\delta}(1-\mu-\qstar)$. In either case the oracle value is $\min\{0,\frac{h+b}{1-\delta}(1-\mu-\qstar)\}$. The regret is continuous in $\mu$, so its maximum on $[0,1]$ exists.

Suppose $1-\mu_n\to\qstar$ and $\Pr(x^{\mathrm{SAA}}=1)\to\rho$. Then $(1-\mu_n)\Pr(x^{\mathrm{SAA}}=1)\to\qstar\rho<1-\delta$, the oracle term tends to zero, and the regret tends to
\begin{equation}\label{eqn:saa_rho_limit}
\frac{h\qstar\rho-b\big((1-\qstar)\rho-\delta\big)^+}{1-\delta}.
\end{equation}
This expression increases up to $\rho=\delta/(1-\qstar)\in(0,1)$ and decreases afterwards: the two slopes are $h\qstar/(1-\delta)>0$ and $-b\delta/(1-\delta)<0$. At the kink the positive part vanishes. Substituting $\qstar=b(1-\delta)/(h+b)$ gives the value $hb\delta/(h+b\delta)$.

For the upper bound, let $\mu_n$ maximize the regret and pass to a subsequence with $\mu_n\to\mu_0\in[0,1]$. The sample fraction of zeros is a binomial average with variance at most $1/(4n)$, so Chebyshev's inequality concentrates it about $1-\mu_n$, uniformly in the parameter. Thus $\Pr(x^{\mathrm{SAA}}=1)\to1$ if $1-\mu_0<\qstar$, and $\Pr(x^{\mathrm{SAA}}=1)\to0$ if $1-\mu_0>\qstar$. In either case the tail formula and the oracle tend to the same limit, so the regret tends to zero. If $1-\mu_0=\qstar$, a further subsequence has $\Pr(x^{\mathrm{SAA}}=1)\to\rho$ for some $\rho\in[0,1]$, and the regret tends to the expression above. Applying this argument along a subsequence that attains the upper limit proves
\[
\limsup_{n\to\infty}\sup_{\mu\in[0,1]}\RegretCVaR{\mathrm{Ber}(\mu)}{\pi^{\mathrm{SAA}}}
\leq\frac{hb\delta}{h+b\delta}.
\]

To attain the bound, let $\Phi$ denote the standard normal distribution function and set
\[
1-\mu_n
=\qstar+\frac{\sqrt{\qstar(1-\qstar)}}{\sqrt{n}}\,
\Phi^{-1}\!\left(\frac{h(1-\delta)}{h+b\delta}\right).
\]
These parameters lie in $(0,1)$ for all sufficiently large $n$. Since $1-\mu_n\to\qstar\in(0,1)$, the centered Bernoulli summands are bounded and $n\mu_n(1-\mu_n)\to\infty$, which verifies the Lindeberg condition. The triangular-array central limit theorem and $\lceil\qstar n\rceil-\qstar n=O(1)$ give
\[
\Pr(x^{\mathrm{SAA}}=0)\to\frac{h(1-\delta)}{h+b\delta}.
\]
Hence $\Pr(x^{\mathrm{SAA}}=1)$ tends to $\delta/(1-\qstar)$, and the regret tends to $hb\delta/(h+b\delta)$. Scaling demand and orders by $T$ multiplies regret by $T$.
\end{proof}

\section{Proof of Theorem~\ref{thm:global_minimax_rate}}\label{prf:global_minimax_rate}

The lower bound restricts demand to a Bernoulli family. It does not restrict the policy to MCC or use the worst-case reduction. Take $T=1$ and write $\mu^*=1-\qstar$. Projecting onto $[0,1]$ lowers the loss pointwise, so $x^\pi$ may be taken in $[0,1]$. It may depend on the whole sample and on independent randomization. Write $\mathbb{E}_\mu$ for expectation under $\mathrm{Ber}(\mu)$. Combined-risk regret is the sum of the two nonnegative components, so each is at most the regret.

\paragraph{A bound on $\mathrm{SR}$.}
The sample-risk premium compares the combined upper tail, which on $\{D=1\}$ keeps only the small orders, with the risk conditional on the realized order, which charges the whole order. The gap between those two averages is the spread of $x^\pi$, and the comparison produces a multiple of its mean absolute deviation from one quantile.

For $\mu\in(\delta,1)$ the mass on $\{D=0\}$ is $1-\mu<1-\delta$, and every loss $hx^\pi$ there exceeds every loss $-bx^\pi$ on $\{D=1\}$. The upper $(1-\delta)$-tail therefore contains every outcome with $D=0$ and still needs probability $\mu-\delta$ from $\{D=1\}$. Future demand is independent of $x^\pi$, and on $\{D=1\}$ the loss $-bx^\pi$ is largest for the smallest orders, so the missing slice is a lower portion of the marginal of $x^\pi$. Let $\tau_\mu=(\mu-\delta)/\mu$ and let $m_\mu$ be a $\tau_\mu$-quantile of $x^\pi$. That portion is the distribution on $[0,1]$ which assigns to each Borel set $A$ the mass
\[
\tau_\mu^{-1}\Pr_\mu(x^\pi\in A,\ x^\pi<m_\mu)+\tau_\mu^{-1}\bigl(\tau_\mu-\Pr_\mu(x^\pi<m_\mu)\bigr)\mathbf{1}_{\{m_\mu\in A\}}.
\]
The first term keeps every order strictly below $m_\mu$, with its probability divided by $\tau_\mu$. The second term puts all of the remaining probability on the single point $m_\mu$. Writing $p=\Pr_\mu(x^\pi<m_\mu)$, the mean of this distribution is
\[
\tau_\mu^{-1}\mathbb{E}_\mu[x^\pi\mathbf{1}_{\{x^\pi<m_\mu\}}]+\tau_\mu^{-1}(\tau_\mu-p)m_\mu.
\]
On $\{x^\pi<m_\mu\}$ one has $x^\pi=m_\mu-(m_\mu-x^\pi)$, and $(m_\mu-x^\pi)^+$ vanishes off that event, so the expectation equals $pm_\mu-\mathbb{E}_\mu(m_\mu-x^\pi)^+$. Substituting and cancelling $pm_\mu$ against $(\tau_\mu-p)m_\mu$ leaves
\[
m_\mu-\mathbb{E}_\mu(m_\mu-x^\pi)^+/\tau_\mu.
\]

The combined-risk $\cvar$ uses this mean as the average order on the part of the tail where $D=1$. Its numerator in the tail average is
\[
(1-\mu)h\,\mathbb{E}_\mu x^\pi-b(\mu-\delta)\Bigl(m_\mu-\frac{\mathbb{E}_\mu(m_\mu-x^\pi)^+}{\tau_\mu}\Bigr).
\]
Since $(\mu-\delta)/\tau_\mu=\mu$, the second summand equals $-b(\mu-\delta)m_\mu+b\mu\,\mathbb{E}_\mu(m_\mu-x^\pi)^+$. The conditional risk of a fixed order $x$ fills the same slice with the constant loss $-bx$, so the expected conditional risk has numerator $\bigl((1-\mu)h-(\mu-\delta)b\bigr)\mathbb{E}_\mu x^\pi$. Subtracting and dividing by $b$ leaves
\[
(\mu-\delta)(\mathbb{E}_\mu x^\pi-m_\mu)+\mu\,\mathbb{E}_\mu(m_\mu-x^\pi)^+.
\]
The identity $\mathbb{E}_\mu x^\pi-m_\mu=\mathbb{E}_\mu(x^\pi-m_\mu)^+-\mathbb{E}_\mu(m_\mu-x^\pi)^+$ replaces the first summand and cancels $(\mu-\delta)\mathbb{E}_\mu(m_\mu-x^\pi)^+$ against part of the second, leaving weight $\delta$ there. Multiplying by $b/(1-\delta)$ gives
\[
\mathrm{SR}_{\delta,n}(\pi,\mathrm{Ber}(\mu))
=\frac{b}{1-\delta}\Bigl[(\mu-\delta)\,\mathbb{E}_\mu(x^\pi-m_\mu)^+
+\delta\,\mathbb{E}_\mu(m_\mu-x^\pi)^+\Bigr].
\]
Both coefficients are positive, so each may be lowered to $\min\{\mu-\delta,\delta\}$. The sum of the two positive parts is $\mathbb{E}_\mu|x^\pi-m_\mu|$, and therefore
\begin{equation}\label{eqn:premium_disperses}
\mathrm{SR}_{\delta,n}(\pi,\mathrm{Ber}(\mu))
\geq\frac{b}{1-\delta}\min\{\mu-\delta,\delta\}\,\mathbb{E}_\mu|x^\pi-m_\mu|.
\end{equation}

\paragraph{An absolute deviation rule for $x^\pi$.}
Fix a compact interval $J$, centered at $\mu^*$ and contained in $(\delta,1)$, on which $\mu(1-\mu)$ and $\min\{\mu-\delta,\delta\}$ are bounded below by positive constants that depend only on $(h,b,\delta)$. Let $\eta>0$ be an upper bound on the policy's worst-case regret over Bernoulli demand. The premium is at most this regret, so equation~\eqref{eqn:premium_disperses} gives $\mathbb{E}_\mu|x^\pi-m_\mu|\leq C\eta$ throughout $J$, for a constant $C$ of the same dependence. A regret of size $\eta$ therefore forces $x^\pi$ to lie within order $\eta$ of the quantile $m_\mu$.

\paragraph{Overlap of nearby parameters.}
Samples from nearby parameters overlap. The Bernoulli relative entropy satisfies
\[
\mu\log\frac{\mu}{\mu'}+(1-\mu)\log\frac{1-\mu}{1-\mu'}
\leq\frac{(\mu-\mu')^2}{\mu'(1-\mu')},
\]
since $\log u\leq u-1$. The display is the divergence of one draw, so a sample of size $n$ has divergence $n$ times as large. Pinsker's inequality, applied to that sample, and contraction under the policy kernel, which is the same at every $\mu$, bound the total variation between the induced distributions of $x^\pi$. Thus there is a constant $c>0$, depending only on $(h,b,\delta)$, such that $|\mu-\mu'|\leq c n^{-1/2}$ with $\mu,\mu'\in J$ implies that this total variation is at most $1/2$. The common part of two such distributions has mass at least $1/2$, so they may be coupled to agree with probability at least $1/2$. On that event $|m_\mu-m_{\mu'}|\leq|x^\pi-m_\mu|+|x^\pi-m_{\mu'}|$. Taking expectations and using $\mathbb{E}|x^\pi-m|\leq C\eta$ gives $|m_\mu-m_{\mu'}|\leq 4C\eta$. Each mean lies within $C\eta$ of its quantile, because $|\mathbb{E}_\mu x^\pi-m_\mu|\leq\mathbb{E}_\mu|x^\pi-m_\mu|$, and the means therefore differ by at most $6C\eta$.

\paragraph{The post-sample term.}
For $\mu\in(\delta,1)$, the conditional risk of a deterministic order $x$ is $\frac{h+b}{1-\delta}(1-\mu-\qstar)x$, and the oracle value is the minimum of zero and the same quantity at the better endpoint. The expected gap is
\begin{equation}\label{eqn:all_policy_ps_bernoulli}
\mathrm{PS}_{\delta,n}(\pi,\mathrm{Ber}(\mu))
=
\begin{cases}
\frac{h+b}{1-\delta}(\mu-\mu^*)\,\mathbb{E}_\mu(1-x^\pi),&\mu>\mu^*,\\[4pt]
\frac{h+b}{1-\delta}(\mu^*-\mu)\,\mathbb{E}_\mu x^\pi,&\mu\leq\mu^*.
\end{cases}
\end{equation}
This term forces the means apart. Set $\theta=4(1-\delta)\eta/(h+b)$. If $\eta$ is smaller than a positive constant $\eta_0(h,b,\delta)$, both points $\mu^*\pm\theta$ lie in $J$. Equation~\eqref{eqn:all_policy_ps_bernoulli} then gives $\mathbb{E}_{\mu^*+\theta}(1-x^\pi)\leq 1/4$ and $\mathbb{E}_{\mu^*-\theta}x^\pi\leq 1/4$, so one mean is at least $3/4$ and the other at most $1/4$. The two parameters differ by $2\theta$. Dividing that distance into steps of length $c n^{-1/2}$ takes at most $1+2\theta\sqrt n/c$ steps, and each step moves the mean by at most $6C\eta$, so
\[
\frac12\leq 6C\eta\Bigl(1+\frac{2\theta\sqrt n}{c}\Bigr).
\]
Since $\theta=4(1-\delta)\eta/(h+b)$, the right-hand side is a constant multiple of $\eta+\eta^2\sqrt n$, with the constant depending only on $(h,b,\delta)$. Decreasing $\eta_0$ if needed makes the multiple of $\eta$ at most $1/4$, and therefore $\eta^2\sqrt n$ is bounded below by a positive constant of the same dependence. Hence $\eta\geq c_0 n^{-1/4}$ for some such $c_0>0$. If instead $\eta\geq\eta_0$, the same conclusion holds for all sufficiently large $n$ after decreasing $c_0$. Every upper bound $\eta$ on a policy's worst-case Bernoulli regret is therefore at least $c_0 n^{-1/4}$. The infimum of those upper bounds is the worst-case Bernoulli regret itself, so that regret is at least $c_0 n^{-1/4}$ for every policy. The worst case over $\mathcal{F}_1$ is at least as large, and hence $V_n^{\mathrm{all}}(1)\geq c_0 n^{-1/4}$. Theorem~\ref{thm:mcc_vanishing} supplies the matching upper bound through uniform windows with $f(n)=n^{3/4}$, which belong to MCC and hence to the class of all policies. This establishes the bounds for every fixed $\delta\in(0,1)$ at $T=1$.

\paragraph{The case of $\delta=0$.}
Take $\delta=0$. On $[0,1]$, the loss is $c(x,D)=hx-(h+b)\min(x,D)$. If $\widehat F_n$ is the empirical CDF, the uniform error of the empirical expected objective is at most
\[
(h+b)\int_0^1|\widehat F_n(t)-F(t)|\,dt.
\]
Empirical optimization incurs regret at most twice this error. Since $\mathbb{E}|\widehat F_n(t)-F(t)|\leq1/(2\sqrt n)$, SAA's worst-case expected regret is at most $(h+b)/\sqrt n$. For the lower bound, $\qstar=b/(h+b)$. Let $\mu_\pm=(1-\qstar)\pm a/\sqrt n$ for a sufficiently small fixed $a>0$. The same relative-entropy bound keeps the total variation between these sampling distributions at most $1/2$ for all large $n$. For any policy with orders in $[0,1]$, the sum of its expected regrets at these parameters is
\[
\frac{(h+b)a}{\sqrt n}\bigl[1-\mathbb{E}_{\mu_+}x^\pi+\mathbb{E}_{\mu_-}x^\pi\bigr]
\geq\frac{(h+b)a}{2\sqrt n},
\]
because the means of $x^\pi$ differ by at most the total variation. Thus at least one regret is at least $(h+b)a/(4\sqrt n)$. This proves the $n^{-1/2}$ minimax order, both over all policies and within MCC.

\paragraph{Scaling.}
Scaling is exact for arbitrary policies and for MCC. A policy on $[0,T]$ induces a policy on $[0,1]$ by applying it to observations multiplied by $T$ and dividing its order by $T$, and conversely. The loss and CVaR are positively homogeneous, so regret is multiplied by $T$. This proves all statements for general $T>0$.

\section{\texorpdfstring{{Numerical evaluation of SAA}}{Numerical evaluation of SAA}}\label{apx:saa_numerics}

{The SAA series in Figures~\ref{fig:saa_performance_linear}, \ref{fig:saa_performance_by_delta_linear}, and~\ref{fig:finite_sample_opt_vs_saa} use the order statistic of rank $k_n=\lceil n\qstar\rceil$, where $\qstar=b(1-\delta)/(h+b)$ is computed from the exact rational values of the stated decimal parameters. When $n\qstar$ is an integer, a floating-point product can fall slightly above that integer, and the ceiling then returns $k_n+1$. The rational value returns $k_n$. The same rule applies at low risk and small $n$: for $\delta=0$ and $n=3,\ldots,9$, and for $\delta=.1$ and $n=3,4,5$, it gives rank $n$.}

{Let $D$ be Bernoulli on $\{0,1\}$ with $p=\Pr(D=0)$, and write $q=\qstar$. Then $x^{\mathrm{SAA}}=1$ when fewer than $k_n$ sample points equal zero, an event of probability $r_n(p)=\Pr\{\mathrm{Bin}(n,p)<k_n\}$, and $x^{\mathrm{SAA}}=0$ otherwise. The combined-risk regret and the expected post-sample regret on this family are
\[
\begin{aligned}
R_n(p)&=\frac{h\min\{p\,r_n(p),1-\delta\}-b((1-p)r_n(p)-\delta)^+}{1-\delta}-\min\{0,(h+b)(p-q)/(1-\delta)\},\\
\mathrm{PS}_n(p)&=r_n(p)\min\{h,(h+b)(p-q)/(1-\delta)\}-\min\{0,(h+b)(p-q)/(1-\delta)\}.
\end{aligned}
\]
We maximize the two expressions separately over $p\in[0,1]$. The search combines a uniform grid, a finer grid around $q$, local scalar maximization, and the kink points $p\in\{0,1,q,1-\delta\}$ together with the roots of $(1-p)r_n(p)=\delta$ and, when they exist, of $p\,r_n(p)=1-\delta$.}

{At $h=1$, $b=9$, and $\delta=.95$, the rank is $45$ at $n=1{,}000$ and $450$ at $n=10{,}000$. The combined-risk maxima are approximately $3.17180$ and $1.73521$. The sample size $n=1{,}000$ is not a plotted row. In Figures~\ref{fig:saa_performance_by_delta_linear} and~\ref{fig:finite_sample_opt_vs_saa} the grid steps from $n=953$ to $n=1{,}001$, where the ranks are $43$ and $46$. Rank $45$ holds for $n=978,\ldots,1{,}000$. The plotted combined-risk value at $n=1{,}001$ is approximately $3.04074$. Figure~\ref{fig:saa_performance_linear} skips $n=1{,}000$ as well, stepping from $n=892$ to $n=1{,}070$. At $n=10{,}000$, Figure~\ref{fig:saa_performance_linear} uses the separate maximum of expected post-sample regret, approximately $0.07127$. Expected post-sample regret at the combined-risk maximizer is approximately $0.00946$; that is a different evaluation, stored in the last column of the corrected performance-check file. At $\delta=0$ and $n=3$, with the same $h$ and $b$, the rank is $3$ and the worst-case expected regret equals $0.69198046875$.}